\documentclass[11pt,twoside, leqno]{article}

\usepackage{mathrsfs}
\usepackage{amssymb}
\usepackage{amsmath}
\usepackage{amsthm}

\allowdisplaybreaks

\def\rr{{\mathbb R}}
\def\rn{{{\rr}^n}}

\def\rrm{{{\rr}^m}}
\def\rnm{{\rn \times \rrm}}
\def\zz{{\mathbb Z}}
\def\nn{{\mathbb N}}
\def\ch{{\mathcal H}}
\def\cc{{\mathbb C}}
\def\cs{{\mathcal S}}
\def\hr{{\mathcal R}}

\def\lfz{{\lfloor}}
\def\rfz{{\rfloor}}

\def\cg{{\mathcal G}}
\def\cl{{\mathcal L}}

\def\cq{{\mathcal Q}}

\def\cm{{\mathcal M}}

\def\cb{{\mathcal B}}
\def\ca{{\mathcal A}}
\def\ccc{{\mathcal C}}
\def\ck{{\mathcal K}}
\def\ct{{\mathcal T}}
\def\cf{{\mathcal F}}

\def\cbgz{{\cb_\gz}}

\def\fz{\infty}
\def\az{\alpha}
\def\supp{{\mathop\mathrm{\,supp\,}}}

\def\loc{{\mathop\mathrm{\,loc\,}}}
\def\esssup{{\mathop\mathrm{\,esssupp\,}}}
\def\lz{\lambda}
\def\dz{\delta}

\def\ez{\epsilon}
\def\bz{\beta}

\def\gz{{\gamma}}
\def\bgz{{\Gamma}}

\def\Oz{{\Omega}}

\def\vz{\varphi}
\def\tz{\theta}
\def\sz{\sigma}

\def\hs{\hspace{0.3cm}}

\def\wh{\widehat}
\def\wz{\widetilde}
\def\ls{\lesssim}

\def\lp{{L^p_w(\rnm)}}
\def\lq{{L^q_w(\rnm)}}

\def\hpaf{{H^{p,\, q,\,\vec s}_{w,\, {\rm fin}}(\rnm;\,\vec{A})}}
\def\hp{{H^p_w(\rnm;\,\vec{A})}}

\def\hpa{{H^{p, \,q, \,\vec s}_w(\rnm;\,\vec{A})}}

\def\suj{{\sum_{j\in\nn}}}

\def\hat{\widehat}

\def\esup{\mathop\mathrm{\,esssup\,}}

\def\esup{\mathop\mathrm{\,esssup\,}}
\def\dsum{\displaystyle\sum}

\def\dint{\displaystyle\int}

\def\gfz{\genfrac{}{}{0pt}{}}

\def\dfrac{\displaystyle\frac}
\def\dsup{\displaystyle\sup}

\def\r{\right}
\def\lf{\left}
\def\la{\langle}
\def\ra{\rangle}

\newtheorem{thm}{Theorem}[section]
\newtheorem{lem}{Lemma}[section]
\newtheorem{prop}{Proposition}[section]
\newtheorem{cor}{Corollary}[section]
\theoremstyle{definition}
\newtheorem{rem}{Remark}[section]
\newtheorem{defn}{Definition}[section]
\numberwithin{equation}{section}

\begin{document}

\arraycolsep=1pt

\title{{\vspace{-4.5cm}\small\hfill\bf Math. Nachr., to appear}\\
\vspace{4cm}\Large\bf Weighted Anisotropic Product Hardy Spaces and
Boundedness of Sublinear Operators \footnotetext{\hspace{-0.35cm}
{\it 2000 Mathematics Subject Classification}. { Primary 42B30;
Secondary 42B20, 42B25, 42B35.}
\endgraf{\it Key words and phrases.} product space, expansive dilation,
weight, Calder\'on reproducing formula, dyadic
rectangle, atom, grand maximal function, Hardy space, quasi-Banach
space, sublinear operator.
\endgraf
The first author was partially supported by the NSF grant
DMS-0653881 and the third (corresponding) author was supported by the National
Natural Science Foundation (Grant No. 10871025) of China.
\endgraf $^\ast$ Corresponding author.}}
\author{Marcin Bownik, Baode Li,
Dachun Yang$\,^\ast$ and Yuan Zhou}

\date{ }

\maketitle

\begin{center}
\begin{minipage}{13cm}\small
{\noindent{\bf Abstract.}
Let $A_1$ and $A_2$ be expansive
dilations, respectively, on ${\mathbb R}^n$ and ${\mathbb R}^m$.
Let $\vec A\equiv(A_1,\,A_2)$ and $\mathcal A_p( \vec A)$ be the class of product
Muckenhoupt weights on ${\mathbb R}^n\times{\mathbb R}^m$ for $p\in(1,\,\infty]$.
When $p\in(1,\,\infty)$ and
$w\in{\mathcal A}_p(\vec A)$, the authors characterize the weighted
Lebesgue space $L^p_w({\mathbb
R}^n\times{\mathbb R}^m)$ via the anisotropic Lusin-area function
associated with $\vec A$. When $p\in(0,\,1]$, $w\in
{\mathcal A}_\infty(\vec A)$,
the authors introduce the weighted anisotropic product Hardy space
$H^p_w({\mathbb R}^n\times{\mathbb R}^m;\,\vec A)$ via
the anisotropic Lusin-area function and establish its
atomic decomposition.
Moreover, the
authors prove that finite atomic norm on a
dense subspace of $H^p_w({\mathbb R}^n\times{\mathbb R}^m;\vec
A)$ is equivalent with the standard infinite atomic
decomposition norm. As an application, the authors prove that if $T$ is
a sublinear operator and maps all atoms
into uniformly bounded elements of a quasi-Banach space
$\mathcal B $, then $T$ uniquely extends to a bounded sublinear
operator from $H^p_w({\mathbb R}^n\times{\mathbb R}^m;\vec A)$ to
$\mathcal B$. The results of this paper improve the existing
results for weighted product Hardy spaces and are new even
in the unweighted anisotropic setting.}
\end{minipage}
\end{center}

\section{Introduction\label{s1}}

\hskip\parindent The theory of Hardy spaces plays an important role
in various fields of analysis and partial differential equations;
see, for example, \cite{clms, fs72, g2, m94, s94, s93, sw60}. One of
the most important applications of Hardy spaces is that they are
good substitutes of Lebesgue spaces when $p\in(0,\,1]$. For example,
when $p\in(0,\,1]$, it is well-known that Riesz transforms are not
bounded on $L^p(\rn)$, however, they are bounded on Hardy spaces
$H^p(\rn)$. There were several efforts of extending classical
function spaces and related operators
arising in harmonic analysis from Euclidean spaces
to other domains and anisotropic settings; see \cite{b1, ct75, ct77,
fs82-1, st87, tr83, tr92, tr06}. Fabes and Rivi\`ere
\cite{fr1,fr2,r71} initiated the study of singular
integrals with mixed homogeneity, and
Calder\'on and Torchinsky \cite{c77, ct75, ct77}
the study of Hardy spaces associated with anisotropic dilations.
Recently, a theory of anisotropic Hardy
spaces and their weighted theory were developed in \cite{b1, blyz}.
Another direction is the development of the theory of Hardy spaces on product domains initiated by Gundy and Stein
\cite{gs79}. In particular, Chang and Fefferman \cite{cf80, cf82}
characterized the classical product Hardy spaces via atoms.
Fefferman \cite{f88}, Krug \cite{k88} and Zhu \cite{z92} established
the weighted theory of the classical product Hardy spaces, and Sato
\cite{s87, s88} established parabolic Hardy spaces on product
domains. It was also proved that the classical product Hardy spaces
are good substitutes of product Lebesgue spaces when $p\in(0,1]$;
see, for example, \cite{f85, f87, f88, s88, st01}.
Recently, the boundedness of singular integrals on
product Lebesgue spaces was further proved to be useful in solving problems from the several complex variables
by Nagel and Stein \cite{ns06}.

On the other hand, to establish the boundedness of operators on
Hardy spaces, one usually appeals to the atomic decomposition
characterization, see \cite{c77, cf85, c74, f87, fs82-1, mc, tw}, which means
that a function or distribution in Hardy spaces can be represented
as a linear combination of functions of an elementary form, namely,
atoms. Then, the boundedness of operators on Hardy spaces can
be deduced from their behavior on atoms or molecules in principle.
However, caution needs to be taken due to an example  constructed in \cite[Theorem
2]{b2}. There exists a linear functional defined on a dense
subspace of $H^1(\rn)$, which maps all $(1,\,\fz,\,0)$-atoms into
bounded scalars, but yet it does not extend to a bounded linear
functional on the whole $H^1(\rn)$. This implies that the uniform
boundedness of a linear operator $T$
on atoms does not automatically guarantee the
boundedness of $T$ from $H^1(\rn)$ to a  Banach space $\cb$.

Recently, there was a flurry of activity addressing the
problem of boundedness of operators on $H^p(\rr^n)$
via atomic decompositions in addition to older contributions;
see \cite{gr, mc, mtw, tw, y93} and the references therein.
Let $p\in (0,1]$, $q\in[1,\,\fz]\cap(p,\,\fz]$ and $s$ be
an integer no less than $\lfloor n(1/p-1)\rfloor$,
where and in what follows, $\lfloor \cdot \rfloor$
denotes the floor function.
Using the Lusin-area function characterization of classical
Hardy spaces, it was proved in \cite{yz} that if a sublinear
operator $T$ maps all smooth $(p,\,2,\,s)$-atoms into uniformly
bounded elements of a quasi-Banach space $\cb$, then $T$
uniquely extends to a bounded sublinear operator from $H^p(\rn)$ to
$\cb$. This result was generalized to the classical product Hardy
spaces in \cite{cyz}. At the same time, Meda, Sj\"ogren and
Vallarino \cite{msv} independently obtained a related result using
the grand maximal function characterization of $H^p(\rn)$.
Precisely, they proved that the
norm of $H^p(\rn)$ can be reached on some dense subspaces of
$H^p(\rn)$ via finite combinations of $(p,\,q,\,s)$-atoms when
$q<\fz$ and continuous $(p,\,\fz,\,s)$-atoms. Their result
immediately implies that if $T$ is a linear operator and maps all
$(p,\,q,\,s)$-atoms with $q<\fz$ or all continuous
$(p,\,\fz,\,s)$-atoms into uniformly bounded elements of a
Banach space $\cb$, then $T$ uniquely extends to a bounded linear
operator from $H^p(\rn)$ to $\cb$. This result was further
generalized to the weighted anisotropic Hardy spaces in \cite{blyz}
and the Hardy spaces on spaces of homogeneous type enjoying the
reverse doubling property in \cite{gly} when $p\le 1$ and near to
$1$. Very recently, Ricci and Verdera \cite{rv} showed that if $p\in
(0,1)$, then the uniform boundedness of a {\it linear}
operator $T$ on all $(p,\,\fz,\,s)$-atoms
does guarantee the boundedness of $T$ from $H^p(\rn)$ to a Banach space $\cb$.

In this paper, we always let $A_1$ and $A_2$ be expansive
dilations, respectively, on ${\mathbb R}^n$ and ${\mathbb R}^m$. Let
$\vec A\equiv(A_1,\,A_2)$ and $\mathcal A_p( \vec A)$ be the class of product
Muckenhoupt weights on $\rnm$ for $p\in(1,\,\fz]$. When
$p\in(1,\,\fz)$ and $w\in{\mathcal A}_p(\vec A)$, we characterize
the anisotropic weighted Lebesgue space $L^p_w({\mathbb
R}^n\times{\mathbb R}^m)$ via the anisotropic Lusin-area function
associated with expansive dilations. For $p\in(0,\,1]$ and $w\in {\mathcal
A}_\infty(\vec A)$ and admissible triplet $(p,\,q,\,\vec s)_w$
 (see Definition \ref{d4.2} below), we introduce the
weighted anisotropic product Hardy space $H^p_w({\mathbb
R}^n\times{\mathbb R}^m;\vec A)$, the atomic one $H^{p,\, q,\,\vec
s}_w({\mathbb R}^n\times{\mathbb R}^m;\vec A)$ and the finite atomic
one $H^{p,\,q,\,\vec s}_{w,\,{\rm fin}}({\mathbb R}^n\times{\mathbb
R}^m;\vec A)$, respectively, via the anisotropic Lusin-area
function, $(p,\,q,\,\vec s)_w$-atoms and finite linear combinations
of $(p,\,q,\,\vec s)_w^\ast$-atoms. We then prove that
$H^p_w({\mathbb R}^n\times{\mathbb R}^m;\vec A)$ coincides with
$H^{p,\,q,\,\vec s}_w({\mathbb R}^n \times{\mathbb R}^m;\vec A)$,
that $H^{p,\,q,\,\vec s}_{w,\, {\rm fin}}({\mathbb R}^n\times{\mathbb
R}^m;\vec A)$ is dense in $H^p_w({\mathbb R}^n\times{\mathbb
R}^m;\vec A)$ and that both the quasi-norms $\|\cdot\|_{H^p_w({\mathbb
R}^n \times{\mathbb R}^m;\,\vec A)}$ and $\|\cdot\|_{H^{p,\,q,\,\vec
s}_{w,\, {\rm fin}}({\mathbb R}^n \times{\mathbb R}^m;\,\vec A)}$
with $\vec s$ being sufficiently large
are equivalent on $H^{p,\,q,\,\vec s}_{w,\, {\rm fin}}({\mathbb R}^n
\times{\mathbb R}^m;\,\vec A)$. As an application, we prove that if
$T$ is a sublinear operator and maps all $(p,\,q,\,\vec
s)^\ast_w$-atoms into uniformly bounded elements of a
quasi-Banach space $\mathcal B $, then $T$ uniquely extends to a
bounded sublinear operator from $H^p_w({\mathbb R}^n\times{\mathbb
R}^m;\vec A)$ to $\mathcal B$.

We point out that the setting in this paper includes the classical
isotropic product Hardy space theory of Gundy and Stein \cite{gs79}
and Chang and Fefferman \cite{cf80, cf82}, the parabolic
product Hardy space theory of Sato \cite{s87, s88} and the weighted
product Hardy space theory of Fefferman \cite{f88}, Krug
\cite{k88} and Zhu \cite{z92}.  Most results of this paper are new
even in the unweighted setting. They also improve the
corresponding results on the isotropic weighted product Hardy
spaces in \cite{f88}, \cite{k88}, and \cite{z92}. The paper is organized as follows.

In Section \ref{s2}, we recall some notation and definitions
concerning expansive dilations, Muckenhoupt weights and maximal
functions, whose basic properties are also presented. Moreover, we
establish discrete Calder\'on reproducing formulae (see Proposition
\ref{p2.5} below) associated to the product expansive dilations for
distributions vanishing weakly at infinity, which were introduced by
Folland and Stein \cite{fs82-1} on homogeneous groups. These
Calder\'on reproducing formulae are crucial tools for this paper.
Another key tool used in this paper are the dyadic cubes of Christ
\cite{c90}, which substitute the role played by dilated balls and
cubes in \cite{b1,b2,bh}, and are used in deriving the atomic
decomposition of product Hardy spaces via the Lusin-area function.
Here we point out that a subtle relation between the dyadic cubes of
Christ \cite{c90} and dilated balls associated to expansive
dilations is established in Lemma \ref{l2.1}(iv) according to the
levels of dyadic cubes. This relation and the concept of the level
of dyadic cubes play an important role in the whole paper,
especially in the choice of dyadic rectangles of $\rnm$; see
\eqref{e4.1} and \eqref{e5.4} below.

In Section \ref{s3}, for $p\in(1,\,\fz)$ and $w\in\ca_p(\vec A)$
(resp. $w\in\ca_p(A)$), with the aid of the theory of one-parameter
vector-valued Calder\'on-Zygmund operators, we characterize the
anisotropic weighted Lebesgue space $L^p_w(\rnm)$ (resp.
$L^p_w(\rn)$) via the anisotropic Lusin-area function associated
with expansive dilations $\vec A$ (resp. $A$); see Theorem
\ref{t3.1} and Theorem \ref{t3.2} below.

In Section \ref{s4}, let $p\in(0,\,1]$, $w\in\ca_\fz(\vec A)$  and
$(p,\,q,\,\vec s)_w$ be admissible. We introduce the Hardy space
$\hp$ and the atomic one $\hpa$, respectively, via the Lusin-area
function and $(p,\,q,\,\vec s)_w$-atoms. Using some ideas from
\cite{cf80,cf82, f88,z92} and the Calder\'on reproducing formulae
established in Proposition \ref{p2.5}, we prove that $\hp$ coincides
with $\hpa$; see Theorem \ref{t4.1} below. We point out that since
we are working on weighted anisotropic product Hardy spaces, when we
decompose a distribution into a sum of atoms, the dual method for estimating norms of atoms in \cite{cf80} does not work any
more in the current setting. Instead, we invoke a method from Fefferman
\cite{f88} with more subtle estimates involving rescaling techniques specific to the anisotropic setting. We also notice that a variant of
Journ\'e's covering lemma for expansive dilations established in
Lemma \ref{l4.4} is crucial to the proof of the imbedding of $\hpa$
into $\hp$. In fact, Lemma \ref{l4.4} plays an important role in
obtaining the boundedness of operators on $\hp$. In particular,
using Lemma \ref{l4.4}, we obtain the boundedness of the anisotropic
grand maximal function from $\hp$ to $L^p_w(\rnm)$; see Proposition
\ref{p4.1} below.

In Section \ref{s5}, we introduce  $\hpaf$ to be the set of all
finite combinations of $(p,\,q,\,\vec s)_w^\ast$-atoms. Via the
Lusin-area function together with the Calder\'on reproducing formula
and by using ideas from \cite{msv}, we prove that $\hpaf$ is dense
in $\hp$ and that the quasi-norm $\|\cdot\|_{H^p_w({\mathbb
R}^n\times{\mathbb R}^m;\,\vec A)}$ is equivalent to
$\|\cdot\|_{H^{p,\,q,\,\vec s}_{w,\, {\rm fin}}({\mathbb R}^n
\times{\mathbb R}^m;\,\vec A)}$ on $H^{p,\,q,\,\vec s}_{w,\, {\rm
fin}}({\mathbb R}^n\times{\mathbb R}^m;\,\vec A)$ with $\vec s$
being sufficiently large; see Theorem \ref{t5.1} below. In fact, by
a careful choice of dyadic rectangles in $\rnm$ (see \eqref{e5.4}
below), we first construct some finite $(p,\,q,\,\vec
s)^\ast_w$-atoms and then by a subtle size estimate on the
complement of the union of chosen rectangles, we prove that the
difference between the original function and the linear combination
of these finite $(p,\,q,\,\vec s)^\ast_w$-atoms is still a
$(p,\,q,\,\vec s)^\ast_w$-atom multiplied by a small constant. We
should point out that while the main idea comes from \cite{msv},
Meda, Sj\"ogren and Vallarino used the grand maximal function
characterization of the classical Hardy space $H^p(\rn)$ to obtain
the desired estimates instead. See also \cite{blyz} for the weighted
anisotropic Hardy space $H^p_w(\rn; A)$. It is not clear if their
approach \cite{msv} also works here, since so far, it is not known
whether $\hp$ can be characterized via the grand maximal function.
Moreover, comparing with the non-product case (see \cite{blyz, gly,
msv}), our results require additional assumptions \eqref{e5.1} and
\eqref{e5.2} on vanishing moments of atoms.

In Section \ref{s6}, we present applications of Theorem
\ref{t5.1}. If $T$ is a sublinear operator defined on
$H^{p,\,q,\,\vec s}_{w,\, {\rm fin}}({\mathbb R}^n\times{\mathbb R}^m;\vec A)$ and maps all
$(p,\,q,\,\vec s)^\ast_w$-atoms into uniformly bounded elements of a quasi-Banach space $\cb$, then $T$ uniquely extends to a
bounded sublinear operator from $\hp$ to $\cb$; see Theorem
\ref{t6.1} bellow. This result is an extension of \cite[Theorem 1.1]{cyz}. Using Theorem \ref{t6.1} and the Journ\'e's covering lemma, we establish a criteria on the boundedness of
certain sublinear operators via their behavior on rectangular
atoms, which extends and complements a result of Fefferman \cite[Theorem 1]{f87}.

We mention that there exist many predictable applications of
our results in the study of boundedness
of sublinear operators on the weighted product Hardy spaces.
For example, in \cite{lbyz}, we establish
the boundedness on these weighted product Hardy spaces
of singular integrals appearing in the work of Nagel and Stein
\cite{ns06}.

We finally make some conventions. Throughout this paper, we always
use $C$ to denote a positive constant which is independent of the
main parameters involved but whose value may differ from line to
line. Constants with subscripts do not change through the whole
paper. Denote by $\nn$ the set $\{1,2,\cdots\}$ and by $\zz_+$ the
set $\nn\cup\{0\}$. We use $f\ls g$ or $g\gtrsim f$ to denote $f\le
Cg$, and if $f\ls g\ls f$, we then write $f\sim g$. Denote by
$M_n(\rr)$ the set of all real $n\times n$ matrices.

\section{Preliminaries\label{s2}}

\hskip\parindent We begin with the following notation and properties
concerning expansive dilations.

\begin{defn}\label{d2.1}
$A\in M_n(\rr)$ is said to be an {\it expansive dilation}, shortly a
{\it dilation}, if $\min_{\lz\in\sz(A)}|\lz|>1$, where $\sz(A)$ is
the set of all eigenvalues of $A$.
\end{defn}

If $A$ is diagonalizable over $\cc$, we take
$\lz_-\equiv\min\{|\lz|,\ \lz\in\sz(A)\}$ and
$\lz_+\equiv\max\{|\lz|,\ \lz\in\sz(A)\}$. Otherwise, let $\lz_-$
and $\lz_+$ be two {\it positive numbers} such that
$$1<\lz_-<\min\{|\lz|,\, \lz\in\sz(A)\}\le\max\{|\lz|,\,
\lz\in\sz(A)\}<\lz_+.$$

Throughout the whole paper, for a fixed dilation $A$, we always let
$b\equiv|\det A|$.

It was proved in \cite[Lemma 2.2]{b1} that for a given dilation $A$,
there exist an {\it open} and {\it symmetric convex} ellipsoid
$\Delta$ and $r\in (1,\,\fz)$ such
that $\Delta\subset r\Delta\subset A\Delta,$ and one can
additionally assume that $|\Delta|=1$, where $|\Delta|$ denotes the
$n$-dimensional Lebesgue measure of the set $\Delta$. Throughout the
whole paper, we set $B_k\equiv A^k\Delta$ for $k\in \zz$ and let
$\sz$ be the {\it minimum integer} such that $2B_0\subset A^\sz
B_0$. Then $B_k$ is open, $B_k\subset rB_k\subset B_{k+1}$ and
$|B_k|=b^k$. Obviously, $\sz\ge 1$. For any subset $E$ of $\rn$, let
$E^\complement \equiv\rn\setminus E$. Then it is easy to prove (see
\cite[p.\,8]{b1}) that for all $k,\,\ell\in\zz$, we have
\begin{eqnarray}
&&B_k+B_\ell\subset B_{\max(k,\,\ell)+\sz},\label{e2.1}\\
&&B_k+(B_{k+\sz})^\complement\subset (B_k)^\complement,\label{e2.2}
\end{eqnarray}
where $E+F$ denotes the algebraic sums $\{x+y:\, x\in E,\,y\in F\}$
of  sets $E,\, F\subset \rn$ (see \cite[p.\,8]{b1}).

Recall that the homogeneous quasi-norm associated with $A$ was
introduced in \cite[Definition 2.3]{b1} as follows.

\begin{defn}\label{d2.2}
A \textit{homogeneous quasi-norm} associated with an expansive
dilation $A$ is a measurable mapping $\rho:\rn\to[0, \fz)$ satisfy
that

(i) $\rho(x)=0$ if and only if $x=0$;

(ii) $\rho(Ax)=b\rho(x)$ for all $x\in\rn$;

(iii) $\rho(x+y)\le H[\rho(x)+\rho(y)]$ for all $x,\, y\in\rn$,
where $H$ is a constant no less than $1$.
\end{defn}

In the standard dyadic case $A=2I_{n\times n}$, $\rho(x)=|x|^n$ is
an example of homogeneous quasi-norms associated with $A$, where and
in what follows, $I_{n\times n}$ always denotes the $n\times n$
\textit{unit matrix} and $|\cdot|$ is the Euclidean norm in $\rn$.

Define the \textit{step homogeneous quasi-norm} $\rho$ associated
with $A$ and $\Delta$ by setting, for all $x\in\rn$, $\rho(x)=b^k$
if $x\in B_{k+1}\setminus B_k$ or else $0$ if $x=0$. It was proved
that all homogeneous quasi-norms associated with a given dilation
$A$ are equivalent (see \cite[Lemma 2.4]{b1}). Therefore, for a
given expansive dilation $A$, in what follows, for convenience, we
always use the step homogeneous quasi-norm $\rho$.

For the step homogeneous quasi-norm $\rho$, from \eqref{e2.1} and
\eqref{e2.2}, it follows that for all $x,\,y\in\rn$,
$\rho(x+y)\le b^\sz\max\lf\{\rho(x),\,\rho(y)\r\}\le
b^\sz[\rho(x)+\rho(y)];$ see \cite[p.\,8]{b1}.

The following inequalities concerning $A$, $\rho$ and the Euclidean
norm $|\cdot|$ established in \cite[Section 2]{b1} are used in the
whole paper: There exists a positive constant $C$ such that
\begin{eqnarray}\label{e2.3}
C^{-1}[\rho(x)]^{\zeta_-}\le |x|\le C[\rho(x)]^{\zeta_+}&\quad
{\rm for\ all}\ \rho(x)\ge 1,\, {\rm and}\\
\label{e2.4} C^{-1}[\rho(x)]^{\zeta_+}\le |x|\le
C[\rho(x)]^{\zeta_-}&\quad {\rm for\ all}\ \rho(x)\le 1,
\end{eqnarray}
where and in what follows $\zeta_+\equiv \ln (\lz_+)/\ln b$ and
$\zeta_-\equiv \ln (\lz_-)/\ln b$, and that
\begin{eqnarray}
\label{e2.5} C^{-1}b^{j\zeta_-}|x|\le |A^jx|\le
Cb^{j\zeta_+}|x|&\quad
{\rm for\ all \ }\ j\ge 0, \, {\rm and}\\
\label{e2.6} C^{-1}b^{j\zeta_+}|x|\le |A^jx|\le C
b^{j\zeta_-}|x|&\quad {\rm for\ all\ }\ j\le 0.
\end{eqnarray}

Moreover, $(\rn,\, \rho,\, dx)$ is a space of homogeneous type in
the sense of Coifman and Weiss \cite{cw71}, where $dx$ is the
$n$-dimensional Lebesgue measure. On such homogeneous spaces, Christ
\cite{c90} provided an analogue of the grid of Euclidean dyadic
cubes as follows.

\begin{lem}\label{l2.1} Let $A$ be a dilation. There exists a
collection
$\cq\equiv\{Q^k_\az\subset \rn:\, k\in\zz,\ \az\in
I_k\}$ of open subsets, where $I_k$ is certain index set, such that

(i) $|\rn\setminus\cup_\az Q^k_\az|=0$ for each fixed $k$ and
$Q^k_\az\cap Q^k_\bz=\emptyset$ if $\az\not=\bz$;

(ii) for any $\az,\, \bz,\, k,\, \ell$ with $\ell\ge k$, either
$Q^k_\az\cap Q^\ell_\bz=\emptyset$ or $Q^\ell_\az\subset Q^k_\bz$;

(iii) for each $(\ell,\, \bz)$ and each $k<\ell$ there exists a
unique $\az$ such that $Q^\ell_\bz\subset Q^k_\az$;

(iv) there exist certain negative integer $v$ and positive integer
$u$ such that for all $Q^k_\az$ with $k\in\zz$ and $\az\in I_k$,
there exists $x_{Q^k_\az}\in Q^k_\az$ satisfying that for
any $x\in Q^k_\az$,
$x_{Q^k_\az}+B_{vk-u}\subset Q^k_\az\subset x+B_{vk+u}.$
\end{lem}

In what follows, for convenience, we call $k$ the {\it level} of the dyadic
cube $Q^k_\az$ with $k\in\zz$ and $\az\in I_k$ and
denote it by $\ell(Q^k_\az)$.
Lemma \ref{l2.1} can be proved by a slight modification of the proof
of \cite[Theorem 11]{c90}. In fact, we only need to choose $\dz$
in the proof of \cite[Theorem 11]{c90} to be $b^v$ with $v$ being
negative integer. We omit the details. From now on, we call
$\{Q^k_\az\}_{k\in\zz,\, \az\in I_k}$ in Lemma \ref{l2.1} dyadic
cubes.

For any locally integrable function $f$,  the {\it Hardy-Littlewood
maximal function} $\cm(f)$ of $f$ is defined by
$$ \cm(f)(x)\equiv\sup_{k\in\zz}\sup_{x\in y+B_k}\frac{1}{|B_k|}
\int_{y+B_k}|f(z)|\,dz, \qquad x\in\rn.$$ It was proved in
\cite[Theorem 3.6]{b1} that $\cm$ is bounded on $L^p(\rn)$ with
$p\in( 1,\, \fz]$ and bounded from $L^1(\rn)$ to $L^{1,\infty}(\rn)$.

We now recall the weight class of Muckenhoupt associated with $A$
introduced in \cite{bh}.

\begin{defn}\label{d2.3}
Let $p\in[1,\,\fz)$, $A$ a dilation and $w$ a non-negative
measurable function on $\rn$. The function $w$ is said to belong to
the {\it weight class of Muckenhoupt} $\ca_p(A)\equiv\ca_p(\rn;\,
A)$, if there exists a positive constant $C$ such that when $p>1$
$$\sup_{x\in\rn}\sup_{k\in\zz}\lf\{ \frac{1}{|B_k|}\int_{x+B_k}w(y)\,dy\r\}
\lf\{ \frac{1}{|B_k|}\int_{x+B_k}[w(y)]^{-1/(p-1)}\,dy\r\}^{p-1}\le
C,$$ and when $p=1$
$$\sup_{x\in\rn}\sup_{k\in\zz}
\lf\{\frac{1}{|B_k|}\int_{x+B_k} w(y)\,dy\r\}\lf\{\esup_{y\in x+B_k}
[w(y)]^{-1}\r\}\le C;$$ and, the minimal constant $C$ as above is
denoted by $C_{p,\, A,\,n}(w)$.

Define $\ca_\fz(A)\equiv\bigcup_{1\le p<\fz} \ca_p(A)$.
\end{defn}

It is easy to see that if $1\le p\le q\le\fz$, then
$\ca_p(A)\subset\ca_q(A)$.

In what follows, for any non-negative local integrable function $w$
and any Lebesgue measurable set $E$, let $w(E)\equiv\int_E
w(x)\,dx$. For $p \in(0,\,\fz)$, denote by $L^p_w(\rn)$ the
\textit{set of all measurable functions} $f$ such that
$$\|f\|_{L^p_w(\rn)}\equiv\lf\{\int_\rn |f(x)|^pw(x)\,dx\r\}^{1/p}<\fz,$$
and $L^\fz_w(\rn)\equiv L^\fz(\rn)$. The space $L^{1,\infty}_w(\rn)$
denotes the \textit{set of all measurable functions} $f$ such that
$$\|f\|_{L^{1,\infty}_w(\rn)}\equiv\sup_{\lz>0}
\lz w(\{x\in\rn:\, |f(x)|>\lz\})<\fz.$$ Moreover, we have the
following conclusions.

\begin{prop}\label{p2.1}
(i) If $p\in[1,\,\fz)$ and $w\in \ca_p(A)$, then there exists a
positive constant $C$ such that for all $x\in\rn$ and $k,\,m\in\zz$
with $k\le m$,
$$C^{-1}b^{(m-k)/p}\le \frac{w(x+B_m)}{w(x+B_k)}\le Cb^{(m-k)p};$$

(ii) If $p\in(1,\,\fz)$, then the Hardy-Littlewood maximal operator
$\cm$ is bounded on $L^p_w(\rn)$ if and only if $w\in\ca_p(A)$; if
$p=1$, then $\cm$ is bounded from $L^1_w(\rn)$ to $L^{1,\infty}_w(\rn)$
if and only if $w\in\ca_1(A)$.
\end{prop}

Proposition \ref{p2.1}(i) is just
\cite[Proposition 2.1(i)]{blyz}. The proof of Proposition \ref{p2.1}(ii) is also
standard; see \cite{st89,gr,g2} for more details.

Let $\cs(\rn)$ be the {\it space of Schwartz functions on $\rn$} as
in \cite[p,\,11]{b1}, namely, the space of all smooth functions $\vz$
satisfying that for all $\az\in(\zz_+)^n$ and $m\in\zz_+$,
$\|\vz\|_{\az,\,m}\equiv\sup_{x\in\rn}[\rho(x)]^m
|\partial^\az\vz(x)|<\fz,$ where and in what follows,
$\az=(\az_1,\cdots,\az_n)$ and
$\partial^{\az}=(\frac{\partial}{\partial
x_1})^{\az_1}\cdots(\frac{\partial}{\partial x_n})^{\az_n}$. It is
easy to see that $\cs(\rn)$ forms a locally convex complete metric
space endowed with the seminorms
$\{\|\cdot\|_{\az,\,m}\}_{\az\in(\zz_+)^n,\,m\in\zz_+}$. From
\eqref{e2.3} and \eqref{e2.4}, it follows that $\cs(\rn)$ coincides
with the classical space of Schwartz functions; see
\cite[p.\,11]{b1}. Moreover, we denote by $\cs_s(\rn)$ {\it the set of
all $\psi\in\cs(\rn)$ satisfying that $\int_\rn \psi(x)x^\gz\,dx=0$
for all $\gz\in (\zz_+)^n$ with $|\gz|\le s$}.
Let $\cs_\infty(\rn)= \cap_{s\in\nn} \cs_s(\rn)$.

The following lemma is a slight improvement of \cite[Lemma 2.2]{blyz}. We omit the details.
\begin{lem}\label{l2.2}
Let $p\in[1,\,\fz]$ and $w\in \ca_p(A)$. Then

(i) if $1/p+1/p'=1$, then $\cs(\rn)\subset
L^{p'}_{w^{-1/(p-1)}}(\rn)$;

(ii)  $L^p_w(\rn)\subset\cs'(\rn)$ and the inclusion is continuous.
\end{lem}

\begin{lem}\label{l2.3}
$\cm(\chi_{B_{k}})(x)\sim \frac{b^k}{b^k+\rho(x)}$ for all $k\in\zz$
and $x\in\rr^{n}$.
\end{lem}

\begin{proof} Let $\sz$ be as in \eqref{e2.1}. If $x\in B_{k+\sz}$,
then
 $$\cm(\chi_{B_{k}})(x)\ge
\frac{|B_k|}{|B_{k+\sz}|}\gtrsim 1\gtrsim \cm(\chi_{B_k})(x),$$
which together with $\rho(x)\ls b^k$ yields the desired estimate
in this case.

Assume now that $x\not\in B_{k+\sz}$. Then $\rho(x)\gtrsim b^k$.
For any $y+B_\ell$ such that
$x\in y+B_\ell$ and $(y+B_\ell)\cap B_k\not=\emptyset$, assume that
$z_0\in (y+B_\ell)\cap B_k$. By \eqref{e2.1}, we have
$x\in z_0+(y-z_0)+B_\ell\subset B_k+B_\ell+B_\ell\subset
B_{\max(\ell+\sz,\,k)+\sz}.$
From this and $x\not\in B_{k+\sz}$, it follows that $\ell+\sz>k$ and
further $x\in B_{\ell+2\sz}$, which implies that $\rho(x)\ls
b^\ell$. Moreover, by the definition of step homogeneous quasi-norm
$\rho$, there exists $s\in\zz$ such that $x\in B_s\setminus
B_{s-1}$, thus we obtain  $B_s\subset B_{\ell+2\sz}$ and
$\rho(x)=|B_s|$. From this, $\rho(x)\ls b^\ell$ and $B_k\subset
B_s$, it follows that
\begin{eqnarray*}
\cm(\chi_{B_k})(x)=\sup_{\gfz{y\in\rn}{\ell+\sz>k}}\sup_{x\in
y+B_\ell} b^{-\ell}\int_{y+B_\ell}\chi_{B_k}(z)\,dz\ls
\frac{|B_k|}{\rho(x)}\sim \frac{|B_k|}{|B_s|}\ls \cm(\chi_{B_k})(x),
\end{eqnarray*}
which together with $\rho(x)\gtrsim b^k$ gives the desired
estimate. This finishes the proof of the Lemma \ref{l2.3}.
\end{proof}

Let $m,\,n\in\nn$. In what follows, for convenience, we often let $n_1 \equiv n$ and $n_2 \equiv m$.
For $i=1,\,2$, let $A_i\in M_{n_i}(\rr)$ be a dilation and
$b_i$, $B^{(i)}_{k_i}$, $\rho_i$, $u_i$ and $v_i$ associated
with $A_i$ as above.

For any locally integrable function $f$ on $\rnm$, the {\it strong
maximal function} $\cm_s(f)$ is defined by setting, for all $x\in\rnm$,
$$\cm_s(f)(x)\equiv\sup_{k_1,\,k_2\in\zz}\sup_{x\in y+B^{(1)}_{k_1}
\times B^{(2)}_{k_2}} \frac{1}{b_1^{k_1}b_2^{k_2}}
\int_{y+B^{(1)}_{k_1}\times B^{(2)}_{k_2}}|f(z)|\,dz.$$
Obviously,  $\cm_s(f)(x)\le \cm^{(1)}[\cm^{(2)}(f)](x)$
for all $x\in\rnm$ and $\cm_s$ is bounded on
$L^p(\rnm)$ for $p\in (1,\, \fz]$, where $\cm^{(i)}$
denotes the Hardy-Littlewood
maximal operator on $\rr^{n_i}$.

\begin{rem}\label{r2.1}
By a slight modification of the proof of Lemma \ref{l2.3},
we also obtain that for all
$k_1,\,k_2\in\zz$ and $x\in\rr^{n_1}\times \rr^{n_2}$,
$\cm_s(\chi_{B^{(1)}_{k_1}\times B^{(2)}_{k_2}})(x)\sim
\prod_{i=1}^2\frac{b_i^{k_i}}{b_i^{k_i}+\rho_i(x_i)}$.
We omit the details here.
\end{rem}

Now we introduce the weight class of Muckenhoupt on $\rnm$ associated
with $A_1$ and $A_2$, which coincides with the isotropic product weights as in \cite{f87} and \cite{s89}
when $A_1=2I_{n\times n}$ and $A_2=2I_{m\times m}$. Among several equivalent ways of introducing product weights \cite[Theorem VI.6.2]{gr} we adopt the following definition.

\begin{defn}\label{d2.4}
For $i=1,\,2$, let $A_i$ be a dilation on $\rr^{n_i}$ and $\vec
A=(A_1,\,A_2)$. Let $p\in(1,\, \fz)$ and $w$ be a non-negative
measurable function on $\rnm$. The function $w$
 is said to be in the weight class of Muckenhoupt
$\ca_p(\vec A)\equiv\ca_p(\rnm,\,\vec A)$, if $w(x_1,\, \cdot)\in
\ca_p(A_2)$ for almost everywhere
$x_1\in\rn$ and $\esssup_{x_1\in\rn}C_{p,\, A_2,\, m}(w(x_1,\,\cdot))<\fz$,
and $w(\cdot,\, x_2)\in \ca_p(A_1)$ for almost everywhere $x_2\in\rrm$
and $\esssup_{x_2\in\rrm} C_{p,\, A_1,\,n}(w(\cdot,\, x_2))<\fz.$
In what follows, let
$$C_{q,\,\vec A,\,n,\,m}(w)
=\max\lf\{\mathop\esssup_{x_1\in\rn}C_{p,\, A_2,\, m}(w(x_1,\,\cdot)),
\ \mathop\esssup_{x_2\in\rrm} C_{p,\, A_1,\,n}(w(\cdot,\, x_2))\r\}.$$

Define $\ca_\fz(\vec A)\equiv\cup_{1< p<\fz} \ca_p(\vec A)$.
\end{defn}

For any $w\in\ca_\fz(\vec A)$, define the critical index of $w$ by
\begin{eqnarray}\label{e2.7}
q_w\equiv\inf \{q\in(1,\,\fz):\, w\in\ca_q(\vec A)\}.
\end{eqnarray}

Obviously, $q_w\in [1,\,\fz)$. If $q_w\in(1,\,\fz)$, then
$w\not\in\ca_{q_w}$. and if $q_w=1$, Johnson and Neugebauer \cite[p.
254]{jn87} gave an example of $w\not\in \ca_1(2I_{n\times n})$ such
that $q_w=1$. It is easy to see that if $1< p\le q\le \fz$, then
$\ca_p(\vec A)\subset\ca_q(\vec A)$. If $w\in \ca_p(\vec A)$ with
$p\in(1,\,\fz)$, then there exists an $\ez\in(0,\,p-1]$ such that
$w\in\ca_{p-\ez}(\vec A)$ by the reverse H\"older inequality.

Throughout the whole paper, for any measurable set $E\subset \rnm$
and $p\in\rr$, we always set $w^p(E)\equiv \int_E [w(x)]^p\, dx.$
Moreover, by the definition of $\ca_p(\vec A)$ and Proposition
\ref{p2.1}, we have the following Proposition. We omit the details.

\begin{prop}\label{p2.2} Let $\vec A$ be as in Definition \ref{d2.4}.

(i) If $p\in(1, \fz]$ and $w\in\ca_p(\vec A)$, there exists a
positive constant $C$ such that for all $x\in\rnm$ and $k_i,\
\ell_i\in\zz$ with $k_i\le \ell_i$,
\begin{eqnarray*}
C^{-1}b_1^{(\ell_1-k_1)/p}b_2^{(\ell_2-k_2)/p} \le
\frac{w(x+B^{(1)}_{\ell_1}\times B^{(2)}_{\ell_2})}
{w(x+B^{(1)}_{k_1}\times B^{(2)}_{k_2})} \le C
b_1^{(\ell_1-k_1)p}b_2^{(\ell_2-k_2)p};
\end{eqnarray*}

(ii) If $p\in(1,\,\fz)$, $w\in \ca_p(\vec A)$ and $q\in (1,\, \fz]$
then the strong maximal operator $\cm_s$ is bounded on $L^p_w(\rnm)$
and moreover, there exists a positive constant $C$ such that for all
$\{f_j\}_{j\in \nn}\subset L^p_w(\rnm)$,
$$\bigg\| \bigg\{\sum_{j\in\nn}[\cm_s(f_j)]^q\bigg\}^{1/q}\bigg\|_{L^p_w(\rnm)}
\le C\bigg\| \bigg\{\sum_{j\in\nn}|f_j|^q\bigg\}^{1/q}\bigg\|_{L^p_w(\rnm)}.$$
\end{prop}

In fact, the vector-valued inequality (ii) can be obtained simply
by iterating the corresponding vector-valued
inequality for the Hardy-Littlewood maximal function in \cite{aj80}.

For $s_1,\,s_2\in\zz_+$, let $\cs_{s_1,\,s_2}(\rnm)$ be {\it the
set of all functions $\psi\in\cs(\rnm)$ satisfying that
$\int_\rn\psi(x_1,\, x_2)x_1^\gz\, dx_1=0$ for all $\gz\in
(\zz_+)^n$, $|\gz|\le s_1$ and $x_2\in\rrm$, and $\int_\rrm
\psi(x_1,\,x_2)x_2^\bz\, dx_2=0$ for all $\bz\in (\zz_+)^m$,
$|\bz|\le s_2$ and $x_1\in\rn$}.
Let $\cs_{\infty}(\rnm)= \cap_{s_1, s_2 \in\nn} \cs_{s_1,s_2}(\rnm)$.

Throughout the whole paper, for a dilation $A$, we always let
$A^\ast$ be its transpose. For functions $\vz$ on $\rn$, $\psi$ on
$\rnm$ and $k, \ k_1,\ k_2\in\zz$, let
$\vz_k(x)\equiv b^{-k}\vz(A^{-k}x)$ for all $x\in\rn$ and
$\psi_{k_1,\,k_2}(x)\equiv b_1^{-k_1}b_2^{-k_2}\psi(A_1^{-k_1}x_1,\,
A_2^{-k_2}x_2)$ for all $x=(x_1,\,x_2)\in\rnm$.

\begin{prop}\label{p2.3}
(i) Let $\vz\in\cs(\rn)$ and $\int_\rn\vz(x)\, dx=1$. For any
$f\in\cs(\rn)$ (or $f\in \cs'(\rn)$), $f\ast\vz_k\to f$ in $\cs(\rn)$
(or $\cs'(\rn)$) as $k\to -\fz$.

(ii) Let $\vz\in\cs(\rnm)$ and $\int_\rnm \vz(x)\, dx=1$. For any
$f\in\cs(\rnm)$ (or $f\in\cs'(\rnm)$), $f\ast\vz_{k_1,\,k_2}\to f$
in $\cs(\rnm)$ (or $\cs'(\rnm)$) as $k_1,\, k_2\to -\fz$.
\end{prop}

In fact, Proposition \ref{p2.3}(i) is just \cite[Lemma 3.8]{b1}.
The proof of Proposition \ref{p2.3}(ii) is similar to that of (i).
We omit the details.

We recall from \cite{fs82-1} that
 $f\in \cs'(\rn)$ is said to {\it vanish weakly at infinity} if
for any $\vz\in\cs(\rn)$, $f\ast\vz_k\to 0$ in $\cs'(\rn)$ as $k\to
\fz$. Denote by $\cs'_{\fz,\,w}(\rn)$ the collection of all $f\in
\cs'(\rn)$ vanishing weakly at infinity. As pointed out in
\cite{fs82-1}, if $f\in L^p(\rn)$ with $p\in[1,\, \fz)$, then
$f\in\cs'_{\fz,\,w}(\rn)$. Similarly, $f\in \cs'(\rnm)$ is said to
{\it vanish weakly at infinity} if for any $\vz^{(1)}\in\cs(\rrm)$
and $\vz^{(2)}\in\cs(\rn)$, $f\ast\vz_{k_1,\,k_2}\to 0$ in
$\cs'(\rnm)$ as $k_1,\,k_2\to \fz$, where
$\vz(x)\equiv\vz^{(1)}(x_1)\vz^{(2)}(x_2)$ for all
$x=(x_1,\,x_2)\in\rnm$. We also denote by $\cs'_{\fz,\,w}(\rnm)$ {\it the
set of all $f\in \cs'(\rnm)$ vanishing weakly at infinity}.

Now we establish the following Calder\'on reproducing formulae.

\begin{lem} \label{l2.4}
Let $A$ be a dilation on $\rr^{n}$ and $A^\ast$ its
transpose. Let $\vz\in\cs(\rr^{n})$ such that $\supp {\hat\vz}$
is compact and bounded away from the
origin and for all $\xi\in\rr^{n}\setminus\{0\}$,
\begin{equation}\label{e2.8}
\sum_{j\in \zz}\hat{\vz}({( A^\ast)}^{j}\xi)=1.
\end{equation}
Then for any $f\in L^2{(\rr^{n})}$, $f=\sum_{j\in\zz} f\ast\vz_{j}$
in $L^2{(\rr^{n})}$. The same holds in $\cs(\rn)$ or $\cs'(\rn)$,
respectively, for $f\in \cs_\infty(\rr^{n})$ or $f\in
\cs'_{\fz,\,w}(\rr^{n})$.
\end{lem}

\begin{proof}
We first prove the lemma for $f\in L^2(\rn)$. Define
$F(\xi)\equiv\sum_{j\in\zz} |\hat\vz((A^\ast)^{j}\xi)|$ for all
$\xi\in\rn$. Obviously, $F(\xi)=F(A^\ast\xi)$ for all $\xi\in\rn$,
which implies that to show $F\in L^\fz(\rn)$,  it suffices to
consider the values of $F$ on $B^\ast_1\setminus B^\ast_0$, where
$B^\ast_0$ is the unit ball associated with the dilation $A^\ast$.
Let $\rho^\ast$ be the homogeneous quasi-norm associated with
$A^\ast$. Since $\widehat{\vz}\in\cs(\rn)$ and $\widehat{\vz}(0)=0$,
we know that $|\widehat{\vz}(\xi)|\ls\rho^\ast(\xi)^{-1}$ for all
$\xi\in\rn\setminus B_0^\ast$ and $|\widehat{\vz}(\xi)|\ls|\xi|$ for
$\xi\in B_1^\ast$. Thus by \eqref{e2.6}, $b>1$ and $\zeta_->0$, for
any $\xi\in B_1^\ast\setminus B_0^\ast$, we have
\begin{eqnarray}\label{e2.9}
F(\xi)\ls\sum_{j\ge0} \rho^\ast((A^\ast)^{j}\xi)^{-1}+
\sum_{j<0}|(A^\ast)^{j}\xi|\ls \sum_{j\ge0}b^{-j}+\sum_{j<0}
b^{j\zeta_-}\ls1.
\end{eqnarray}
Thus, $F\in L^\fz(\rn)$. By this, the Lebesgue dominated convergence
theorem and \eqref{e2.8}, for $f\in L^2(\rn)$, we have $\hat
f=\sum_{j\in\zz} \hat \vz((A^\ast)^{j}\cdot)\hat f$ in $L^2(\rn)$,
and thus $f=\sum_{j\in\zz}\vz_{j}\ast f $ in  $L^2(\rn)$.

Now let us prove the lemma for $f\in\cs_\fz(\rn)$ (or
$f\in\cs'_{\fz,\,w}(\rn)$). Set $\phi\equiv \sum_{j=0}^\fz \vz_{j}$.
Since $\vz\in\cs(\rn)$ and $\vz_j(x)=b^{-j}\vz(A^{-j}x)$, then
$\phi$ is well-defined pointwise on $\rn$. We claim that
$\phi\in\cs(\rn)$ and $\int_\rn \phi (x)\, dx=1$. Assuming the claim
for the moment, by Proposition \ref{p2.3}, we have
$f\ast\phi_{-N}\to f$ in $\cs(\rn)$ (or $\cs'(\rn)$) as $N\to \fz$.
On the other hand, by H\"older's inequality, for $f\in\cs_\fz(\rn)$ (or
by $f\in\cs'_{\fz,\,w}(\rn)$),
we obtain that $f\ast\phi_N\to 0$ in $\cs(\rn)$ (or $\cs'(\rn)$) as
$N\to \fz$.
Therefore, for $ f\in\cs_\fz(\rn)$ (or $f\in\cs'_{\fz,\,w}(\rn)$),
we have $(f\ast\phi_{-N} -f\ast\phi_N) \to f$ in $\cs(\rn)$ (or
$\cs'(\rn)$) as $N\to \fz$. Moreover, observing that $\phi_k
=\sum_{j=0}^\fz (\vz_j)_k =\sum_{j=k}^\fz \vz_j$, and thus
$\sum_{j=-N}^N \vz_j= \phi_{-N}-\phi_{N+1}, $ we obtain the lemma
for $f\in\cs_\fz(\rn)$ (or $f\in\cs'_{\fz,\,w}(\rn)$).

Let us now prove the above claim. Let
$G(\xi)\equiv \sum_{j=0}^\fz \hat\vz((A^\ast)^j\xi)$ for all
$\xi\in\rn$. Then it suffices to prove that $G\in\cs(\rn)$,
$\phi=\cf^{-1}G$ and $\int_\rn \phi (x)\, dx=1$,
where $\cf^{-1}$ denotes the inverse Fourier
transform.

Since $\supp\hat\vz$ is compact, we may assume that $\supp
\hat\vz\subset B^\ast_{k_0}$ for certain $k_0\in\zz$. Then for any
$j\in\zz_+$, we have $\supp\hat\vz((A^\ast)^j\cdot)\subset
B_{k_0-j}^\ast\subset B^\ast_{k_0}$, which implies that $\supp
G\subset B^\ast_{k_0}$. To prove $G\in\ccc^\fz(\rn)$, for any
$\az\in(\zz_+)^n$ and $\xi\in\rn$, set
$F_\az(\xi)\equiv\sum_{j\in\zz}|\partial^\az
[\hat\vz((A^\ast)^j\xi)]|$. Let us now show $F_\az\in L^\fz(\rn)$.
Notice that for all $\xi\in\rn$,
\begin{eqnarray*}
&F_\az(A^\ast\xi)&=\sum_{j\in\zz}|\partial^\az
[\hat\vz((A^\ast)^{j+1}\xi)]| =\sum_{j\in\zz}|\partial^\az
[\hat\vz((A^\ast)^j\xi)]|=F_\az(\xi),
\end{eqnarray*}
which implies that to verify $F_\az\in L^\fz(\rn)$, we only need to
consider the values of $F_\az$ on $B^\ast_1\setminus B^\ast_0$. By
(2.19) in \cite{bh}, $\vz\in\cs(\rn)$ and $\rho^\ast(\xi)\sim 1$, we
have
\begin{eqnarray*}
|\partial^\az\wh\vz((A^\ast)^{j}\xi)|\ls
b^{j|\az|\zeta_+}|(\partial^\az\hat\vz)((A^\ast)^j\xi)|\ls
b^{j|\az|\zeta_+}\frac1{\rho^\ast((A^\ast)^{j}
\xi)^{(1+|\az|\zeta_+)}}\ls b^{-j}
\end{eqnarray*}
when $j>0$, and
$|\partial^\az\wh\vz((A^\ast)^j\xi)|\ls b^{j|\az|\zeta_-}$ when
$j\le 0$. From this, $b>1$ and $\zeta_->0$, by \eqref{e2.6}, it
follows that
$F_\az(\xi)\ls \sum_{j\le 0} b^{j|\az|\zeta_-}+\sum_{j>0}b^{-j}\ls 1,$
and hence $F_\az\in L^\fz(\rn)$. Notice that $\partial^\az
G(\xi)=\sum_{j=0}^\fz\partial^\az [\hat\vz((A^\ast)^j\xi)]$ for all
$\xi\in\rn$. Thus, $G\in \ccc^\fz(\rn)$. From this and $\supp
G\subset B^\ast_{k_0}$, we deduce $G\in\cs(\rn)$.

Moreover, by the proof of $\supp G\subset B^\ast_{k_0}$,
it is easy to see that $\sum_{j=0}^\fz |\wh\vz((A^\ast)^j\xi)|
\subset B^\ast_{k_0}$, which together with H\"older's
inequality and Minkowski's inequality implies that
\begin{eqnarray*}
&\dint_\rn \sum_{j=0}^\fz |\wh\vz((A^\ast)^j\xi)|\, d\xi&\ls
|B_{k_0}^\ast|^{1/2}\lf\{\int_\rn \lf(\sum_{j=0}^\fz
|\hat\vz((A^\ast)^j\xi)|\r)^2\,d\xi \r\}^{1/2}\\
&&\ls b^{k_0/2}\sum_{j=0}^\fz\lf(\int_\rn
|\hat\vz((A^\ast)^j\xi)|^2\, d\xi \r)^{1/2}\ls
b^{k_0/2}\sum_{j=0}^\fz b^{-j/2}\ls 1.
\end{eqnarray*}
Then by Fubini's theorem, we obtain $\cf^{-1} G
=\sum_{j\in\zz}\cf^{-1}[\hat\vz((A^\ast)^j\cdot)] =\phi$ and hence,
$\phi\in \cs(\rn)$.

Let $e_1=(1,\,0,\,\cdots,\,0)$. Since $\hat\vz\in\cs(\rn)$, by
\eqref{e2.8}, we obtain
\begin{eqnarray*}
\int_\rn \phi(x)\, dx=\hat \phi(0)=\lim_{k\to
-\fz}\hat\phi((A^\ast)^k e_1)=\lim_{k\to-\fz}\sum_{j=0}^\fz
\hat\vz((A^\ast)^{j+k}e_1)=\sum_{j\in\zz}\hat\vz((A^\ast)^j e_1)=1,
\end{eqnarray*}
which completes the proof of our claim and hence the proof of Lemma
\ref{l2.4}.
\end{proof}

\begin{rem}\label{r2.2}
From the proof of Lemma \ref{l2.4}, it is easy to see that if
$\vz\in\cs(\rr^{n})$ and $\hat{\vz}(0)=0$, then
$F(\xi)\equiv\sum_{j\in\zz} |\widehat {\vz}((A^\ast)^{j}\xi)|$
for all $\xi\in\rn$ is bounded on $\rn$.
\end{rem}

Using Lemma \ref{l2.4}, we have the following Calder\'on reproducing
formulae.

\begin{prop}\label{p2.4}
Let $s\in\zz_+$ and $A$ be a dilation on $\rn$. There exist $\tz,\,
\psi\in \cs(\rn)$ such that

(i) $\supp \tz\subset B_0,\, \int_{\rn} x^\gz\tz(x)\,dx=0$
for all $\gz\in(\zz_+)^{n}$ with $|\gz|\le s$,
$\hat\tz(\xi)\ge C>0$ for $\xi$ in certain annulus, where
$C$ is a positive constant;

(ii) $\supp\widehat\psi$ is compact and bounded away
from the origin;

(iii $\sum_{j\in\zz}\widehat\psi((A^\ast)^{j}\xi)\widehat\tz
((A^\ast)^j\xi)=1$ for all $\xi\in\rn\setminus\{0\}$.

\noindent Then for all $f\in L^2{(\rn)}$, $f=\sum_{j\in\zz}
f\ast\psi_{j}\ast\tz_{j}$ in $L^2{(\rn)}$. The same holds in
$\cs(\rn)$ or $\cs'(\rn)$, respectively, for any $f\in
\cs_\infty(\rn)$ or $f\in \cs'_{\fz,\,w}(\rn)$.
\end{prop}

We point out that the existences of such $\tz$ and $\psi$ in
Proposition \ref{p2.4} were proved in
Theorem 5.8 of \cite{bh}. The conclusions of
Proposition \ref{p2.4} then just follow from
Lemma \ref{l2.4} by taking $\vz=\tz\ast\psi$.
Moreover, we also need the following variant on $\rnm$ of
Lemma \ref{l2.4}.

\begin{lem} \label{l2.5}
Let $i=1,\,2$, $A_i$ be a dilation on $\rr^{n_i}$ and
$\vz^{(i)}\in\cs(\rr^{n_i})$ such that $\supp {\hat{\vz^{(i)}}}$ is
compact and bounded away from the origin and for all
$\xi_i\in\rr^{n_i}\setminus\{0\}$, \eqref{e2.8} holds with $A$
replaced by $A_i$, $\vz$ by $\vz^{(i)}$ and $\xi$ by $\xi_i$. Set
$\vz(x)\equiv\vz^{(1)}(x_1)\vz^{(2)}(x_2)$ for all
$x=(x_1,\,x_2)\in\rnm$. Then  for any $f\in L^2{(\rnm)}$,
$f=\sum_{j_1,\, j_2\in\zz} f\ast\vz_{j_1,\, j_2}$ in $L^2{(\rnm)}$.
The same holds in $\cs(\rnm)$ or $\cs'(\rnm)$, respectively, for any
$f\in \cs_\fz(\rnm)$ or $f\in\cs'_{\fz,\,w}(\rnm)$.
\end{lem}

\begin{proof}
We first prove the lemma for $f\in L^2(\rnm)$. For
$\vz=\vz^{(1)}\vz^{(2)}$, by \eqref{e2.9}, we  obtain
that for all $\xi=(\xi_1,\,\xi_2)\in\rnm$,
$$F(\xi)\equiv \sum_{j_1,\,j_2\in\zz}\hat\vz((A_1^\ast)^{j_1}\xi_1,\,
(A_2^\ast)^{j_2}\xi_2)=\sum_{j_1\in\zz}\wh{\vz^{(1)}}((A_1^\ast)^{j_1}
\xi_1)\sum_{j_2\in\zz}\wh{\vz^{(2)}}((A_2^\ast)^{j_2} \xi_2)$$ is
bounded on $\rnm$. Then from this and the fact that
$\sum_{j_1,\,j_2\in\zz} \hat\vz((A_1^\ast)^{j_1}\xi_1,\,
(A_2^\ast)^{j_2}\xi_2)=1$ for any $\xi\in(\rnm)\setminus\{(0,0)\}$,
similarly to Lemma \ref{l2.4}, we deduce the desired formula for
$f\in L^2(\rnm) $.

For $f\in\cs_\fz(\rnm)$ or $f\in\cs'_{\fz,\,w}(\rnm)$, observing
that in the proof of Lemma \ref{l2.4}, we have shown
 that $\phi^{(i)}\equiv \sum_{j_i=0}^\fz
\vz_j\in\cs(\rr^{n_i})$ and $\int_{\rr^{n_i}}\phi(x_i)\, dx_i=1$
for $i=1,\,2$, which imply that
$\phi\equiv \phi^{(1)}\phi^{(2)}\in\cs(\rnm)$ and
$\int_\rnm \phi(x)\, dx=1$. Then, similarly to the proof of
Lemma \ref{l2.4}, we obtain the desired formulae, which completes
the proof of Lemma \ref{l2.5}.
\end{proof}

By Lemma \ref{l2.5}, we have the following proposition.

\begin{prop}\label{p2.5}
Let $s_i\in\zz_+$ and $A_i $ be a dilation on $\rr^{n_i}$ for
$i=1,\,2$. Suppose that $\tz^{(i)},\ \psi^{(i)}\in \cs(\rr^{n_i})$
satisfy conditions (i) through (iii) of Proposition \ref{p2.4} on
$\rr^{n_i}$. Set $\tz(\xi)\equiv\tz^{(1)}(\xi_1)\tz^{(2)}(\xi_2)$
and $\psi(\xi)\equiv\psi^{(1)}(\xi_1)\psi^{(2)}(\xi_2)$ for all
$\xi=(\xi_1,\,\xi_2)\in\rnm$. Then for any $f\in L^2{(\rnm)}$,
$f=\sum_{j_1,\, j_2\in\zz} f\ast\psi_{j_1,\, j_2}\ast\tz_{j_1,\,
j_2}$ in $L^2{(\rnm)}$. The same holds in $\cs(\rnm)$ or
$\cs'(\rnm)$, respectively, for any $f\in \cs_\fz(\rnm)$ or $f\in
\cs'_{\fz,\,w}(\rnm)$.
\end{prop}

\section{Weighted anisotropic Littlewood-Paley theory}\label{s3}

\hskip\parindent We begin with the one parameter Lusin-area
function.

\begin{defn}\label{d3.1}
Let $A$ be a dilation on $\rn$. Suppose $\vz\in\cs(\rn)$ such that
$\hat\vz(0)=0$. For all $f\in\cs'(\rn)$ and $x\in\rn$, define the
{\it anisotropic Lusin-area function} of $f$ by
$$S_\vz(f)(x)\equiv\lf\{\sum_{k\in\zz} b^{-k}\int_{B_k}|f\ast\vz_k
(x-y)|^2\,dy\r\}^{1/2}.$$
\end{defn}

By the Plancherel formula and Remark \ref{r2.2}, we have
\begin{eqnarray}\label{e3.1}
\|S_\vz(f)\|_{L^2(\rn)}^2&&=\sum_{k\in\zz}
b^{-k}\int_{B_k} \int_{\rn}|f\ast\vz_k(x-y)|^2\, dx\,dy\\
&&=\sum_{k\in\zz}\int_\rn |\hat f(\xi)|^2|\hat\vz_k(\xi)|^2\, d\xi
\ls\|\hat f\|^2_{L^2(\rn)}\ls\|f\|^2_{L^2(\rn)},\nonumber
\end{eqnarray}
which implies that $S_\vz$ is bounded on $L^2(\rn)$. Moreover, we
have the following theorem.

\begin{thm}\label{t3.1}
Let $A$ be a dilation on $\rn$, $p\in(1,\, \fz)$, $w\in\ca_p(A)$,
and $\tz,\,\psi$ be as in Proposition \ref{p2.4}. Suppose
$\vz\equiv\tz$ or $\psi$. Then $f\in L^p_w(\rn)$ if and only if
$f\in\cs'_{\fz,\,w}(\rn)$ and $S_\vz(f)\in L^p_w(\rn)$. Moreover,
for all $f\in L^p_w(\rn)$, $\|f\|_{L^p_w(\rn)}\sim
 \|S_\vz(f)\|_{L^p_w(\rn)}$.
\end{thm}

The proof of Theorem \ref{t3.1} will be given later.
Similarly, we can introduce the product Lusin-area function as
follows.

\begin{defn}\label{d3.2}
Let $A_i$ be a dilation on $\rr^{n_i}$ and
$\vz^{(i)}\in\cs(\rr^{n_i})$ with $\widehat{\vz^{(i)}}(0)=0$ for
$i=1,\,2$. Set $\vz(x)\equiv\vz^{(1)}(x_1)\vz^{(2)}(x_2)$ for all
$x=(x_1,\,x_2)\in\rnm$. For all $f\in\cs'(\rnm)$ and $x\in\rnm$,
define the {\it anisotropic product Lusin-area function} of $f$ by
$$ \vec S_\vz(f)(x)\equiv\lf\{\sum_{k_1,\, k_2\in\zz}
b^{-k_1}_1b^{-k_2}_2\int_{B^{(1)}_{k_1}\times
B^{(2)}_{k_2}}|\vz_{k_1,\, k_2}\ast f(x-y)|^2\,dy\r\}^{1/2}.$$
\end{defn}

Then by the Plancherel formula and Remark \ref{r2.2}, similarly to
\eqref{e3.1}, we know that $\vec S_\vz$ is bounded on $L^2(\rnm)$.
Moreover, we have the following product version of Theorem \ref{t3.1} which will be proved later.

\begin{thm}\label{t3.2}
Let $A_i$ be a dilation on $\rr^{n_i}$ for $i=1,\,2$, $p\in(1,\,
\fz)$, $w\in\ca_p(\vec A)$ and $\tz,\,\psi$ be as in Proposition
\ref{p2.5}. Suppose $\vz\equiv\tz$ or $\psi$. Then $f\in
L^p_w(\rnm)$ if and only if $f\in\cs'_{\fz,\,w}(\rnm)$ and $\vec
S_\vz (f)\in L^p_w(\rnm)$. Moreover, for all $f\in L^p_w(\rnm)$,
  $\|f\|_{L^p_w(\rnm)}\sim\|\vec S_\vz(f)\|_{L^p_w(\rnm)}$.
\end{thm}

\begin{rem}\label{r3.1}
For convenience, we can also rewrite $\vec S_\vz(f)$ as
\begin{eqnarray*}
\vec S_\vz(f)(x)=\lf\{\iint_{\bgz(x)}|\vz_{t_1,\,t_2}\ast
f(y)|^2\,dy\,\frac{d\sz(t_1)\,d\sz(t_2)}{b_1^{t_1}b_2^{t_2}}\r\}^{1/2},
\end{eqnarray*}
where $\bgz(x)\equiv \{(y,\ t):\, y\in x+B^{(1)}_{t_1}\times
B^{(2)}_{t_2},\, t=(t_1,\ t_2)\in\rr^2\}$ and $\sz$ is the integer
counting measure on $\rr$, i.\,e., for all $E\subset\rr$, $\sz(E)$
is the number of integers contained in $E$.
\end{rem}

Theorems \ref{t3.1} and \ref{t3.2} will be proved by viewing the Lusin-area
function as the vector-valued Calder\'on-Zygmund operator and
applying a duality argument. In fact, we will verify that the kernel of
Lusin-area function satisfies the standard conditions of
vector-valued Calder\'on-Zygmund operator, and then we will apply a
well-known result on the boundedness of vector-valued
Calder\'on-Zygmund operator in $L^p_w(\rn)$ with $p\in (1,\,\fz)$, see Proposition \ref{p3.1}.

To this end, we first recall the theory of vector-valued
Calder\'on-Zygmund operator. Let $\cb$ be a complex Banach space
with norm $\|\cdot\|_{\cb}$ and $\cb^\ast$ its dual space with norm
$\|\cdot\|_{\cb^\ast}$. A function $f:\rn\to \cb$ is called
$\cb$-measurable, if there exists a measurable subset $\Oz$ of $\rn$
such that $|\rn\setminus \Oz|=0$, the values of $f$ on $\Oz$ are
contained in some separable subspace $\cb_0$ of $\cb$, and for every
$u^\ast\in \cb^\ast$, the complex valued map $x\to \langle u^\ast,\
f(x)\rangle$ is measurable. From this definition and theorem in
\cite[p.131]{y95}, it follows that the function $x\to \|f(x)\|_\cb$
on $\rn$ is measurable. For Banach spaces $\cb_1,\, \cb_2$, define
by $L(\cb_1,\, \cb_2)$ the space of all the bounded linear operators
from $\cb_1$ to $\cb_2$.

For all $p\in (0,\,\fz]$, define by $L^p(\rn,\,\cb)$  the space of
all $\cb$-measurable functions $f$ on $\rn$ satisfying
$$\|f\|_{L^p(\rn,\,\cb)}\equiv \lf\{\int_\rn \|f(x)\|_\cb^p\,
dx\r\}^{1/p}<\fz$$
with a usual modification made when $p=\fz$.
Define by $L^{\fz}_c(\rn,\, \cb)$ the space of $f\in
L^\fz(\rn,\,\cb)$ with compact support.

The proof of the following proposition is presented in Appendix at
the end of the paper.

\begin{prop}\label{p3.1}
Let $A$ be a dilation on $\rn$, and $\cb_1$ and $\cb_2$ be Banach
spaces. Assume that $\ct$ is a linear operator bounded from
$L^2(\rn,\,\cb_1)$ to $L^2(\rn,\,\cb_2)$. Moreover, assume that
there exists a continuous vector-valued function $\ck$:
$\rn\setminus \{0\}\to L(\cb_1,\,\cb_2)$ such that for all $f\in
L_{c}^\fz(\rn,\,\cb_1)$ and $x\not\in\supp f$,
$$\ct(f)(x)=\int_\rn \ck(x-y)f(y)\,dy.$$
If there exist positive constants $C$ and $\ez$ such that for all
$y\in\rn\setminus\{0\}$,
\begin{eqnarray}\label{e3.2}
\|\ck(y)\|_{L(\cb_1,\,\cb_2)}\le \frac C{\rho(y)},
\end{eqnarray}
and for all $x,\,y\in\rn\setminus\{0\}$ with $\rho(x-y)\le b^{-2\sz}
\rho(y)$,
\begin{eqnarray}\label{e3.3}
\|\ck(y)-\ck(x)\|_{L(\cb_1,\,\cb_2)}\le
C\frac{\rho(x-y)^\ez}{\rho(y)^{1+\ez}},
\end{eqnarray}
then for all $p\in (1,\, \fz)$ and $w\in\ca_p(A)$, $\ct$ is bounded
from $L^p_w(\rn,\, \cb_1)$ to $L^p_w(\rn,\, \cb_2)$.
\end{prop}

Now we turn to the proofs of Theorems \ref{t3.1} and
\ref{t3.2}.

\begin{proof}[Proof of Theorem \ref{t3.1}]
Let $f\in L^p_w(\rn)$. By Lemma \ref{l2.2}, $f\in \cs'(\rn)$. To
show that $f$ vanishes weakly at infinity, for any $\vz\in \cs(\rn)$
and $k\in\zz_+$, by H\"older's inequality, we obtain
$|\langle f,\, \vz_k \rangle |\le \|f\|_{L^p_w(\rn)}
\|\vz_k\|_{L^{p'}_{w^{-1/(p-1)}}(\rn)}.$ Moreover, by the
definition of $\ca_p(A)$ and Proposition \ref{p2.1}(i), we have that
for $j\in\zz_+$,
\begin{eqnarray*}
 w^{-1/(p-1)}(B_j)=\int_{B_j} [w(x)]^{-1/(p-1)}\,dx\ls
[w(B_j)]^{-1/(p-1)}|B_j|^{p'}\ls b^{j[p'-1/(p(p-1))]},
\end{eqnarray*}
From this and $\vz\in\cs(\rn)$, it follows that
\begin{eqnarray*}
&&\int_{\rn}|\vz_k(x)|^{p'} [w(x)]^{-1/(p-1)}\, dx\\
&&\quad\ls b^{-kp'}w^{-1/(p-1)}(B_k)+ b^{-kp'} \sum_{j=k}^\fz
\int_{B_{j+1}\setminus
B_j}[b^{-k}\rho(x)]^{-p'} [w(x)]^{-1/(p-1)}\, dx \\
&&\quad\ls \sum_{j=k}^\fz b^{-jp'}w^{-1/(p-1)}(B_j)\ls\sum_{j=k}^\fz
b^{-j/[p(p-1)]}\ls b^{-k/[p(p-1)]},
\end{eqnarray*}
which implies that $f$ vanishes weakly at infinity and hence,
$f\in\cs'_{\fz,\,w}(\rn)$.

We now prove the boundedness of $S_\vz$ on $L_w^p(\rn)$ with
$p\in(1,\,\fz)$. Let
$$
\ch\equiv \{F=\{f_k\}_{k\in\zz}:\, f_k \ {\rm is\ a\
measurable\ function\ on } \ B_k\ {\rm \ for \ any} \ k\in\zz\ {\rm
and}\ \|F\|_{\ch}<\fz\},$$ where $\|F\|_{\ch}\equiv
\{\sum_{k\in\zz}b^{-k}\int_{B_k} |f_k(y)|^2\, dy \}^{1/2}$.
Obviously, $\ch$ is a Hilbert space. For all
$x\in\rn\setminus\{0\}$, set
$\ck(x)\equiv\{\vz_k(x-z):\,
k\in\zz,\,z\in B_k\} \in L(\cc,\,\ch),$
 and for all $f\in L^\fz_c(\rn)$ and
$x\not\in \supp f$, define $\ct:
 L^\fz_c(\rn)\to \ch$ by
$$\ct(f)(x)= \int_\rn \ck(x-y)f(y)\,dy=
\{\vz_k\ast f(x-z): \, z\in B_k,\, k\in\zz\}.$$ Then
$\|\ct(f)(x)\|_\ch=S_\vz(f)(x)$ for all $x\in\rn$. From this and
\eqref{e3.1}, $\ct$ is bounded from $L^2(\rn)$ to $L^2(\rn,\,\ch)$.
To obtain the boundedness of $S_\vz$ on $L^p_w(\rn)$, it suffices to
prove $\ck$  satisfies \eqref{e3.2} and \eqref{e3.3}.

To see \eqref{e3.2}, for $z\in B_k$ and $y\in\rn\setminus \{0\}$,
let $j_0\in\zz$ such that $\rho(y)=b^{j_0}$. By Definition
\ref{d2.2}(iii), $\rho(y)\le b^\sz[\rho(z)+\rho(y-z)]\le
b^\sz[b^k+\rho(y-z)]$, which implies that $b^{j_0-k}\ls
1+b^{-k}\rho(y-z).$ Then for all $y\in\rn\setminus\{0\}$, we obtain
\begin{eqnarray*}
\|\ck(y)\|_{L(\cc,\,\ch)}^2&=&\|
\{\vz_k(y-\cdot)\}_{k\in\zz}\|^2_\ch=\sum_{k\in\zz} b^{-k}\int_{B_k}
|\vz_k(y-z)|^2\, dz\\
&\ls& \sum_{k\in\zz}b^{-k}\int_{B_k}\frac{b^{-2k}}
{[1+b^{-k}\rho(y-z)]^4}\, dz\ls \sum_{k\le j_0}
b^{-2k}b^{-4(j_0-k)}+\sum_{k>j_0}b^{-2k}\\
&\ls& b^{-2j_0}\sim[\rho(y)]^{-2}.
\end{eqnarray*}
which gives \eqref{e3.2}.

To show \eqref{e3.3},
 let $y,\, x\in\rn$ with $y\ne 0$ and
$\rho(x-y)\le b^{-2\sz}\rho(y)$. Without loss of generality, we may
assume that $\rho(x-y)=b^{j_0}$ and $\rho(y)=b^{j_0+j_1+2\sz}$ for
certain $j_0\in\zz$ and $j_1\in\zz_+$. Write
\begin{eqnarray*}
&&\|\ck(x)-\ck(y)\|_{L(\cc,\,\ch)}^2\\
&&\hs=\|\{\vz_k(y-\cdot)-\{\vz_k(x-\cdot)\}_{k\in\zz}\|_\ch^2\\
&&\hs=\sum_{k\in\zz} b^{-3k}\int_{B_k}
|\vz(A^{-k}(y-z))-\vz(A^{-k}(x-z))|^2\, dz\\
&&\hs\le \sum_{k\in\zz}b^{-3k}\int_{B_k}\sup_{\xi\in
B_{j_0}}|\nabla\vz(A^{-k}(y-z-\xi))|^2 |A^{-k}(x-y)|^2\,dz\\
&&\hs\ls
\bigg[ \sum_{k< j_0} + \sum_{j_0\le k\le j_0+j_1} +
\sum_{j_0+j_1<k} \bigg]b^{-3k}\\
&&\hs\hs\times\int_{B_k}\sup_{\xi\in
B_{j_0}}[1+\rho(A^{-k}(y-z-\xi))
]^{-4}|A^{-k}(x-y)|^2\, dz
\equiv {\mathrm I_1}+{\mathrm I_2}+{\mathrm I_3}.
\end{eqnarray*}

To estimate ${\mathrm I_1}$,  since $\rho(A^{-k}(x-y))=b^{j_0-k}>1$
for $k< j_0$, by \eqref{e2.3}, we obtain
\begin{eqnarray}\label{e3.4}
|A^{-k}(x-y)|\ls [\rho(A^{-k}(x-y))]^{\zeta_+}=b^{\zeta_+(j_0-k)}.
\end{eqnarray}
Observing that for $y\in B^\complement_{j_0+j_1+2\sz},\, z\in B_k,\,
j_1\ge 0,\, j_0>k$ and $\xi\in B_{j_0}$, by \eqref{e2.1} and
\eqref{e2.2}, we have
$A^{-k}(y-z-\xi)\in B^\complement_{j_0-k+j_1+2\sz}+B_0+B_{j_0-k}
\subset B^\complement_{j_0-k+j_1+\sz},$
which implies that $\rho(A^{-k}(y-z-\xi))\ge b^{j_0-k+j_1+\sz}$.
From this, \eqref{e3.4}, $\zeta_+<1$, $\rho(x-y)=b^{j_0}$ and
$\rho(y)=b^{j_0+j_1+2\sz}$, it follows that
\begin{eqnarray*}
{\mathrm I_1}\ls \sum_{k<j_0}b^{-2k} b^{-4(j_0-k+
j_1)}b^{2\zeta_+(j_0-k)}\ls b^{-2j_0-4j_1} \ls
\frac{[\rho(y-x)]^2}{[\rho(y)]^4}.
\end{eqnarray*}

To estimate ${\mathrm I_2}$, since  $\rho(A^{-k}(y-x))=b^{j_0-k}\le
1$ for $k\ge j_0$, by \eqref{e2.4}, we obtain
\begin{eqnarray}\label{e3.5}
|A^{-k}(x-y)|\ls [\rho(A^{-k}(x-y))]^{\zeta_-}\sim
b^{\zeta_-(j_0-k)}.
\end{eqnarray}
Moreover, observing that for $j_0\le k \le  j_0+j_1,\,
y\in B^\complement_{j_0+j_1+2\sz},\, z\in B_k,\, \xi\in B_{j_0}$ and
$j_1\ge 0$, by $\eqref{e2.1}$ and $\eqref{e2.2}$, we still have that
$\rho(A^{-k}(y-z-\xi))\ge b^{j_0-k+j_1+\sz}.$ From this,
\eqref{e3.5}, $\rho(x-y)=b^{j_0}$ and $\rho(y)=b^{j_0+j_1+2\sz}$, it
follows that
\begin{eqnarray*}
&{\mathrm I_2}&\ls \sum_{j_0\le k\le j_0+j_1}b^{-2k} b^{-4(j_0-k+
j_1)}b^{2\zeta_-(j_0-k)} \ls b^{-2(j_0+j_1) }b^{-2\zeta_-j_1}\ls
\frac{[\rho(y-x)]^{2\zeta_-}}{[\rho(y)]^{2(1+\zeta_-)}}.
\end{eqnarray*}

To estimate ${\mathrm I_3}$,  by \eqref{e3.5}, $\rho(x-y)=b^{j_0}$,
$\rho(y)=b^{j_0+j_1+2\sz}$ and $j_1\ge 0$,  we have
\begin{eqnarray*}
&{\mathrm I_3}&\ls \sum_{k>j_0+j_1} b^{-2k}b^{2(j_0-k)\zeta_-} \ls
b^{-2(j_0+j)}
b^{-2j_1\zeta_-}\ls\frac{[\rho(y-x)]^{2\zeta_-}}{[\rho(y)]^{2(1+\zeta_-)}}
\end{eqnarray*}
Combining the estimates of ${\mathrm I_1}$, ${\mathrm I_2}$ and
${\mathrm I_3}$, finishes the proof of \eqref{e3.3}. Thus, by
Proposition \ref{p3.1}, we obtain the boundedness of $S_\vz$ on
$L^p_w(\rn)$ for $p\in (1,\,\fz)$.

Conversely, let $f\in\cs'_{\fz,\,w}(\rn)$ and $S_\psi(f)\in
L^p_w(\rn)$ with $p\in (1,\,\fz)$. Set $\widetilde \tz(x)\equiv
\tz(-x)$ for all $x\in\rn$. For any $h\in\cs(\rn)$ with
$\|h\|_{L^{p'}_{w^{-p'/p}}(\rn)}\le 1$, by Proposition \ref{p2.4},
the boundedness of $S_{\wz\tz}$ on $L^p_w(\rn)$ with $p\in
(1,\,\fz)$ and H\"older's inequality, we have
\begin{eqnarray*}
&\lf| \la f,\, h\ra\r|&=\lf| \sum_{k\in\zz}\int_\rn f\ast\psi_k\ast
\tz_k(x) h(x)\, dx\r|=\lf| \sum_{k\in\zz}\int_\rn
f\ast\psi_k(x)h\ast\wz\tz_k(x)  \, dx\r|\\
&&=\lf| \sum_{k\in\zz}b^{-k}\int_\rn \int_{x+B_k}
f\ast\psi_k(x)h\ast\wz\tz_k(x)\, dy  \, dx\r|\\
&&=\lf| \int_\rn\sum_{k\in\zz}b^{-k}\int_{y+B_k}
f\ast\psi_k(x)h\ast\widetilde\tz_k(x)\, dx  \, dy\r|\\
&&\le \int_\rn
\lf\{\sum_{k\in\zz}b^{-k}\int_{y+B_k}|f\ast\psi_k(x)|^2\,
dx\r\}^{1/2} \lf\{\sum_{k\in\zz}b^{-k}\int_{y+B_k}
|h\ast\wz\tz_k(x)|^2\,dx\r\}^{1/2}\, dy\\
&&\le
\|S_\psi(f)\|_{L^p_w(\rn)}\|S_{\wz\tz}(h)\|_{L^{p'}_{w^{-p'/p}}(\rn)}
\ls \|S_{\psi}(f)\|_{L^p_w(\rn)}\|h\|_{L^{p'}_{w^{-p'/p}}(\rn)}\\
&&\ls \|S_{\psi}(f)\|_{L^p_w(\rn)},
\end{eqnarray*}
which together with the density of $\cs(\rn)$ in
$L^{p'}_{w^{-p'/p}}(\rn)$ and
$(L^{p'}_{w^{-p'/p}}(\rn))^\ast=L^p_w(\rn)$ implies that $f\in
L^p_w(\rn)$ and $\|f\|_{L^p_w(\rn)}\ls \|S_\psi(f)\|_{L^p_w}$.
Similarly, for
 $f\in\cs'_{\fz,\,w}(\rn)$ and $S_\tz(f)\in
L^p_w(\rn)$, we have $f\in L^p_w(\rn)$ and $\|f\|_{L^p_w(\rn)}\ls
\|S_\tz(f)\|_{L^p_w}$.
 This finishes the proof of Theorem \ref{t3.1}.
\end{proof}

\begin{proof}[Proof of Theorem \ref{t3.2}]
We shall only prove that $\vec S_\vz$ is bounded on $L^p_w(\rnm)$.
This is because the proofs of the other conclusions are similar to
those of Theorem \ref{t3.1}.

Let $\ch_i$ be the space $\ch$ as in the proof of Theorem \ref{t3.1}
with $B_k$ and $b$ replaced, respectively, by $B^{(i)}_k$ and $b_i$
with $i=1,\,2$. Let $\ch_1 \otimes \ch_2$ be the set of all sequences
$F=\{f_{k_1,\,k_2}\}_{k_1,\,k_2\in\zz}$
such that each $f_{k_1,k_2}$ is measurable on $B^{(1)}_{k_1} \times B^{(2)}_{k_2}$ and
\begin{eqnarray*}
\|F\|_{\ch_1 \otimes \ch_2}
&& \equiv \lf\{\sum_{k_1\in\zz}\sum_{k_2\in\zz}b_1^{-k_1}b_2^{-k_2}
\int_{B^{(2)}_{k_2}}\int_{B^{(1)}_{k_1}}
|f_{k_1,\,k_2}(y_1,\,y_2)|^2\,dy_1\, dy_2\r\}^{1/2} \\
&&
=
\lf\{\sum_{k_2\in\zz}b_2^{-k_2}\int_{B^{(2)}_{k_2}}
\|f_{\cdot,\, k_2}(\cdot,\, y_2)\|_{\ch_1}^2\, dy_2\r\}^{1/2}.
\end{eqnarray*}
The last equation is the consequence of the fact that $\ch_1 \otimes \ch_2$ can be thought of as a collection of measurable $\ch_1$-valued functions $\{f_{\cdot,\, k_2}(\cdot,\, y_2)\}_{k_2\in\zz}$ defined almost everywhere for $y_2\in B^{(2)}_{k_2}$.
Clearly, $\ch_1,\, \ch_2,\, \ch_1 \otimes \ch_2$ are Hilbert spaces.
Here and in what follows, we always let
$$\vz_{k_1}^{(1)}\ast_1g(x_1,\,
x_2)\equiv\int_{\rr^{n_1}}\vz_{k_1}^{(1)}(x_1-y_1)g(y_1,\,x_2)\,
dy_1$$ and  $$\vz_{k_2}^{(2)}\ast_2g(x_1,x_2)
\equiv\int_{\rr^{n_2}}\vz_{k_2}^{(2)}(x_2-y_2)g(x_1,\,y_2)\, dy_2.$$

For any $x_2\in\rr^{n_2}\setminus \{0\}$, define
$\ck^{(2)}(x_2):\, \ch_1\to\ch_1 \otimes \ch_2$ by tensoring
$$\ck^{(2)}(x_2)\equiv\{\vz_{k_2}^{(2)}(x_2-z_2):
\, k_2\in\zz,\,z_2\in B_{k_2}^{(2)}\}.$$
As in the proof of
Theorem \ref{t3.1}, we know that $\ck^{(2)}$ satisfies \eqref{e3.2}
and \eqref{e3.3} with $\cb_1=\ch_1$ and $\cb_2=\ch_1 \otimes \ch_2$.
Moreover, for any 
$F(\cdot)=\{F_{k_1}(y_1,\,\cdot):\, y_1\in
B_{k_1}^{(1)}\}_{k_1\in\zz} \in L^\fz_c(\rr^{m},\,\ch_1)$, define
\begin{eqnarray*}
\ct(F)(x_2)&&\equiv\ck^{(2)}\ast_2F(x_2) \equiv
\lf\{(\vz^{(2)}_{k_2}\ast_2F)(x_2-y_2):\, y_2\in B^{(2)}_{k_2},\,
k_2\in\zz\r\}\\
&&=\lf\{(\vz^{(2)}_{k_2}\ast_2F_{k_1})(y_1,\,x_2-y_2):\, y_1\in
B^{(1)}_{k_1},\, y_2\in B^{(2)}_{k_2},\, k_1,\, k_2\in\zz\r\}.
\end{eqnarray*}
Denote by $\cf_2$ the Fourier transform on the second variable. By
the Plancherel formula and Remark \ref{r2.2}, we have
\begin{eqnarray*}
&&\|\ct(F)\|^2_{L^2(\rrm,\,\ch_1 \otimes \ch_2)}\\
&&\hs=\int_{\rrm}
\sum_{k_1\in\zz}\sum_{k_2\in\zz}b_1^{-k_1}b_2^{-k_2}
\int_{B^{(2)}_{k_2}}\int_{B^{(1)}_{k_1}}
|\vz^{(2)}_{k_2}\ast_2 F_{k_1}(y_1,\,x_2-y_2)|^2\,dy_1\, dy_2\,dx_2\\
&&\hs= \sum_{k_1\in\zz}b_1^{-k_1} \int_{B^{(1)}_{k_1}}
\int_\rrm\sum_{k_2\in\zz} |\widehat{\vz^{(2)}_{k_2}}(\xi_2)|^2
| {\cf_2 F}_{k_1}(y_1,\,\xi_2)|^2d\,\xi_2\,dy_1\\
&&\hs\ls \int_\rrm \sum_{k_1\in\zz}b_1^{-k_1}\int_{B^{(1)}_{k_1}}\sum_{k_2\in\zz}
|(F_{k_1}(y_1,\,y_2))|^2\, dy_1\, dy_2 \ls\|F\|^2_{L^2(\rrm,\,\ch_1)}.
\end{eqnarray*}
Therefore by Proposition \ref{p3.1}, for any $p\in(1,\,\fz)$ and
$w\in\ca_p(A_2)$, $\ct$ is bounded from
$L^p_w(\rrm,\,\ch_1)$ to $L^p_w(\rrm,\,\ch_1 \otimes \ch_2)$.

Let $f\in L^\fz_c(\rnm)$. For any $x_1\in\rr^{n}$ and $x_2\in\rrm$,
set
$$F_{x_1}(x_2)\equiv
\lf\{(\vz^{(1)}_{k_1}\ast_1 f)(x_1-y_1,\, x_2):\, y_1\in
B^{(1)}_{k_1},\, k_1\in\zz\r\} \in \ch_1.$$
Then $F_{x_1} \in L^\fz_c(\rr^{m},\,\ch_1)$ and we have
$$\ct(F_{x_1})(x_2)\equiv
\lf\{(\vz_{k_1,\,k_2}\ast F)(x_1-y_1,\,x_2-y_2):\, y_1\in
B^{(1)}_{k_1},\, k_1\in\zz,\, y_2\in B^{(2)}_{k_2},\, k_2\in\zz\r\},
$$
and $\vec S_\vz(f)(x_1,\,x_2)=\|\ct(F_{x_1})(x_2)\|_{\ch_1 \otimes \ch_2}$.
Recall that by Definition \ref{d2.4}, for almost all $x_1$ (or
$x_2$), $w(x_1,\,\cdot)\in \ca_p(A_2)$ (or $w(\cdot,\,x_2)\in
\ca_p(A_1)$) and the weighted constants are uniformly bounded. Then,
by Theorem \ref{t3.1} for $S_{\vz^{(1)}}$, we have
\begin{eqnarray*}
&\|\vec S_\vz(f)\|_{L^p_w(\rnm)}^p&=\int_\rn\lf\{\int_\rrm
\|\ct(F_{x_1})(x_2)\|^p_{\ch_1 \otimes \ch_2}w(x_1,\,x_2)\, dx_2\r\}\,
dx_1\\
&&\ls \int_\rn\int_\rrm \|F_{x_1}(x_2)\|^p_{\ch_1}w(x_1,\,
x_2)\,dx_2\, dx_1\\
&& \sim\int_\rrm\lf\{\int_\rn [S_{\vz^{(1)}}(f(\cdot,\,x_2))(x_1)]^p
w(x_1,\, x_2)\, dx_1\r\}\, dx_2\ls \|f\|^p_{L^p_w(\rnm)},
\end{eqnarray*}
which completes the proof of Theorem \ref{t3.2}.
\end{proof}

\section{Weighted anisotropic product Hardy spaces}\label{s4}

\hskip\parindent We begin with the notion of weighted anisotropic
product Hardy spaces.

\begin{defn}\label{d4.1}
Let $p\in (0,\ 1]$, $w\in \ca_\fz(\vec A)$ and $q_w$ be as in
\eqref{e2.7}, $\psi$ be as in Proposition \ref{p2.5}. Define the
{\it weighted anisotropic product Hardy space} by
\begin{eqnarray*}
&\hp\equiv & \{f\in \cs'_{\fz,\,w}(\rnm): \\
&&\hs\hs\hs\hs\hs\|f\|_\hp\equiv\|\vec S_\psi(f)\|_\lp<\fz\}.
\end{eqnarray*}
\end{defn}

Notice that if $p\in (q_w,\,\fz)$, where $q_w$ is as in
\eqref{e2.7}, then by Theorem \ref{t3.2}, we obtain
$\hp=L^p_w(\rnm)$ with equivalent norms. If $p\in(1,\, q_w]$, the
element of $\hp$ may be a distribution, and hence, $\hp\ne \lp$; see
\cite[p.\,86]{st89} for one parameter case. For applications
considered in this paper, we concentrate only on $\hp$ with
$p\in(0,\, 1]$.

To define atomic Hardy spaces, we introduce the following notation
and notions. Let $A_i$ be a dilation on $\rr^{n_i}$, and
$\cq^{(i)}$, $\ell(Q_i),\,v_i,\,u_i$ be the same as in Lemma \ref{l2.1}
corresponding to $A_i$ for $i=1,\,2$. Let $\hr\equiv\cq^{(1)}\times
\cq^{(2)}$. For $R\in\hr$, we always write $R=R_1\times R_2$ with
$R_i\in\cq^{(i)}$ and call $R$ dyadic rectangle. For $(k_1,\
k_2)\in\zz\times\zz$, define $\hr_{k_1,k_2}\equiv\{R\in\hr:\,
\ell(R_1)=k_1,\,\ell(R_2)=k_2\}$. For $R\in\hr$, let
\begin{eqnarray}\label{e4.1}
R_+\equiv\{(y,\ t):\, y\in R,\ t=(t_1,\ t_2)\in\rr^2,\ t_i\sim
v_i\ell(R_i)+u_i,\ i=1,\ 2\},
\end{eqnarray}
where and in what follows, $t_i\sim v_i\ell(R_i)+u_i$ always means
\begin{eqnarray}\label{e4.2}
v_i\ell(R_i)+u_i+\sz_i\le t_i<v_i(\ell(R_i)-1)+u_i+\sz_i,
\end{eqnarray}
and $\sz_i$ is as in \eqref{e2.1} and \eqref{e2.2} associated with
$A_i$ for $i=1,\,2$. Note that the inequality \eqref{e4.2} is seemingly reversed since $v_i$'s are negative.

Assume that $\Oz$ is an open set of $\rnm$. A dyadic rectangle
$R\subset\Oz$ is said to be maximal in $\Oz$ if for any rectangle
$S\subset \Oz$ satisfying that $R\subset S$, then $S=R$. Denote by
$m(\Oz)$ the family of all maximal dyadic rectangles contained in
$\Oz$. We choose a positive integer $c_0>2$ such that
$b_1^{-c_0u_1}b_2^{-c_0u_2}\le (b_1^{-2u_1}b_2^{-2u_2}/2)$ and set
\begin{eqnarray}\label{e4.3}
\wz \Oz\equiv\{x\in\rnm,\hs \cm_s(\chi_\Oz)(x)>b_1^{-c_0u_1}
b_2^{-c_0u_2} \}.
\end{eqnarray}
\begin{defn}\label{d4.2}
Let $w\in \ca_\fz(\vec A)$ and $q_w$ be as in \eqref{e2.7}. The
triplet $(p,\ q,\ \vec s)_w$ is said to be {\it admissible} if
$p\in(0,\, 1]$, $q\in [2,\, \fz)\cap(q_w,\, \fz)$ and $s_i\ge\lfloor
(\frac{q_w}p-1)\zeta_{i,-}^{-1}\rfloor$, where $\zeta_{i,-}$ is
defined as in \eqref{e2.3}, $i=1,\ 2$.

A function $a$ is said to be a $(p,\ q,\, \vec s)_w$-{\it atom}
associated to an open set  $\Oz$ of $\rnm$ with $w(\Oz)<\fz$ if

\begin{enumerate}
\item[(I)] $a$ can be written as $a=\sum_{R\in m(\wz\Oz)}a_R$ in
$\cs'(\rnm)$, where $a_R$  satisfies that
\begin{enumerate}
    \item[(i)] $a_R$ is supported on $R''=R''_1\times R''_2$, where
$R''_i\equiv x_{R_i}+B^{(i)}_{v_i(\ell(R_i)-1)+u_i+3\sz_i}$ for
$i=1,\ 2$.
    \item[(ii)] $\int_\rn a_R(x_1,\ x_2)x_1^\az\, dx_1=0$ for all
$|\az|\le s_1$ and  almost all $\ x_2\in \rrm$, and

$\int_\rrm a_R(x_1,\ x_2)x_2^\bz\,dx_2=0$ for all $|\bz|\le s_2$ and
almost all $\ x_1\in \rn$.
\end{enumerate}
Here $a_R$ is called a particle associated with the rectangle $R$.
\item[(II)] $\|a\|_\lq\le [w(\Oz)]^{1/q-1/p}$ and  $\sum_{R\in
m(\wz\Oz)}\|a_R\|_{L^q_w(\rnm)}^q\le [w(\Oz)]^{1-q/p}$.
\end{enumerate}
\end{defn}

\begin{defn}\label{d4.3}
Let $p\in(0,\, 1]$, $w\in \ca_\fz(\vec A)$ and $q_w$ be as in
\eqref{e2.7} and $(p,\ q,\,\vec s)_w$ be an admissible triplet. The
{\it weighted atomic anisotropic product Hardy space} $\hpa$ is
defined to be the collection of all $f\in\cs'(\rnm)$ of the form
$f=\suj \lz_ja_j$ in $\cs'(\rnm)$, where $\suj|\lz_j|^p<\fz$ and
$\{a_j\}_{j\in\nn}$ are $(p,\, q,\, \vec s)_w$-atoms. For $f$ in
$\hpa$, the norm $f$ on $\hpa$ is defined by
$$\|f\|_\hpa\equiv\inf\lf\{\lf(\suj|\lz_i|^p\r)^{1/p}\r \},$$
where the infimum is taken over all the above decompositions of $f$.
\end{defn}

\begin{rem}\label{r4.1}
a) We remark here that the restriction $q\in [2,\fz)$ in Definition
\ref{d4.2} seems reasonable, since we use the Lusin-area function to
introduce $\hp$. Moreover, from the known result on classical
product Hardy spaces, we know that $\{s_i\}_{i=1,\,2}$ in Definition
\ref{d4.2} are best possible.

b) Notice that if $(p,\,q,\,\vec s )_w$ and $(p,\,r,\,\vec t)_w$ are
admissible, $q\le r$ and $s_i\le t_i$ for $i=1,\,2$, then a $(p,\,r
,\,\vec t )_w$-atom is a $(p,\,q,\,\vec s )_w$-atom. Thus, the space
$H^{p,\,r,\,\vec t}_w(\rnm)\subset\hpa$.

\end{rem}

The main result of this section is as follows.

\begin{thm}\label{t4.1}Let $w\in \ca_\fz(\vec A)$ and $q_w$ be as in
\eqref{e2.7}. If $(p,\ q,\,\vec s)_w$ is an admissible triplet, then
$\hp=\hpa$ with equivalent norms.
\end{thm}

From Theorem \ref{t4.1}, we immediately deduce that the definition
of the Hardy space $\hp$ in Definition \ref{d4.1}
is independent of the choice of $\psi$ as
in Proposition \ref{p2.5}.

Since this proof of Theorem \ref{t4.1} is quite complicated, we will use several lemmas. Precisely, by choosing $s_i$ such that
$s_i\ge \lfz(q_w/p-1)\zeta_{i,\,-}\rfz$ and $(s_i+1)\zeta_{i,-}>1$
for $i=1,\,2$, we first prove in Lemma \ref{l4.1} bellow that
$\hp\subset \hpa$. Conversely, for all admissible $(p,\ q,\,\vec
s)_w$, in Lemma \ref{l4.3}, we prove
$$[\hpa\cap\cs'_{\fz,\,w}(\rn)]\subset \hp$$
by using Journ\'e covering lemma established in Lemma \ref{l4.4}
below, and in Lemma \ref{l4.5}, we further show that
$\hpa\subset\cs'_{\fz,\,w}(\rn)$. Combing Lemmas \ref{l4.1},
\ref{l4.3}, \ref{l4.4} and Remark \ref{r4.1} b), then finishes the
proof of Theorem \ref{t4.1}.

\begin{lem}\label{l4.1}
Let $w\in \ca_\fz(\vec A)$ and $q_w$ be as in \eqref{e2.7}.  If
$(p,\ q,\,\vec s)_w$ is an admissible triplet and
$(s_i+1)\zeta_{i,-}>1$ for $i=1,\,2$, then there exists a positive
constant $C$ such that
$\|f\|_\hpa\le C\|f\|_\hp$ for all $f\in\hp$.
\end{lem}

\begin{proof}
To prove this lemma, we borrow some ideas from Fefferman \cite{f86,
f88}. The whole proof is divided into 8 steps. In Step 1, we use the
Calder\'on reproducing formula from Proposition \ref{p2.5} to
decompose $f$ into a sum of functions $\{e_R\}_R$ essentially
supported in rectangles and recombine these functions (according to
the size of the intersection between their corresponding rectangles
and the level sets of the Lusin-area function) to obtain the
particles $\{a_P\}_P$ and atoms $\{a_k\}_k$; see \eqref{e4.6},
\eqref{e4.7} and \eqref{e4.8}. In Step 2 through Step 5, we show
that $\{a_k\}_k$ are $(p,\,q,\,\vec s)_w$-atoms. The crucial step is
to estimate the size of these atoms in Step 3. Here we use the
method from Fefferman \cite{f88} instead of the dual method used in
\cite{cf80} via a subtle inequality \eqref{e4.10}. Step 6 through
Step 8 is devoted to proving the inequality \eqref{e4.10}, which
when $n=m=1$ was established in \cite{cf82, f88}. To obtain
\eqref{e4.10} here, in Step 6, we conclude its proof to the proofs
of the inequalities \eqref{e4.17} and \eqref{e4.18}, which are
given, respectively, in Step 7 and Step 8. To prove \eqref{e4.17}, a
main technique used here is to scale the longer sides of considered
rectangles to $1$ via the anisotropic dilation invariance of the
Lebesgue measure so that we can obtain a desired decreasing factor;
see $|\ell(R_i)-\ell(P_i)|$ in \eqref{e4.17}.

We now start to prove Lemma \ref{l4.1}
by letting $\psi$ be as in Proposition \ref{p2.5} and $f\in\hp$.

{\bf Step 1. Decompose $f$ by the Calder\'on reproducing formula.}

For $k\in\zz$, set
$\Oz_k\equiv\{x\in\rnm:\, \vec S_\psi(f)(x)>2^k\}$ and
$$\hr_k\equiv\{R\in\hr:\, |R\cap\Oz_k|>|R|/2,\ |R\cap\Oz_{k+1}|\le
|R|/2\}.$$
Then for each $R=R_1\times R_2\in\hr$, there exists a
unique $k\in\zz$ such that $R\in\hr_k$. Thus,
\begin{eqnarray}\label{e4.4}
\bigcup_{R\in\hr}R=\bigcup_{k_1,\,k_2\in \,\zz} \bigcup_{R\in
R_{k_1,\,k_2}}R =\bigcup_{k\in\,\zz}\bigcup_{R\in\hr_k}R.
\end{eqnarray}

Moreover, for all $R\in\hr_k$ and all $x\in R$, by Lemma
\ref{l2.1}(iv), we obtain
\begin{eqnarray*}
\cm_s(\chi_{\Oz_k})(x) &&\ge
\frac{1}{b_1^{v_1\ell(R_1)+u_1}b_2^{v_2\ell(R_2)+u_2}}
\int_{x_R+B^{(1)}_{v_1\ell(R_1)+u_1}\times
B^{(2)}_{v_2\ell(R_2)+u_2}}\chi_{\Oz_k}(y)\, dy\\
&& \ge b_1^{-2u_1}b_2^{-2u_2}\frac{|\Oz_k\cap
R|}{|R|}> b_1^{-c_0u_1}b_2^{-c_0u_2},
\end{eqnarray*}
which implies that
\begin{eqnarray}\label{e4.5}
\bigcup_{R\in \hr_k}R\subset \wz \Oz_k,
\end{eqnarray}
where $\wz\Oz_k$ is as in \eqref{e4.3}.

Let $\tz^{(i)}$ and $\psi^{(i)}$ be as in Proposition \ref{p2.4}
such that each $\tz^{(i)}$ has the vanishing moments up to degree
$s_3 \equiv 2\max(s_1,s_2)+1$, where
$s_i\ge\lfloor(q_w/p-1)\zeta_{i,\,-}^{-1}\rfloor $ and
$(s_i+1)\zeta_{i,\,-}>1$, $i=1,\,2$. Set
$\tz\equiv\tz^{(1)}\tz^{(2)}$ and $\psi\equiv\psi^{(1)}\psi^{(2)}$.
Then by Proposition \ref{p2.5}, Lemma \ref{l2.1}(i) and
\eqref{e4.4}, for all $x\in\rnm$, we have
\begin{eqnarray*}
f(x)&&=\sum_{k_1,\, k_2\in\zz}\tz_{k_1,\,k_2}\ast\psi_{k_1,\,k_2}\ast f(x)\\
&&=\sum_{k_1,\, k_2\in\zz}\sum_{\gfz{m_1\sim v_1k_1+u_1}{m_2\sim
v_2k_2+u_2}}\int_\rnm\tz_{m_1,\,m_2}(x-y)\psi_{m_1,\,m_2}\ast f(y)\,dy\nonumber\\
&&=\sum_{k_1,\,
k_2\in\zz}\sum_{R\in\hr_{k_1,\,k_2}}\sum_{\gfz{m_1\sim
v_1k_1+u_1}{m_2\sim v_2k_2+u_2}}\int_R \tz_{m_1,\,m_2}(x-y)
\psi_{m_1,\,m_2}\ast f(y)\,dy\nonumber\\
&&=\sum_{k\in\zz}\sum_{R\in\hr_k}\iint_{R_+} \tz_{t_1,\,t_2}(x-y)
\psi_{t_1,\,t_2}\ast f(y)\,dy\,d\sz(t_1)\,d\sz(t_2)\nonumber
\end{eqnarray*}
in $\cs'(\rnm)$, where $R_+$ is as in \eqref{e4.1} and $\sz$ is the
counting measure on $\rr$.

Set $\lz_k\equiv 2^k[w(\Oz_k)]^{1/p}$ and $ a_k\equiv
\lz_k^{-1}\sum_{R\in\hr_k}e_R$, where for all $x\in\rnm$,
\begin{eqnarray}\label{e4.6}
e_R(x)\equiv \iint_{R_+} \tz_{t_1,\,t_2}(x-y) \psi_{t_1,\,t_2}\ast
f(y)\,dy\,d\sz(t_1)\,d\sz(t_2).
\end{eqnarray}
It is not hard to show that $e_R\in\cs_{s_1,\,s_2}(\rnm)$. Let
$m(\wz\Oz_k)$ be the set of all maximal dyadic rectangles contained
in $\wz\Oz_k$. For each $R\in\hr_k$, by \eqref{e4.5}, there exists
at least one maximal dyadic rectangle in $m(\wz\Oz_k)$ containing
$R$; if there exists only one such maximal dyadic rectangle, we then
denote it by $R^\ast$; if there exist more than one such cubes, we
denote the one which has the ``longest'' side in the $\rn$
``direction'' by $R^\ast$. We point out that $R^\ast$ is unique by
the choice. For each $P\in m(\wz\Oz_k)$, let
\begin{eqnarray}\label{e4.7}
a_P\equiv \lz_k^{-1}\sum_{R\in\hr_k,\,R^\ast= P}e_R,
\end{eqnarray}
and then $a_k=\sum_{P\in m(\wz\Oz_k)}a_P$ in $\cs'(\rnm)$. Moreover,
we rewrite $f$ as
\begin{eqnarray}\label{e4.8}
\hs f=\sum_{k\in\zz}\lz_k a_k=\sum_{k\in\zz}\lz_k\sum_{P\in
m(\wz\Oz_k)} a_P= \sum_{k\in\zz}\lz_k\sum_{P\in
m(\wz\Oz_k)}\sum_{R\in\hr_k,\, R^\ast=P}\lz_k^{-1}e_R
\end{eqnarray}
in $\cs'(\rnm)$.

Then we have
\begin{equation*}
\sum_{k\in\zz}\lz_k^p=\sum_{k\in\zz}2^{pk}w(\Oz_k)\le \|\vec
S_\psi(f)\|^p_\lp=\|f\|_\hp^p.
\end{equation*}
By this and \eqref{e4.8}, to conclude the proof of Lemma \ref{l4.1},
we must show that each $a_k$ is a fixed multiple of a $(p,\ q,\ \vec
s)_w$-atom associated with $\Oz_k$.

{\bf Step 2. Show $\supp a_P\subset P''\equiv P_1''\times P_2''$.}

If $x\in\supp a_P$, by \eqref{e4.7}, $a_P(x)\not=0$ implies that
there exists $R\in\hr_k$ such that $R^\ast=P$ and $e_R(x)\not= 0$.
Recall that for all $t_1,\, t_2\in\zz$ and $(x_1,\, x_2)\in\rnm$,
\begin{eqnarray*}
\tz_{t_1,\,t_2}(x_1,\,x_2)=b_1^{-t_1}b_2^{-t_2}\tz^{(1)}_{t_1}(A_1^{-t_1}x_1)
\tz^{(2)}_{t_2}(A_2^{-t_2}x_2)
\end{eqnarray*}
and $\supp\tz^{(i)}\subset B^{(i)}_0$. If $e_R(x_1,\,x_2)\not=0$, by
\eqref{e4.6}, there exists $(y,\,(t_1,\, t_2))\in R_+$ such that
$A_{i}^{t_i}(x_i-y_i)\in B_0^{(i)}$. Moreover, by \eqref{e4.2}, we
have $t_i< v_i[\ell(R_i)-1]+u_i+\sz_i$. Therefore, by Lemma
\ref{l2.1}(iv) and \eqref{e2.1}, we further have
\begin{eqnarray*}
&x_i\in y_i+B^{(i)}_{t_i}&\subset
x_{R_i}+B^{(i)}_{v_i\ell(R_i)+u_i}+B^{(i)}_{v_i(\ell(R_i)-1)+u_i
+\sz_i}\\
&&\subset x_{R_i}+B^{(i)}_{v_i(\ell(R_i)-1)+u_i+2\sz_i}\equiv R_i'.
\nonumber
\end{eqnarray*}
Thus,
\begin{eqnarray}\label{e4.9}
\supp e_R\subset R'\equiv R_1'\times R_2'.
\end{eqnarray}
Since $ R_i\subset P_i$, by Lemma \ref{l2.1}(iv) and \eqref{e2.1},
we obtain
\begin{eqnarray*}
R_i'&&=x_{R_i}+B^{(i)}_{v_i(\ell(R_i)-1)+u_i+2\sz_i}
\subset
x_{R_i}-x_{P_i}+x_{P_i}+B^{(i)}_{v_i(\ell(P_i)-1)+u_i+2\sz_i}\\
&&\subset
x_{P_i}+B^{(i)}_{v_i\ell(P_i)+u_i}+B^{(i)}_{v_i(\ell(P_i)-1)
+u_i+2\sz_i}\subset x_{P_i}+B^{(i)}_{v_i(\ell(P_i)-1)+u_i+3\sz_i}\equiv
P_i''.
\end{eqnarray*}
From this and \eqref{e4.7}, we obtain $\supp a_P\subset P''\equiv
P_1''\times P_2''$.

{\bf Step 3. Prove $\|a_k\|_{L_w^q(\rnm)}\ls[w(\Oz_k)]^{1/q-1/p}$.}

To this end, we need the following key lemma which will be shown in Steps 6--8 bellow.

\begin{lem}\label{l4.2}
Let $\tz,\,\psi$ be as in Proposition \ref{p2.5}, $\cg$ any set of
dyadic rectangles in $\rnm$, and $e_R$ as in \eqref{e4.6} for any
$R\in\cg$. Then, there exists a positive constant $C$ such that for
all $x\in\rnm$,
\begin{eqnarray}\label{e4.10}
\bigg[\vec S_\tz\bigg(\sum_{R\in\cg}e_R\bigg)(x)\bigg]^2\le C
\sum_{R\in\cg}[\cm_s(c_R\chi_R)(x)]^2,
\end{eqnarray}
where
\begin{eqnarray*}c_R=\lf\{\iint_{R_+}|\psi_{t_1,\,t_2}\ast f(y)|^2\,dy\,\frac{d\sz(t_1)\,d\sz(t_2)}{b_1^{t_1}b_2^{t_2}}\r\}^{1/2}.
\end{eqnarray*}
\end{lem}

Assuming Lemma \ref{l4.2} for the moment, since $q>q_w$, we have
$w\in \ca_q(\vec A)$. By this, Theorem \ref{t3.2}, Lemma \ref{l4.2}
with $\cg=\hr_k$ and Proposition \ref{p2.2}(ii), we have
\begin{align*}
\|a_k\|_\lq&\ls \lz_k^{-1}\bigg\|\vec S_\tz\bigg(\sum_{R\in\hr_k}
e_R\bigg)\bigg\|_\lq\\
&\ls\lz_k^{-1}\bigg\|\bigg\{\sum_{R\in\hr_k}[\cm_s\lf(c_R\chi_R
\r)]^2\bigg\}^{1/2}\bigg\|_{L^q_w(\rnm)}\\
&\ls\lz_k^{-1}\bigg\|\bigg[\sum_{R\in\hr_k}c_R^2\chi_R\bigg]^{1/2}\bigg\|_\lq.
\end{align*}
Since for all $R\in\hr_k$, $|R\cap\Oz_{k+1}|\le |R|/2$ and $R\subset
\wz\Oz_k$ by Lemma \ref{l2.1}(iv) and \eqref{e4.5}, then for all
$x\in R$, we have
\begin{eqnarray*}
&&\cm_s\lf(\chi_{R\cap(\wz\Oz_k\setminus\Oz_{k+1})}\r)(x)\gtrsim
\frac{1}{|R|}\int_R\chi_{R\cap(\wz\Oz_k\setminus\Oz_{k+1})}(y)\,dy
\gtrsim \frac{|R|-|R|/2}{|R|}\gtrsim \chi_R(x).
\end{eqnarray*}
From this and Proposition \ref{p2.2}(ii), it follows that
\begin{align}\label{e4.11}
\|a_k\|_\lq &\ls \lz_k^{-1}\bigg\|
\bigg\{\sum_{R\in\hr_k}\lf[\cm_s\lf(c_R\chi_{R\cap(\wz\Oz_k\setminus
\Oz_{k+1})}\r)\r]^2\bigg\}^{1/2}\bigg\|_\lq\\
&\ls\lz_k^{-1}\bigg\| \bigg(\sum_{R\in\hr_k}c_R^2
\chi_{R\cap(\wz\Oz_k\setminus\Oz_{k+1})} \bigg)^{1/2}\bigg\|_\lq.\nonumber
\end{align}

Moreover, fix $x\in\rnm$. If $R\in\hr_k$ and $x\in R$, then for any
$(y,\, t)\in R_+$, by Lemma \ref{l2.1}(iv) and \eqref{e2.1},
$x_i-y_i\in B^{(i)}_{v_i\ell(R_i)+u_i+\sz_i}\subset B^{(i)}_{t_i},$
which together with Remark \ref{r3.1} and the disjointness of $R_+$
implies that
\begin{align}\label{e4.12}
&\sum_{R\in \hr_k} c_R^2\chi_{R\cap
(\wz\Oz_k\setminus\Oz_{k+1})}(x)\\
&\hs=\sum_{R\in \hr_k} \iint_{R_+}|\psi_{t_1,\,t_2}\ast
f(y)|^2\,dy\,\frac{d\sz(t_1)\,d\sz(t_2)}{b_1^{t_1}b_2^{t_2}}\chi_{R\cap
(\wz\Oz_k\setminus\Oz_{k+1})}(x)\nonumber\\
&\hs\ls [\vec S_\psi(f)(x)]^2\chi_{\wz\Oz_k\setminus\Oz_{k+1}}(x)\ls
2^{2k}\chi_{\wz\Oz_k\setminus\Oz_{k+1}}(x).\nonumber
\end{align}
Notice that $w(\wz\Oz_k)\ls w(\Oz_k)$ by $w\in\ca_q(\vec A)$ and
Proposition \ref{p2.2}(ii). From these estimates, we deduce
\begin{eqnarray}\label{e4.13}
 \|a_k\|_\lq
 \ls 2^{-k}[w(\Oz_k)]^{-1/p}2^k [w(\wz\Oz_k)]^{1/q}
 \ls [w(\Oz_k)]^{1/q-1/p}.
\end{eqnarray}

{\bf Step 4. Prove $\sum_{P\in m(\wz \Oz_k)}\|a_P\|^q_{L^q_w(\rnm)}
\ls[w(\Oz_k)]^{1-q/p}$.} Similarly to the proof of \eqref{e4.13}, by
Theorem \ref{t3.2}, Lemma \ref{l4.2} with $\cg=\{R\in\hr_k:\,
R^\ast=P\}$, the monotonicity of $\ell^{q/2}$ with $q\ge 2$,
$\eqref{e4.12}$ and $w(\wz\Oz_k)\ls w(\Oz_k)$, we have
\begin{eqnarray}\label{e4.14}
&&\sum_{P\in m(\wz \Oz_k)}\|a_P\|^q_{L^q_w(\rnm)}\\
&&\hs\ls \lz^{-q}_k\sum_{P\in m(\wz \Oz_k)}\bigg\| \vec
S_\tz\bigg(\sum_{R\in\hr_k, \,R^\ast=P}e_R\bigg)\bigg\|^q_{L^q_w(\rnm)}\nonumber\\
&&\hs \ls \lz_k^{-q}\sum_{P\in m(\wz\Oz_k)}\bigg\|\bigg\{
\sum_{R\in\hr_k,
\,R^\ast=P}c_R^2\chi_{R\cap(\wz\Oz_k\setminus\Oz_{k+1})}\bigg\}^{1/2}
\bigg\|^q_{L^q_w(\rnm)}\nonumber\\
&&\hs\ls \lz_k^{-q}\lf\|\vec
S_\psi(f)\chi_{\wz\Oz_k\setminus\Oz_{k+1}}
\r\|^q_{L^q_w(\rnm)}\nonumber\\
&&\hs\ls 2^{-qk}[w(\Oz_k)]^{-q/p}w(\wz\Oz_k) 2^{q(k+1)}\ls
[w(\Oz_k)]^{1-q/p}.\nonumber
\end{eqnarray}

{\bf Step 5. Show the vanishing moments of $a_P$.}

By \eqref{e4.12} and $w(\wz\Oz_k)\ls w(\Oz_k)$, we have
\begin{align}\label{e4.15}
&\lz_k^{-q}\int_\rnm\bigg\{\iint_{\bigcup_{R\in\hr_k}R_+}
|f\ast\psi_{t_1,\,t_2}(y)|^2\\
&\hs\quad\quad\times\chi_{R\cap(\wz\Oz_k\setminus\Oz_k+1)}(x)\,
dy
\frac{d\sz(t_1)\,d\sz(t_2)}
{b_1^{t_1}b_2^{t_2}}\bigg\}^{q/2}w(x)\,dx \nonumber\\
&\hs=\lz_k^{-q}\bigg\| \bigg(\sum_{R\in\hr_k}c_R^2
\chi_{R\cap(\wz\Oz_k\setminus\Oz_{k+1})}
\bigg)^{1/2}\bigg\|^q_{L^q_w(\rnm)}\ls
[w(\Oz_k)]^{1-q/p}<\fz.\nonumber
\end{align}
Take any $N\in\nn$ and let
$\hr_{k,\,N}\equiv\{R\in\hr_k:\, |\ell(R_i)|>N,\, i=1,\,2\}.$
Replacing $a_k$ by
$\lz_k^{-1}\sum_{R\in\hr_{k,\,N}}e_R$, similarly to the estimate of
\eqref{e4.11}, we obtain
\begin{align*}
&\bigg\|\lz_k^{-1}
\dsum_{R\in\hr_{k,\,N}}e_R\bigg\|^q_{L^q_w(\rnm)}\\
&\ls\lz_k^{-q}\bigg\|
\bigg(\sum_{R\in\hr_{k,\,N}}c_R^2
\chi_{R\cap(\wz\Oz_k\setminus\Oz_{k+1})}
\bigg)^{1/2}\bigg\|^q_{L^q_w(\rnm)}\\
&\sim
\lz_k^{-q}\int_\rnm\bigg\{\iint_{\bigcup_{R\in\hr_{k,\,N}}R_+}
|f\ast\psi_{t_1,\,t_2}(y)|^2\\
&\hs\hs\times\chi_{R\cap(\wz\Oz_k\setminus\Oz_k+1)}(x)
dy\frac{d\sz(t_1)\,d\sz(t_2)}
{b_1^{t_1}b_2^{t_2}}\bigg\}^{q/2} w(x)\,dx.
\end{align*}
Then by \eqref{e4.15} and Lebesgue dominated convergence theorem, we
have
$$\bigg\|\lz_k^{-1}\sum_{R\in\hr_{k,\,N}}e_R\bigg\|_{L^q_w(\rnm)}\to
0,$$ as $N\to \fz$, which implies that
$a_P=\lz_k^{-1}\sum_{R\in\hr_k,\,R^\ast=P}e_R$ converges in
$L^q_w(\rnm)$, and thus for almost everywhere $x_2\in\rrm$,
$a_P(\cdot,\, x_2)\in L^q_{w(\cdot,\, x_2)}(\rn)$. Moreover, recall
that $\tz$ has vanishing moments
$s_1\ge\lfloor(q_w/p-1)\zeta_{1,\,-}^{-1}\rfloor$ in the first
variable and so is $e_R$. Let $h_1(x_1)\equiv
x_1^\az\chi_{P_1''}(x_1)$ with $|\az|\le s_1$ and $\wz q\in\rr_+$
such that $q^{-1}+(\wz q)^{-1}=1$. Obviously, $h_1\in L^{\wz
q}_{w^{-\wz q/q}(\cdot,\, x_2)}(\rn)$. Then by the fact that $\supp
a_P(\cdot,\, x_2)\subset P_1''$, $(L^{\wz q}_{w^{-\wz q/q}(\cdot,\,
x_2)}(\rn))^\ast=L^q_{w(\cdot,\, x_2)}(\rn)$ and $\supp e_R(\cdot,\,
x_2)\subset P_1''$, we have
\begin{eqnarray*}
\dint_\rn a_P(x_1,\,x_2)x_1^\az\, dx_1&&=\langle a_P(\cdot,\,
x_2),\, h_1\rangle=\sum_{R^\ast=P,\,
R\in\hr_k}\langle e_R(\cdot,\, x_2),\, h_1\rangle\\
&&=\sum_{R^\ast=P,\, R\in\hr_k}\int_\rn e_R(x_1,\, x_2)x_1^{\az_1}\,
dx_1=0.
\end{eqnarray*}
Thus, $a_P$ has vanishing moments up to order $s_1$ in the first
variable. By symmetry, $a_P$  has vanishing moments up
to order $s_2$ in the second variable.

Combining Steps 3 through 5 shows that $a_k$ is a fixed multiple
of a $(p,\,q,\,\vec s)_w$-atom associated with $\Oz_k$. To finish
the proof of Lemma \ref{l4.1}, we still need to show Lemma
\ref{l4.2}.

{\bf Step 6. Proof of Lemma \ref{l4.2}.}

For $P\in\hr$, let $P_+$ be as in \eqref{e4.1}. For all $x\in\rnm$,
by Remark \ref{r3.1}, we have
\begin{align}\label{e4.16}
&\bigg[\vec
S_\tz\bigg(\sum_{R\in\cg}e_R\bigg)(x)\bigg]^2\\
&\hs=\iint_{\bgz(x)}\bigg|\tz_{t_1,\,t_2}
\ast\bigg(\sum_{R\in\cg}e_R\bigg)(y)\bigg|^2\,
dy\,\frac{d\sz(t_1)\,d\sz(t_2)}{b_1^{t_1}b_2^{t_2}}\nonumber\\
&\hs\le \sum_{P\in\hr,\ P_+\cap\bgz(x)\not=\emptyset}\iint_{P_+}\bigg[
\sum_{R\in\cg}|e_R \ast \tz_{t_1,\,t_2}(y)|\bigg]^2\,dy \,
\frac{d\sz(t_1)\,d\sz(t_2)}{b_1^{t_1}b_2^{t_2}}.\nonumber
\end{align}

For any $(y,\, t)\in P_+$ with $P_+\cap \Gamma(x)\ne \emptyset$, we
will prove in Step 7 that if $P'\cap R'=\emptyset$, $e_R\ast
\vz_{t_1,\,t_2}(y)\equiv 0$, or else,
\begin{eqnarray}\label{e4.17}
|e_R\ast\tz_{t_1,\,t_2}(y)|\ls c_R\cm_s(\chi_R)(x)
\prod_{i=1}^2b_i^{(s_i+1)v_i|\ell(R_i)-\ell(P_i)|\zeta_{i,-}}.
\end{eqnarray}
For any $P\in\hr$, we will show in Step 8 that
\begin{eqnarray}\label{e4.18}
\sum_{R\in\mathcal R,\ R'\cap P'\ne\emptyset} \prod_{i=1}^2
b_i^{(s_i+1)v_i|\ell(R_i)-\ell(P_i)|\zeta_{i,-}}\ls 1.
\end{eqnarray}
Assuming that \eqref{e4.17} and \eqref{e4.18} for the moment, for
any $(y,\, t)\in P_+$ and $P_+\cap \Gamma(x)\ne \emptyset$, by
\eqref{e4.17}, the Cauchy-Schwarz inequality and \eqref{e4.18}, we
obtain
\begin{align*}
\bigg(\sum_{R\in\cg}|e_R\ast\tz_{t_1,\,t_2}(y)|\bigg)^2
&\ls \bigg\{
\sum_{R\in\cg,\, R'\cap
P'\not=\emptyset}c_R\cm_s(\chi_R)(x)\prod_{i=1}^2
b_i^{(s_i+1)v_i|\ell(R_i)-\ell(P_i)|\zeta_{i,-}}\bigg\}^2
\\
& \ls \sum_{R\in\cg,\,R'\cap P'\not=\emptyset}
c_R^2[\cm_s(\chi_R)(x)]^2\prod_{i=1}^2
b_i^{(s_i+1)v_i|\ell(R_i)-\ell(P_i)|\zeta_{i,-}}.
\end{align*}
From this, \eqref{e4.16} and \eqref{e4.18}, it follows that
\begin{align*}
&\bigg[\vec S_\tz\bigg(\sum_{R\in\cg}e_R\bigg)(x)\bigg]^2\\
&\hs\ls \sum_{P_+\cap \bgz(x)\not=\emptyset}\iint_{P_+}
\sum_{\gfz{R\in\cg}{R'\cap
P'\not=\emptyset}}c_R^2[\cm_s(\chi_R)(x)]^2 \prod_{i=1}^2
b_i^{(s_i+1)v_i|\ell(R_i)-\ell(P_i)|\zeta_{i,-}}\,dy \,
\frac{d\sz(t_1)\,d\sz(t_2)}{b^{t_1}_1b^{t_2}_2}\\
&\hs\ls\sum_{R\in\cg} c^2_R[\cm_s(\chi_R)(x)]^2
\bigg\{\sum_{\gfz{P\in\mathcal R}{R'\cap P'\not=\emptyset} } \prod_{i=1}^2
b_i^{(s_i+1)v_i|\ell(R_i)-\ell(P_i)|\zeta_{i,-}}\bigg\}
\ls \sum_{R\in\cg}c_R^2[\cm_s(\chi_R)(x)]^2,
\end{align*}
which yields \eqref{e4.10}. To finish the proof of Lemma \ref{l4.2},
we still need to show \eqref{e4.17} and \eqref{e4.18}.

{\bf Step 7. Show \eqref{e4.17}.}

Consider first the trivial case when $R'\cap P'=\emptyset$. In this
case we claim that for $(y,\,(t_1,\,t_2))\in P_+$,  we have $e_R\ast
\tz_{t_1,\,t_2}(y)= 0$. By \eqref{e4.9}, we have
$$e_R\ast\tz_{t_1,\,t_2}(y_1,\,y_2)=\iint_{R'}e_R(z_1,\,z_2)
\tz_{t_1,\,t_2}(y_1-z_1,\,y_2-z_2)\,dz_1\,dz_2.$$  Recall that
$$\tz_{t_1,\,t_2}(y_1-z_1,\,y_2-z_2)=b_1^{-t_1}b_2^{-t_2}\tz^{(1)}
(A_1^{-t_1}(y_1-z_1))\tz^{(2)}(A_2^{-t_2}(y_2-z_2)),$$ and
$\supp\tz^{(i)}\subset B^{(i)}_0$ for $i=1,\,2$. Moreover, since
$(y,\,(t_1,\,t_2))\in P_+$, by \eqref{e4.1}, \eqref{e4.2} and Lemma
\ref{l2.1}(iv), we obtain $y_i\in P_i\subset
x_{P_i}+B^{(i)}_{v_i\ell(P_i)+u_i}$ and
$t_i<v_i(\ell(P_i)-1)+u_i+\sz_i$ for $i=1,\,2$. Therefore, if $(y,\,
(t_1,\,t_2))\in P_+$ and $\tz_{t_1,\,t_2}(y_1-z_1,\,y_2-z_2)\not=0$,
then by \eqref{e2.1}, we have
\[
z_i\in y_i+ B^{(i)}_{t_i}\subset
x_{P_i}+B^{(i)}_{v_i\ell(P_i)+u_i}+B^{(i)}_{v_i(\ell(P_i)-1)+u_i+\sz_i}\subset x_{P_i}+B^{(i)}_{v_i(\ell(P_i)-1)+u_i+2\sz_i}= P_i'.
\]
Thus, for all $(y,\,(t_1,\,t_2))\in P_+$, we have
\begin{eqnarray}\label{e4.19}
e_R\ast\tz_{t_1,\,t_2}(y)=\int_{R'\cap P'}e_R(z)
\tz_{t_1,\,t_2}(y-z)\,dz,
\end{eqnarray}
and if $P'\cap R'=\emptyset$, we obtain
$e_R\ast\tz_{t_1,\,t_2}(y)=0$.

We now consider the non-trivial case $R'\cap P'\not=\emptyset$. We
shall establish $\eqref{e4.17}$ by considering the following four
subcases.

\vspace{0.2cm}

{\it Case I.} $\ell(R_1)\ge\ell(P_1)$ and $\ell(R_2)\ge \ell(P_2)$.
Let $i=1,\,2$. We first observe that for any $(y,\,(t_1,\,t_2))\in P_+$
 and $z_i\in R_i' \equiv
x_{R_i}+B^{(i)}_{v_i[\ell(R_i)-1]+u_i+2\sz_i}$, by
$t_i\ge v_i\ell(P_i)+u_i+\sz_i$, we have
\[
z_i'\equiv A_i^{-t_i}z_i\in A_i^{-t_i}x_{R_i}+
B^{(i)}_{v_i[\ell(R_i)-1]+u_i+2\sz_i-t_i}
\subset A_i^{-t_i}x_{R_i}+B^{(i)}_{v_i [\ell(R_i)-1-\ell(P_i)]
+\sz_i}\equiv \bar R_i.
\]
Let $\bar R\equiv \bar R_1\times \bar R_2$.
Then for any $z_i'\in \bar R_i$, since $-v_i,\, \sz_i>0$ and
$v_i[\ell(R_i)-\ell(P_i)]\le0$, by \eqref{e2.5} and
\eqref{e2.4}, we obtain
\begin{eqnarray}\label{e4.20}
|z_i'-A_i^{-t_i}x_{R_i}|=|A_i^{-v_i+\sz_i}[A_i^{v_i-\sz_i}
(z_i'-A_i^{-t_i}x_{P_i})]|\ls
b_i^{v_i[\ell(R_i)-\ell(P_i)]\zeta_{i,\,-}}.
\end{eqnarray}
On the other hand, by the Cauchy-Schwarz inequality,
$\tz\in\cs(\rnm)$ and Lemma \ref{l2.1}(iv), we have
\begin{eqnarray}\label{e4.21}
&|e_R(x)|^2&\le
c_R^2\iint_{R_+}|\tz_{t_1,\,t_2}(x_1-y_1,\,x_2-y_2)|^2
b_1^{t_1}b_2^{t_2}\, dy_1\,dy_2\, d\sz(t_1)\,d\sz(t_2) \\
&&\ls c_R^2\sum_{t_1\sim v_1\ell(R_1)+u_1} \sum_{t_2\sim
v_2\ell(R_2)+u_2} |R|b_1^{-t_1}b_2^{-t_2} \ls c_R^2.\nonumber
\end{eqnarray}
Let
\[
\mathcal P^{(i)}_{w_i}(z_i) = \sum_{|\az_i|\le s_i}\frac 1{\az_i!}
\partial^{\az_i} \tz^{(i)}
(w_i)(z_i-w_i
)^{\az_i}
\]
be the Taylor polynomial of $\tz^{(i)}$ about $w_i \in \rr^{n_i}$ of
degree $s_i$. For any $(y,\,(t_1,\,t_2))\in P_+$, since
$e_R\in\cs_{s_1,\,s_2}(\rnm)$, $\tz=\tz^{(1)}\tz^{(2)}$ and
$\tz^{(i)}\in\cs_{s_i}(\rr^{n_i})$ for $i=1,\,2$, by \eqref{e4.9},
Taylor's remainder theorem, \eqref{e4.20} and \eqref{e4.21}, we
obtain
\begin{align}\label{e4.22}
&|e_R\ast \tz_{t_1,\,t_2}(y)|\\
&\hs=\bigg|\int_{\bar R} e_R(A_1^{t_1}z_1,\, A_2^{t_2}z_2)
\prod_{i=1}^2\tz^{(i)}(A_i^{-t_i}y_i-z_i)\, dz\bigg|\nonumber\\
&\hs=\bigg|\int_{\bar R} e_R(A_1^{t_1}z_1,\, A_2^{t_2}z_2)
\prod_{i=1}^2 (\tz^{(i)}(A_i^{-t_i}y_i-z_i) -
\mathcal P^{(i)}_{A_i^{-t_i}y_i-A_i^{-t_i}x_{R_i}}(A_i^{-t_i}y_i-z_i))dz \bigg|\nonumber
\\
&\hs\le \int_{\bar R_1 \times \bar R_2} |e_R(A_1^{t_1}z_1,\, A_2^{t_2}z_2)|
\prod_{i=1}^2 |A_i^{-t_i}x_{R_i}-z_i|^{s_i+1} \, dz_1 dz_2
\nonumber\\
&\hs\ls  c_R \prod_{i=1}^2
b_i^{v_i[\ell(R_i)-\ell(P_i)]}
b_i^{(s_i+1)v_i[\ell(R_i)-\ell(P_i)]\zeta_{i,\,-}}.\nonumber
\end{align}

Observing that since $\ell(P_i)\le \ell(R_i)$ and $P_i'\cap
R_i'\not=\emptyset$ for $i=1,\, 2$, by \eqref{e2.1} and Lemma
\ref{l2.1}(iv), it is easy to see
\begin{eqnarray}\label{e4.23}
R_i'\subset P_i'''\equiv
x_{P_i}+B^{(i)}_{v_i(\ell(P_i)-1)+u_i+4\sz_i},
\end{eqnarray}
and hence $R'\subset P'''$. Moreover, for any $x\in\rnm$ and $
\bgz(x)\cap P_+\not=\emptyset$, by \eqref{e2.1} and Lemma
\ref{l2.1}(iv), we obtain
\begin{eqnarray}\label{e4.24}
x\in P'.
\end{eqnarray}

By Lemma \ref{l2.1}(iv), we have that $b_i^{v_i\ell(R_i)}\sim
|R'_i|$ and $b_i^{v_i\ell(P_i)}\sim |P_i'''|$. By this,
\eqref{e4.23}, \eqref{e4.24}, Lemma \ref{l2.1}(iv) and Remark
\ref{r2.1}, we have that for any $x\in\rnm$ and $\bgz(x)\cap
P_+\not=\emptyset$,
\begin{equation}\label{e4.25}
\prod_{i=1}^2
b_i^{v_i[\ell(R_i)-\ell(P_i)]}\ls\cm_s(\chi_{R'})(x)\ls\cm_s(\chi_R)(x).
\end{equation}
Combining this and \eqref{e4.22} yields $\eqref{e4.17}$.

\vspace{0.2cm}

{\it Case II.} $\ell(R_1)< \ell(P_1)$ and $\ell(R_2)<\ell(P_2)$. In
this case, for any $z_i\in
P_i'=x_{P_i}+B^{(i)}_{v_i[\ell(P_i)-1]+u_i+2\sz_i}$, we have
\begin{eqnarray}\label{e4.26}
z_i'\equiv A_i^{-v_i\ell(R_i)-u_i}z_i\in
A_i^{-v_i\ell(R_i)-u_i}x_{P_i}+B^{(i)}_{v_i
[\ell(P_i)-1-\ell(R_i)]+2\sz_i}\equiv \wz P_i.
\end{eqnarray}
Let $\wz P\equiv\wz P_1\times\wz P_2$. For any $z_i'\in \wz P_i$,
since $\ell(P_i)>\ell(R_i)$ and $-v_i,\sz_i >0$, by \eqref{e2.5} and
\eqref{e2.4}, similarly to the estimate of \eqref{e4.20}, we obtain
\begin{eqnarray}\label{e4.27}
|z_i'-A_i^{-v_i\ell(R_i)-u_i}x_{P_i}|\ls b_i^{v_i[\ell(P_i)
-\ell(R_i)]\zeta_{i,\,-}}.
\end{eqnarray}

Let $\wz e_R(z)\equiv e_R(A_1^{v_1\ell(R_1)+u_1}z_1,\,
A_2^{v_2\ell(R_2)+u_2}z_2)$. For any $z\in\rnm$ and $(\az_1,
\az_2)\in(\zz_+)^{n_1}\times(\zz_+)^{n_2}$, we have
\begin{eqnarray}\label{e4.28}
|\partial^{\az_1}_1\partial^{\az_2}_2\wz e_R(z)|\ls c_R.
\end{eqnarray}
Indeed,  for any $\gz_i\ge v_i\ell(R_i)+u_i+\sz_i$ and
$z_i\in\rr^{n_i}$, an application of chain rule yields
$$\|\partial^{\az_i}[\tz^{(i)}(A_i^{v_i\ell(R_i)+u_i-\gz_i}\cdot)]\|_{\infty}\ls 1.$$
Hence \eqref{e4.28} follows by the Cauchy-Schwarz inequality,
\eqref{e4.2} and Lemma \ref{l2.1}(iv), similarly to the estimate of
\eqref{e4.21},
\begin{align*}
|\partial^{\az_1}_1\partial^{\az_2}_2\wz e_R(z)|^2
&=\bigg|\iint_{R_+}\partial^{\az_1}_{z_1}
\partial^{\az_2}_{z_2}[\tz_{\gz_1,\,\gz_2}
(A_1^{v_1\ell(R_1)+u_1}\cdot-y_1,\,
A_2^{v_2\ell(R_2)+u_2}\cdot-y_2)](z_1,z_2)\\
&\hs \times (\psi_{\gz_1,\, \gz_2}\ast f)(y)
\,dy\,d\sz(\gz_1)\,d\sz(\gz_2)\bigg|^2\\
&\ls c_R^2\iint_{R_+} b_1^{-\gz_1}b_2^{-\gz_2}\, dy_1\,dy_2\,
d\sz(\gz_1)\,d\sz(\gz_2)\ls c_R^2.
\end{align*}
Without loss of generality we can assume that
\begin{equation}\label{e4.29}
b_1^{(s_1+1)v_1(\ell(P_1)-\ell(R_1))\zeta_{1,\,-}}
\le
b_2^{(s_2+1)v_2(\ell(P_2)-\ell(R_2))\zeta_{2,\,-}},
\end{equation}
since the other case is dealt in the same way. Let
\[
\mathcal P_{w_1}(z_1,z_2) = \sum_{|\az_1|\le s_3}
\frac {\partial_{z_1}^{\az_1} \tilde e_R(w_1,z_2)}{\az_1!}
(z_1-w_1)^{\az_1}
\]
be the Taylor polynomial of $\tilde e_R(\cdot,z_2)$ in the first
variable about $w_1 \in \rr^{n_1}$ and degree $s_3$. For any $(y,\,
(t_1,\,t_2))\in \bgz(x)\cap P_+$, by \eqref{e4.19} and
\eqref{e4.26}, the change of variables, and our hypothesis that each
$\tz^{(i)}$ has vanishing moments up to degree
$s_3=2\max(s_1,s_2)+1$ we have
\begin{align*}
|e_R\ast\tz_{t_1,\,t_2}(y)|
&= \bigg|\int_{\wz P}
\bigg\{\wz e_R(z_1,\,z_2)- \mathcal P_{A_1^{-v_1\ell(R_1)-u_1}x_{P_1}}(z_1,z_2)
\bigg\}\\
&\hs\times\prod_{i=1}^2
\tz^{(i)}_{t_i}\lf(y_i-A_i^{v_i\ell(R_i)+u_i}z_i\r)
b_i^{v_i\ell(R_i)+u_i} \, dz\bigg|\nonumber
\\
&\ls c_R \int_{\wz P}
|z_1-A_1^{-v_1\ell(R_1)-u_1}x_{P_1}|^{s_3+1}
\prod_{i=1}^2 |\tz^{(i)}_{t_i}(y_i-A_i^{v_i\ell(R_i)+u_i}z_i)|
b_i^{v_i\ell(R_i)}\, dz
\nonumber
\\
& \ls c_R b_1^{(s_3+1)v_1[\ell(P_1)-\ell(R_1)]\zeta_{1,\,-}}
\prod_{i=1}^2 b_i^{v_i[\ell(P_i)-\ell(R_i)]}
b_i^{-t_i+v_i\ell(R_i)}\nonumber\\
&\ls
c_R\prod_{i=1}^2
b_i^{(s_i+1)v_i\lf(\ell(P_i)-\ell(R_i)\r)\zeta_{i,\,-}}.\nonumber
\end{align*}
Indeed, the first estimate is a consequence of Taylor's remainder
theorem and \eqref{e4.28}, the second follows from \eqref{e4.27},
and the last follows from \eqref{e4.29} and $b_i^{t_i}\sim
b_i^{v_i\ell(P_i)}$ for $i=1,\,2$.

Since $\ell(R_1)< \ell(P_1)$ and $\ell(R_2)<\ell(P_2)$, by
\eqref{e4.23} and symmetry, we obtain $P'\subset R'''$. From this,
\eqref{e4.24}, Remark \ref{r2.1} and Lemma \ref{l2.1}(iv), it
follows that for $x\in P'$, $1=\cm_s(\chi_{R'''})(x)\ls
\cm_s(\chi_R)(x)$; see also \eqref{e4.25}. Then, combining this and
\eqref{e4.29} yields \eqref{e4.17}.

{\it Case III.} $\ell(R_1)\ge\ell(P_1)$ and $ \ell(R_2)<\ell(P_2)$.
In this case define
$$e_R^{(2)}(z_1,\,z_2)\equiv
e_R(z_1,\,A_2^{v_2\ell(R_2)+u_2}z_2).$$ For any $z\in\rnm$ and
$\az_2\in(\zz_+)^{n_2}$, similarly to the estimate of \eqref{e4.28},
we obtain
\begin{eqnarray}\label{e4.30}
|\partial^{\az_2}_{2} e_R^{(2)}(z_1,z_2)| \ls c_R.
\end{eqnarray}
Let $\bar R_1\equiv A_1^{-t_1}x_{R_1} +B^{(1)}_{v_1[\ell(R_1)-1]
+u_1+2\sz_1-t_1}$ and
$$\wz P_2\equiv A_2^{-v_2
\ell(R_2)-u_2}x_{P_2}+B^{(2)}_{v_2 [\ell(P_2)-1-\ell(R_2)]
+2\sz_2}.$$
Let $\mathcal P^{(1)}_{w_1}$ be the Taylor polynomial
of $\tz^{(1)}$ about $w_1\in\rr^{n_1}$ of degree $s_1$, and let
\[
\mathcal P_{w_2}(z_1,z_2) = \sum_{|\az_2|\le s_2}
\frac {\partial_{z_2}^{\az_2} e^{(2)}_R(z_1,w_2)}{\az_2!}
(z_2-w_2)^{\az_2}
\]
be the Taylor polynomial of $e^{(2)}_R(z_1,\cdot)$ in
the second variable about $w_2 \in \rr^{n_2}$ of degree $s_2$.
For any $(y,\,(t_1,\,t_2))\in\bgz(x)\cap P_+$, by \eqref{e4.19}, the
change of variables, and vanishing moment conditions, we have
\begin{align*}
& e_R\ast  \tz_{t_1,\,t_2}(y)\\
&\hs=\int_{\bar R_1\times\wz P_2} e_R^{(2)}(A_1^{t_1}z_1,\,
z_2)\tz^{(1)}(A_1^{-t_1}y_1-
z_1)\tz_{t_2}^{(2)}(y_2-A_2^{v_2\ell(R_2)+u_2}z_2)
b_2^{v_2\ell(R_2)+u_2}\,dz\\
&\hs=b_2^{v_2\ell(R_2)+u_2} \int_{\bar R_1\times\wz P_2}
\bigg(e_R^{(2)}(A_1^{t_1}z_1,\, z_2)
- \mathcal P_{A_2^{-v_2 \ell(R_2)-u_2}x_{P_2}}(A_1^{t_1}z_1,z_2)\bigg)\\
&\hs\hs\times\bigg(  \tz^{(1)}(A_1^{-t_1}y_1-z_1) -
\mathcal P^{(1)}_{A_1^{-t_1}y_1-A_1^{-t_1}x_{R_1}}(A_1^{-t_1}y_1-z_1) \bigg)\\
&\hs\hs\times\tz_{t_2}^{(2)}(y_2-A_2^{v_2\ell(R_2)+u_2}z_2)
\,dz_1 \,dz_2.
\end{align*}
The last equation is a consequence of Fubini's theorem with the
inside integration over the $z_2$ variable. Consequently, Taylor's
remainder theorem, \eqref{e4.20} for $i=1$, and \eqref{e4.27} for
$i=2$ yields
\begin{align}\label{e4.31}
& |e_R\ast  \tz_{t_1,\,t_2}(y)|\\
&\hs\ls c_R b_1^{v_1(\ell(R_1)-\ell(P_1))}b_2^{v_2
[\ell(P_2)-\ell(R_2)]}b_2^{-t_2+v_2\ell(R_2)}\prod_{i=1}^2
b_i^{(s_i+1)v_i|\ell(P_i)-\ell(R_i)|\zeta_{i,\,-}}\nonumber\\
&\hs\ls c_R b_1^{v_1[\ell(R_1)-\ell(P_1)]}\prod_{i=1}^2
b_i^{(s_i+1)v_i|\ell(P_i)-\ell(R_i)|\zeta_{i,\,-}}.\nonumber
\end{align}

Moreover, observing that $\ell(R_1)\ge\ell(P_1)$ and
$\ell(R_2)<\ell(P_2)$, by  $\eqref{e4.23}$ and symmetry, we obtain
that $R_1'\subset P_1'''$ and $P_2'\subset R_2'''$. From this,
$b_i^{v_i\ell(R_i)}\sim |R'_i|$, $b_i^{v_i\ell(P_i)}\sim |P_i'''|$,
$\eqref{e4.24}$, Remark \ref{r2.1} and Lemma \ref{l2.1}(iv), it
follows that
\begin{eqnarray*}
&b_1^{v_1[\ell(R_1)-\ell(P_1)]}&\sim
\dfrac{|R_1'|}{|P_1'''|}\sim\frac{|R_1'\cap
P_1'''|}{|P_1'''|}\frac{|R_2'''\cap P_2'|} {|P_2'|}\ls
\frac{|R_1'''\cap P_1'''|}{|P_1'''|}\frac{|R_2'''\cap
P_2'''|}{|P_2'''|}\\
&&\sim \frac{|R'''\cap P'''|}{|P'''|}\ls \cm_s(\chi_{R'''})(x)\ls
\cm_s(\chi_R)(x);
\end{eqnarray*}
see also \eqref{e4.25}. Combining this and \eqref{e4.31} yields
\eqref{e4.17}.

\vspace{0.2cm}

{\it Case $IV$.} Finally, the case $\ell(R_1)<\ell(P_1)$ and
$\ell(R_2)\ge\ell(P_2)$ follows from Case III by the symmetry. This
completes the proof of the crucial estimate \eqref{e4.17}.

{\bf Step 8. Verify \eqref{e4.18}}

Let $\sharp E$ be the cardinality of the set $E$. For $i=1,\,2$, by
\eqref{e4.23} and Lemma \ref{l2.1}, we have
\begin{eqnarray*}
\sharp\{R_i\in\cq^{(i)}:\, R'_i\cap P'_i\ne\emptyset,\, |R_i|\le
|P_i|,\, \ell(R_i)=k_i\}\ls \frac{|P_i'''|}{|R_i|}\sim
b_i^{-v_i|k_i-\ell(P_i)|}
\end{eqnarray*}
and
\begin{eqnarray*}
\sharp\{R_i\in\cq^{(i)}:\, R'_i\cap P'_i\ne\emptyset,\, |R_i|\ge
|P_i|,\, \ell(R_i)=k_i\}=1.
\end{eqnarray*}
Then by this and $(s_i+1)\zeta_{i,-}-1>0$, we obtain
\begin{eqnarray*}
\sum_{R'\cap P'\ne\emptyset} \prod_{i=1}^2
b_i^{(s_i+1)v_i|\ell(R_i)-\ell(P_i)|\zeta_{i,-}}&&
\ls\prod_{i=1}^2\sum_{k_i\in\zz}\sum_{\gfz{\ell(R_i)=k_i} {R_i'\cap
P_i'\ne\emptyset}}
b_i^{(s_i+1)v_i|\ell(R_i)-\ell(P_i)|\zeta_{i,\,-}}\\
&&\ls \prod_{i=1}^2\sum_{k_i\in\zz}
b_i^{v_i[(s_i+1)\zeta_{i,-}-1]|k_i-\ell(P_i)|}\ls 1,
\end{eqnarray*}
which shows \eqref{e4.18} and hence, completes the proof of Lemma
\ref{l4.1}.
\end{proof}

We now prove the converse of Lemma \ref{l4.1}.

\begin{lem}\label{l4.3}
Let the assumptions be as in Theorem \ref{t4.1}. Then there exists a
positive constant $C$ such that for all
$f\in\hpa\cap\cs'_{\fz,\,w}(\rnm)$,
$\|f\|_\hpa\le C\|f\|_\hp.$
\end{lem}

To prove Lemma \ref{l4.3}, we need a variant of the Journ\'e's
covering lemma established in \cite{j86, p86};
see also \cite{clmp06} for some different variants. We first recall some
notation and definitions. Let $\Oz\subset \rnm$ be an open set.
Denote by $m_i(\Oz)$ the family of all dyadic rectangles
$R\subset\Oz$ which are maximal in the $\rr^{n_i}$ ``direction'',
where $i=1,\ 2$. Recall that $n_1=n$ and $n_2=m$. Let
$\eta_0\in(0,\,1)$. For $R=R_1\times R_2\in m_1(\Oz)$, let $\widehat
R_2\equiv \widehat R_2(R_1)$ be the ``longest'' dyadic cube
containing $R_2$ such that
$|(R_1\times \widehat R_2) \cap\Oz|>\eta_0 |R_1\times
\widehat R_2|;$
and for $R=R_1\times R_2\in m_2(\Oz)$, let $\widehat R_1\equiv
\widehat R_1(R_2)$ be the ``longest'' dyadic cube containing
$R_1$ such that
\begin{eqnarray}\label{e4.32}
|(\widehat R_1\times R_2) \cap\Oz|>\eta_0 |\widehat
R_1\times R_2|.
\end{eqnarray}
For $R_i\in\cq^{(i)}$ and $j_i\in\nn$, we denote by $(R_i)_{j_i}$
the unique dyadic cube in $\cq^{(i)}$ containing $R_i$ with
$\ell((R_i)_{j_i})=\ell(R_i)-j_i$. Obviously, $(R_i)_0=R_i$. Also,
let $h: [0,\infty) \to [0,\infty)$ be an increasing function such that $\sum_{j=0}^\fz j
h(C_0\dz^j_0)<\fz$, where $C_0\equiv\max\{b_1^{2u_1-1},\,
b_2^{2u_2-1}\}$ and $\dz_0\equiv\max\{b_1^{v_1},\, b_2^{v_2}\}$.

 The following result is a variant of the well-known Journ\'e's
 covering lemma in \cite{p86} and is adapted to expansive dilations.

\begin{lem}\label{l4.4}
Let $A_i$ be a dilation on $\rr^{n_i}$ for $i=1,\,2$,
$w\in\ca_\fz(\vec A)$ and $q_w$ be as in \eqref{e2.7}. Let $\eta_0\in (0,1)$.
Then there exists a positive constant $C$, only depending on $n$, $m$, $\eta_0$
and $C_{q,\,\vec A,\,n,\,m}(w)$ with $q\in (q_w,\fz)$,
such that for all open sets $\Oz\subset\rnm$ with $w(\Oz)<\fz$,
\begin{eqnarray}\label{e4.33}
\sum_{R=R_1\times R_2\in m_1(\Oz)}w(R)h\bigg(\frac{|R_2|}{|\widehat
R_2|}\bigg)\le C w(\Oz)
\end{eqnarray}
and
\begin{eqnarray}\label{e4.34}
\sum_{R=R_1\times R_2\in m_2(\Oz)}w(R)h\bigg(\frac{|R_1|}{|\wh
R_1|}\bigg)\le C w(\Oz).
\end{eqnarray}
\end{lem}
\begin{proof}
Since the proofs for \eqref{e4.33} and \eqref{e4.34}  are similar,
we only show \eqref{e4.34}.

Let $R_1\in\cq^{(1)}$ such that $R_1\times R_2\in  m_2(\Oz)$
for certain $R_2$. Notice that for any given $R_1\in\cq^{(1)}$,
there may exist more than one $P_2\in\cq^{(2)}$ such
that $R_1\times P_2\in  m_2(\Oz)$. Based on this, for any $j_1\in\nn$,
we define
\begin{eqnarray}\label{e4.35}
A_{R_1,\,{j_1}}\equiv \{P_2\in \cq^{(2)}: \,R_1\times P_2\in
m_2(\Oz),\ \hat R_1\equiv\widehat R_1(P_2)=(R_1)_{{j_1}-1}\}.
\end{eqnarray}
If $A_{R_1,\,{j_1}}\ne\emptyset$,
for each $R_2\in A_{R_1,\,{j_1}}$,
then by Lemma \ref{l2.1}(iv), we have
$$x_{R_1}+B^{(1)}_{v_1\ell(R_1)-u_1}\subset R_1\subset x_{R_1}+B^{(1)}
_{v_1\ell(R_1)+u_1}$$ and
$x_{\wh R_1}+B^{(1)}_{v_1\ell(\wh R_1)-u_1}\subset
\wh R_1\subset x_{\wh R_1}+B^{(1)} _{v_1\ell(\wh R_1)+u_1}.$ From
this, it follows that

\begin{eqnarray}\label{e4.36}
b_1^{-2u_1}b_1^{v_1({j_1}-1)}\le \frac{|R_1|}{|\widehat R_1|} \le
b_1^{2u_1} b_1^{v_1({j_1}-1)}.
\end{eqnarray}
Let $\wz C\equiv b_1^{2u_1-1}$. By \eqref{e4.36} and the
disjointness of $\{R_2:\, R_1\times R_2\in m_2(\Oz)\}$, we have
\begin{eqnarray*}
&& \sum_{R=R_1\times R_2\in
m_2(\Oz)}w(R)h\bigg(\frac{|R_1|}{|\wh R_1|}\bigg)\\
&&\hs =\sum_{\{R_1:\, R_1\times R_2\in
m_2(\Oz)\}}\sum_{j_1\in\nn,\,A_{R_1,\,{j_1}}\ne\emptyset}
\sum_{R_2\in A_{R_1,\,{j_1}}} w(R_1\times
R_2)h\bigg(\frac{|R_1|}{|\wh R_1|}\bigg)
\nonumber\\
&&\hs \le\sum_{j_1\in\nn}h(\wz C b_1^{v_1 j_1})
\sum_{\{R_1:\, R_1\times R_2\in m_2(\Oz),\,
A_{R_1,\,{j_1}}\ne\emptyset\}}
 w\bigg(R_1\times\bigcup_{R_2\in
A_{R_1,\,{j_1}}}R_2\bigg)\nonumber.
\end{eqnarray*}

Set $E_{R_1}\equiv \bigcup_{R_1\times R_2\subset \Oz}R_2$. For any
$j_1\in\nn$ and any given $R_1\in\cq^{(1)}$ satisfying
$A_{R_1,\,{j_1}}\ne\emptyset$,
if $x_2\in \cup_{R_2\in A_{R_1,\, {j_1}}} R_2$, then there
exists dyadic cube $R_2\in\cq^{(2)}$ such that $R_1\times R_2\in
m_2(\Oz)$,
 $x_2\in R_2$ and $\hat
R_1=(R_1)_{{j_1}-1}$ by $\eqref{e4.35}$. By \eqref{e4.32} and the
maximality of $\hat R_1$, we have
$|((R_1)_{{j_1}-1}\times R_2 )\cap \Oz|>\eta_0|(R_1)_{{j_1}-1}\times R_2|$
and
$|((R_1)_{j_1}\times R_2 )\cap \Oz|\le
\eta_0 |(R_1)_{j_1}\times R_2|,$ which implies that
$|((R_1)_{j_1}\times R_2 )\cap ((R_1)_{j_1}\times E_{(R_1)_{j_1}})|\le\eta_0 |
(R_1)_{j_1}\times R_2|,$ namely,
$|(R_1)_{j_1}\times (R_2\cap E_{(R_1)_{j_1}})|\le\eta_0 |(R_1)_{j_1}\times R_2|.$
Therefore, $|R_2\cap E_{(R_1)_{j_1}}|\le\eta_0 |R_2|,$ and hence, $
|R_2\cap (E_{(R_1)_{j_1}})^\complement|>(1-\eta_0) |R_2|, $ where
$(E_{(R_1)_{j_1}})^\complement\equiv(\rrm\setminus E_{(R_1)_{j_1}})$.
From this and $R_2\subset E_{R_1}$, it follows that for $x_2\in
R_2$,
$\cm^{(2)}(\chi_{E_{R_1}\setminus E_{(R_1)_{j_1}}})(x_2)>1-\eta_0,$
where $\cm^{(2)}$ is the Hardy-Littlewood maximal operator  with respect
to the second variable, namely, on $\rrm$. Thus, for any $j_1\in\nn$,
we obtain
\begin{eqnarray}\label{e4.37}
\bigcup_{R_2\in A_{R_1,\,{j_1}}}R_2\subset K\equiv \lf\{x_2\in \rrm:\,
\cm_2(\chi_{E_{R_1}\setminus E_{(R_1)_{j_1}}})(x_2)>1-\eta_0\r\}.
\end{eqnarray}
Since $w\in\ca_\fz(\vec A)$ implies that there exists
$q\in(1,\,\fz)$ such that $w\in\ca_q(\vec A)$. Then by Definition
\ref{d2.4}, for almost all $x_1\in\rn$, we obtain that $w(x_1,\,
\cdot)\in \ca_q(A_2)$ and the weighted constants are uniformly
bounded. By this, \eqref{e4.37} and Proposition \ref{p2.1}(ii), we
have
\begin{eqnarray}
&w\lf(R_1\times\lf(\bigcup_{R_2\in A_{R_1,\, {j_1}}}R_2\r)\r)&\le
w(R_1\times K)\ls w\lf(R_1\times (E_{R_1}\setminus
E_{(R_1)_{j_1}})\r).\label{e4.38}
\end{eqnarray}
For $i= 1,\,\ldots,\,j_1$, by the disjointness of sets
$\{(R_1)_{i-1}\times (E_{(R_1)_{i-1}}\setminus
E_{(R_1)_i})\subset\Oz:\, R_1\in \cq^{(1)}\},$
we have
$$\sum_{\{R_1:\, R_1\times R_2\in m_2(\Oz),
\,A_{R_1,\,{j_1}}\ne\emptyset\}}w\lf((R_1)_{i-1} \times
(E_{(R_1)_{i-1}}\setminus E_{(R_1)_i})\r)\le w(\Oz).$$ By this,
$R_1\subset (R_1)_{i-1}$ for $i\in \nn$ and \eqref{e4.38}, we obtain
\begin{eqnarray*}
&& \sum_{R=R_1\times R_2\in
m_2(\Oz)}w(R)h\bigg(\frac{|R_1|}{|\wh R_1|}\bigg)\\
&&\hs\ls\sum_{{j_1}\in\nn}h(\wz C b_1^{v_1 j_1})
\sum_{\{R_1:\, R_1\times R_2\in m_2(\Oz),
\,A_{R_1,\,{j_1}}\ne\emptyset\}}
w(R_1\times (E_{R_1}\setminus E_{(R_1)_{j_1}}))\\
&&\hs \ls \sum_{{j_1}\in\nn}h(\wz C b_1^{v_1 j_1})\sum_{\{R_1:\, R_1\times R_2\in m_2(\Oz),
\,A_{R_1,\,{j_1}}\ne\emptyset\}}\sum_{i=1}^{j_1}
w((R_1)_{i-1}\times
(E_{(R_1)_{i-1}}\setminus E_{(R_1)_i}))\\
&&\hs\ls w(\Oz)\sum_{{j_1}\in\nn}{j_1} h(\wz C b_1^{v_1 j_1})\ls w(\Oz),
\end{eqnarray*}
which completes the proof of Lemma \ref{l4.4}.
\end{proof}

\begin{proof}[Proof of Lemma \ref{l4.3}]

We prove Lemma \ref{l4.3} by the following 7 steps.

{\bf Step 1. Reduce to the uniform estimates on atoms.}

Let $\psi$ be as in Proposition \ref{p2.5}. It suffices to prove
that for all $(p,\,q,\,\,\vec s)_w$-atoms $a$,
 \begin{eqnarray}\label{e4.39}
\|\vec S_\psi(a)\|_{L^p_w(\rnm)}\ls1.
\end{eqnarray}
In fact, for any $f\in\hpa$, there exist
$\{\lz_k\}_{k\in\nn}\subset\cc$ and $(p,\,q,\,\vec s)_w$-atoms
$\{a_k\}_{k\in\nn}$ such that $f=\sum_{k\in\nn}\lz_k a_k$ in
$\cs'(\rnm)$ and $\sum_{k\in\nn}|\lz_k|^p\ls \|f\|^p_\hpa.$ By this,
$\psi\in\cs(\rnm)$, Minkowski's inequality, Fatou's lemma, and the
monotonicity of the $\ell^p$-norm with $p\in (0,1]$ and
\eqref{e4.39}, we have
$$\|\vec S_\psi(f)\|^p_{L^p_w(\rnm)}\le \sum_{k\in\nn}
\lz_k^p\|\vec S_\psi(a)\|^p_{L^p_w(\rnm)}\ls \|f\|_\hpa^p.$$

Let us now show \eqref{e4.39} by Step 2 through Step 7.

{\bf Step 2. Estimate $\vec S_\psi(a)$ on a ``finite"
expansion of the support of $a$.}

Assume that $a$ is a $(p,\ q,\ \vec s)_w$-atom associated with an
open set $\Oz$ satisfying $w(\Oz)<\fz$ as in Definition \ref{d4.2}.
Let $\wz\Oz$ be as in \eqref{e4.3} and $\eta_0\equiv
b_1^{v_1-5\sz_1}b_2^{v_2-5\sz_2}$. Obviously, $\eta_0\in(0,1)$.
For each $R=R_1\times R_2\in m(\wz\Oz)$, let $\wh R_1$ be the ``longest'' dyadic cube
containing $R_1$ such that
$|(\wh R_1\times R_2)\cap\wz\Oz|>\eta_0|\wh R_1\times R_2|.$
For $\wz\Oz$, we define
$\wz\Oz'\equiv \{x\in\rnm:\, \cm_s(\chi_{\wz\Oz})(x)>b_1^{-2u_1}
b_2^{-2u_2}\eta_0\}.$
Similarly, we define $\wz\Oz''$ and $\wz\Oz'''$ by replacing
$\wz\Oz$ in the definition of $\wz\Oz'$, respectively, by $\wz\Oz'$
and $\wz\Oz''$. Obviously, $(\wh R_1\times R_2)\subset\wz\Oz'$. For
any given $\wz R_1\times R_2\in m_1(\wz\Oz')$ and $\wz R_1\supset
R_1$, let $\wh R_2$ be the ``longest'' dyadic cube containing $R_2$
such that
$|(\wz R_1\times \wh R_2)\cap \wz\Oz'|>\eta_0|\wz R_1\times \wh R_2|.$
Set $\bar R^\ast\equiv \bar R_1^\ast\times \bar R_2^\ast\equiv
(x_{R_1}+B^{(1)}_{v_1(\ell(\wh R_1)-1)+u_1+5\sz_1})\times
(x_{R_2}+B^{(2)}_{v_2(\ell(\wh R_2)-1)+u_2+5\sz_2}).$ Then we have
\begin{eqnarray}\label{e4.40}
w\bigg(\bigcup_{R\in m(\wz\Oz)}\bar R^\ast\bigg)\ls w(\Oz).
\end{eqnarray}
In fact, to prove \eqref{e4.40}, let $R^\sharp\equiv (x_{\wh
R_1}+B^{(1)}_{v_1\ell(\wh R_1)-u_1})\times (x_{\wh
R_2}+B^{(2)}_{v_2\ell(\wh R_2)-u_2})$. By Lemma \ref{l2.1}(iv) and
\eqref{e2.1}, $R^\sharp \subset(\wh R_1\times \wh R_2)\subsetneq
\bar R^\ast$ and $R^\sharp\subset \wz\Oz''$ which is deduced from
the fact that $\wz R_1\times\wh R_2\subset\wz\Oz''$ and $\wh
R_1\subset \wz R_1$. For any $R\in m(\wz\Oz)$ and $x\in\bar R^\ast$,
$$\cm_s(\chi_{\wz\Oz''})(x)\ge \frac{1}{|\bar R^\ast|}
\int_{\bar R^\ast}\chi_{\wz\Oz''}(y) \, dy>\frac{|R^\sharp|}{|\bar
R^\ast|}=b_1^{-2u_1}b_2^{-2u_2}\eta_0,$$
which implies that $\cup_{R\in m(\wz\Oz)}\bar
R^\ast\subset \wz\Oz'''$. From this, $w\in\ca_q(\vec A)$ and the
boundedness of $\cm_s$ on $L^q_w(\rnm)$
(see Proposition \ref{p2.2}(ii)), it follows that
$$w\bigg(\bigcup_{R\in m(\wz\Oz)}\bar R^\ast\bigg)\le w(\wz\Oz''')
\ls w(\Oz).$$ Thus, \eqref{e4.40} holds.

Then for $w\in\ca_\fz(\vec A)$, $p\in(0,\,1]$ and
$q\in[2,\,\fz)\cap(q_w,\,\fz)$, by H\"older's inequality, Theorem
\ref{t3.2}, \eqref{e4.40} and Definition \ref{d4.2}(II), we obtain
\begin{eqnarray}\label{e4.41}
\iint_{\cup_{R\in m(\wz\Oz)}\bar R^\ast}[\vec S_\psi(a)(x)]^pw(x)\,
dx &&\le \bigg[w\bigg(\bigcup_{R\in m(\wz\Oz)}\bar
R^\ast\bigg)\bigg]^{1-p/q}\|
\vec S_\psi(a)\|^p_{L^q_w(\rnm)}\\
&&\ls [w(\Oz)]^{1-p/q}\|a\|^p_{L^q_w(\rnm)}\ls 1.\nonumber
\end{eqnarray}

{\bf Step 3. Estimate $\vec S_\psi(a)$ on the complement of a ``
finite'' expansion of the support of $a$.}

Set $\bar R_1\equiv x_{R_1}+B^{(1)}_{v_1(\ell(R_1)-1)+u_1+5\sz_1}$
and $\bar R_2\equiv x_{R_2}+B^{(2)}_{v_2(\ell(R_2)-1)+u_2+5\sz_2}$.
Then by $a=\sum_{R\in m(\wz\Oz)}a_R$ in $\cs'(\rnm)$ as in
Definition \ref{d4.2} and the monotonicity of the $\ell^p$-norm with
$p\in (0,1]$, we obtain
\begin{align}\label{e4.42}
&\iint_{\lf( \cup_{R\in m(\wz\Oz)}\bar R^\ast\r)^\complement}
[\vec S_\psi(a)(x)]^p w(x)\, dx\\
&\hs\le \sum_{R\in
m(\wz\Oz)}\iint_{(\bar R^\ast)^\complement}
[\vec S_\psi(a_R)(x)]^p w(x)\, dx\nonumber\\
&\hs\le\sum_{R\in m(\wz\Oz)}\bigg[\iint_{(\bar
R_1^\ast)^\complement\times \bar R_2}+\iint_{(\bar
R_1^\ast)^\complement\times (\bar R_2)^\complement}+\iint_{\bar
R_1\times (\bar R_2^\ast)^\complement}\nonumber\\
&\hs\hs+\iint_{(\bar R_1)^\complement\times (\bar
R_2^\ast)^\complement} \bigg] [\vec S_\psi(a_R)(x)]^p w(x)\, dx
\equiv \sum_{R\in m(\wz\Oz)}({\mathrm K_1}+{\mathrm K_2}+{\mathrm
K_3}+{\mathrm K_4}).\nonumber
\end{align}

{\bf Step 4. Pointwise estimate of $\vec S_\psi(a_R)$ on $(\bar
R_1)^\complement\times \bar R_2^\ast$.}

Let $\gz_1(R)\equiv \ell(\wh R_1)-\ell(R_1)$, $\bar
R_{1,\,k_1}^\ast\equiv x_{R_1}+B^{(1)}_{v_1(\ell(\wh R_1)-1-k_1)
+u_1+5\sz_1}$ for $k_1\in\nn$, and $\bar {R}_{1,0}^\ast\equiv\bar
R_1^\ast$. We will prove in this step that for all $k_1\in\zz_+$
and $x=(x_1,\,x_2)$ with $x_1\in
\bar R_{1,\,k_1+1}^\ast\setminus\bar R_{1,\,k_1}^\ast$ and $x_2\in
\bar R_2$,
\begin{eqnarray}\label{e4.43}
\vec S_\psi(a_R)(x)&&\ls b_1^{[k_1-\gz_1(R)]v_1(s_1+1)\zeta_{1,\,-}}
b_1^{-v_1(\ell(\wh R_1)-k_1)} \int_{R_1''}
S_{\psi^{(2)}}(a_R(z_1,\,\cdot)) (x_2)\, dz_1,
\end{eqnarray}
where $S_{\psi^{(2)}}$ is the Lusin-area function with respect to
the second variable and $s_1$ as in Definition \ref{d4.2}.

Let $L(R_1)\equiv v_1[\ell(R_1)-1]+u_1+3\sz_1$. We now estimate
$a_R\ast\psi_{j_1,\,j_2}(x-y)$ by considering two cases, where $x$ is as in
\eqref{e4.43}, $j_1,\,j_2\in\zz$ and $y\in B^{(1)}_{j_1}\times
B^{(2)}_{j_2}$.

{\it Case I.} $j_1>L(R_1)$. For any $z_1\in R_1''\equiv
x_{R_1}+B^{(1)}_{L(R_1)}$, we have
$$z_1'\equiv A_1^{-j_1}z_1\in A_1^{-j_1}R_1''= A_1^{-j_1}x_{R_1}+
B^{(1)}_{L(R_1)-j_1}\equiv \wz R_1''.$$ Then, by $j_1>L(R_1)$ and
\eqref{e2.4}, we have
\begin{eqnarray}\label{e4.44}
&|z_1'-A_1^{-j_1}x_{R_1}|\ls b_1^{[L(R_1)-j_1] \zeta_{1,\,-}}.
\end{eqnarray}
Let
\[
\mathcal P_{w_1}(z_1)\equiv\sum_{|\az_1|\le s_1}\frac 1{\az_1!}
\partial^{\az_1} \psi^{(1)}
(w_1)(z_1-w_1)^{\az_1}
\]
be the Taylor polynomial of $\psi^{(1)}$ about $w_1 \in \rr^{n_1}$ of degree $s_1$.
Since $\supp a_R\subset R''$ and $a_R$ has vanishing moments up to
order $s_1$ for the first variable, by Taylor's remainder theorem,
we have
\begin{align*}
& |a_R\ast\psi_{j_1,\,j_2}(x-y)|\\
&\hs=\bigg| \int_{R_1''}(a_R\ast_2\psi^{(2)}_{j_2}(x_2-y_2))(z_1)
b_1^{-j_1}\psi^{(1)}(A_1^{-j_1}(x_1-y_1-z_1))
\, dz_1 \bigg|\nonumber\\
&\hs=\bigg|\int_{\wz R_1''}(a_R\ast_2\psi^{(2)}(x_2-y_2))
 (A_1^{j_1}z_1) \lf[\psi^{(1)}- \mathcal P_{A_1^{-j_1}(x_1-y_1-x_{R_1})}\r]
(A_1^{-j_1}(x_1-y_1)-z_1) \, dz_1\bigg|
\nonumber\\
&\hs\ls \mathcal D
\int_{\wz
R_1''}|(a_R\ast_2\psi^{(2)}_{j_2}
(x_2-y_2))(A_1^{j_1}z_1)|
|(A_1^{-j_1}x_{R_1}-z_1)|^{s_1+1}\, dz_1,
\nonumber
\end{align*}
where
\[
\mathcal D \equiv \sup_{|\az_1|=s_1+1}
\sup_{\xi_1\in A_1^{-j_1}x_{R_1} +B^{(1)}_{L(R_1)-j_1}}
|\partial^{\az_1}\psi^{(1)}(A_1^{-j_1}(x_1-y_1)-\xi_1)|.
\]
Since $A_1^{-j_1}x_1-A_1^{-j_1}x_{R_1}\not\in
B^{(1)}_{L(R_1)+v_1[\gz_1(R)-k_1]+2\sz_1-j_1},$
by \eqref{e2.2}, we know $A_1^{j_1}x_1-\xi_1\not\in
B^{(1)}_{L(R_1)+v_1[\gz_1(R)-k_1]+\sz_1-j_1}$. Thus,
if $j_1\le L(R_1)+v_1[\gz_1(R)-k_1]$, by $A_1^{-j_1}y_1\in
B_0^{(1)}$ and \eqref{e2.2}, we have
$A_1^{-j_1}(x_1-y_1)-\xi_1\not\in
B^{(1)}_{L(R_1)+v_1[\gz_1(R)-k_1]-j_1}.$
This together with $\psi^{(1)}\in\cs(\rn)$ and \eqref{e2.3} yields
that
\begin{eqnarray}\label{e4.45}
&\mathcal D&\ls
\sup_{\xi_1\in A_1^{-j_1}x_{R_1} +B^{(1)}_{L(R_1)-j_1}}
[1+\rho_1(A_1^{-j_1}(x_1-y_1)-\xi_1)]^{-N_1}\\
&&\ls [1+b_1^{L(R_1)+v_1(\gz_1(R)-k_1)-j_1}]^{-N_1} \nonumber
\end{eqnarray}
for any given $N_1>0$. The same estimate also holds trivially for $j_1>
L(R_1)+v_1[\gz_1(R)-k_1]$ since $\mathcal D \ls 1$. Combining
\eqref{e4.44} through \eqref{e4.45} yields
\begin{align}\label{e4.46}
&|a_R\ast\psi_{j_1,\,j_2}(x-y)|
\ls I(j_1)\int_{R_1''}|(a_R\ast_2
\psi^{(2)}_{j_2}(x_2-y_2))(z_1)|\, dz_1,\\
&\hs\text{where }\quad  I(j_1)\equiv b_1^{-j_1}
\lf[1+b_1^{L(R_1)+v_1(\gz_1(R)-k_1)-j_1}\r]^{-N_1}
b_1^{(s_1+1)[L(R_1)-j_1]
\zeta_{1,\,-}}.\nonumber
\end{align}
Observe also that by choosing $N_1>s_1+2$, which implies that
$N_1>(s_1+1)\zeta_{1,\,-}+1$, we have
\begin{equation}\label{e4.47}
\sum_{j_1>L(R_1)} I(j_1)^2
\ls b_1^{-2v_1(\ell(\wh R_1)-k_1)}
b_1^{2[k_1-\gz_1(R)]v_1(s_1+1)\zeta_{1,\,-}}.
\end{equation}

{\it Case II.} $j_1\le L(R_1)$. In this case, for
$x_1-x_{R_1}\not\in B^{(1)}_{L(R_1)+v_1(\gz_1(R)-k_1) +2\sz_1}$ and
$z_1-x_{R_1}\in B^{(1)}_{L(R_1)}$, since $-v_1,\,-\gz_1(R),\,k_1\ge
0$ and $j_1\le L(R_1)$, by \eqref{e2.2}, we obtain
$A_1^{-j_1}(x_1-z_1)\not\in B^{(1)}_{L(R_1)+
v_1(\gz_1(R)-k_1)+\sz_1-j_1}.$ From this, $A_1^{-j_1}y_1\in
B^{(1)}_0$ and \eqref{e2.2}, we deduce
$A_1^{-j_1}(x_1-y_1-z_1)\not\in
B^{(1)}_{L(R_1)+v_1(\gz_1(R)-k_1)-j_1}$, and hence,
\begin{eqnarray}\label{e4.48}
\rho_1 (A_1^{-j_1}(x_1-y_1-z_1) )\ge
b_1^{L(R_1)+v_1(\gz_1(R)-k_1)-j_1}.
\end{eqnarray}
Choosing $N_1\ge s_1+2$,   we have $b_1b_1^{(s_1+1)\zeta_{1,\,-}}\le
b_1^{N_1}$. Since $\supp a_R\subset R''$ and
$\psi^{(1)}\in\cs(\rrm)$, by \eqref{e4.48} and H\"{o}lder's
inequality, we have
\begin{eqnarray}\label{e4.49}
&&|a_R\ast\psi_{j_1,\,j_2}(x-y)|\\
&&\hs=\lf|\int_{R_1''} (a_R\ast_2\psi_{j_2}^{(2)}(x_2-y_2))(z_1)
b_1^{-j_1}\psi_1^{(1)}(A_1^{-j_1}(x_1-y_1-z_1))\,dz_1\r|
\nonumber\\
&&\hs\ls b_1^{-j_1}b_1^{-N_1[L(R_1)+v_1(\gz_1(R)-k_1)-j_1]}
\int_{R_1''} |(a_R\ast_2\psi_{j_2}^{(2)}
(x_2-y_2))(z_1) |\, dz_1\nonumber\\
&&\hs\ls I(j_1)\int_{R_1''}\lf|\lf(a_R\ast_2\psi_{j_2}^{(2)}
(x_2-y_2)\r)(z_1)\r|\, dz_1,\nonumber
\\
&&\hs\hs\text{where } \quad I(j_1) \equiv b_1^{-v_1(\ell(\wh R_1)-k_1)}
b_1^{-[L(R_1)+v_1(\gz_1(R)-k_1)-j_1](s_1+1)\zeta_{1,\,-}}.
\nonumber
\end{eqnarray}
Observe also that we have
\begin{equation}\label{e4.50}
\sum_{j_1\le L(R_1)}
[I(j_1)]^2 \ls b_1^{-2v_1(\ell(\wh R_1)-k_1)}
b_1^{2[k_1-\gz_1(R)]v_1(s_1+1)\zeta_{1,\,-}}.
\end{equation}
Therefore, \eqref{e4.43} follows by \eqref{e4.46}, \eqref{e4.47},
\eqref{e4.49}, \eqref{e4.50}, and Minkowski's inequality
\begin{eqnarray*}
&&\vec S_\psi(a_R)(x)\\
&&\hs\ls \bigg\{\bigg(\sum_{j_1 \in \zz} [I(j_1)]^2
\bigg)
\sum_{j_2\in\zz}b_2^{-j_2}\int_{B^{(2)}_{j_2}}
\bigg[\int_{R''_1} | (a_R\ast_2\psi^{(2)}_{j_2}(x_2-y_2) )(z_1) |\,
dz_1\bigg]^2\, dy_2\bigg\}^{1/2}\\
&&\hs\ls b_1^{-v_1(\ell(\wh R_1)-k_1)}
b_1^{[k_1-\gz_1(R)]v_1(s_1+1)\zeta_{1,\,-}}
\int_{R_1''}S_{\psi^{(2)}}(a_R(z_1,\,\cdot)) (x_2) dz_1.
\end{eqnarray*}

{\bf Step 5. Estimate for ${\mathrm K_1}$.}

Since $s_1\ge \lfloor(q_w/p-1)\zeta_{1,\,-}^{-1}\rfloor$, there
exists $r\in(q_w,\, q]$ such that $p(s_1+1)\zeta_{1,\,-}+p-r>0$ and
$w\in\ca_r(\vec A)$. Recall that $\cm^{(1)}$ denotes the
Hardy-Littlewood maximal operator on $\rn$. Then, by \eqref{e4.43},
$\supp a_R\subset R''$, $w\in\ca_r(\vec A)$, the
$L^r_w(\rn)$-boundedness of $\cm^{(1)}$, Theorem \ref{t3.1} and
H\"{o}lder's inequality, we obtain
\begin{eqnarray*}
&&b_1^{-[k_1-\gz_1(R)]v_1(s_1+1)\zeta_{1,\,-}}\bigg[\int_{\bar R_{1,\,k_1+1}^\ast\setminus\bar R_{1,\,k_1}^\ast}
\int_{\bar R_2} [\vec S_\psi(a_R)(x)]^rw(x)\, dx
\bigg]^{1/r}\\
&&\hs \ls
b_1^{-v_1(\ell(\wh R_1)-k_1)}
 \bigg\{\int_{\bar R_{1,\,k_1+1}^\ast\setminus\bar
R_{1,\,k_1}^\ast} \int_{\bar
R_2}\bigg[\int_{R_1''}S_{\psi^{(2)}}(a_R(z_1,\,\cdot))
(x_2)\, dz_1 \bigg]^r w(x)\,dx\bigg\}^{1/r}\\
&&\hs \ls
\bigg\{\int_{\bar R_{1,\,k_1+1}^\ast\setminus\bar
R_{1,\,k_1}^\ast} \int_{\bar
R_2}\bigg[\cm^{(1)}\lf(S_{\psi^{(2)}}(a_R)(x_2)\r)(x_1)\bigg]^rw(x)\,dx\bigg\}^{1/r}\\
&&\hs \ls
\bigg\{\int_\rn\int_{\bar
R_2}\lf(S_{\psi^{(2)}}(a_R(x_1,\,\cdot))(x_2)\r)^rw(x)
\,dx\bigg\}^{1/r}\\
&&\hs \ls
\bigg\{\int_{R''}|a_R(x)|^rw(x)\,dx\bigg\}^{1/r}
\ls
\|a_R\|_{L^q_w(\rnm)}[w(R)]^{1/r-1/q}.
\end{eqnarray*}
From this, $v_1<0$, $p(s_1+1)\zeta_{1,\,-}+p-r>0$ and
$w\in\ca_r(\vec A)$,  H\"older's inequality and Lemma
\ref{l2.1}(iv), it follows that
\begin{align}\label{e4.51}
{\mathrm K_1}&=\sum_{k_1=0}^\fz \int_{\bar
R_{1,\,k_1+1}^\ast\setminus \bar R_{1,\,k_1}^\ast}\int_{\bar
R_2}[\vec S_\psi(a_R)(x)]^p w(x)\, dx\\
&\le \sum_{k_1=0}^\fz  \lf[w\lf(\bar R_{1,\, k_1+1}^\ast\times \bar
R_2\r)\r]^{1-p/r}\bigg[\int_{\bar R_{1,\,k_1+1}^\ast\setminus \bar
R_{1,\,k_1}^\ast}\int_{\bar R_2}[\vec S_\psi(a_R)(x)]^r w(x)
\, dx\bigg]^{p/r} \nonumber\\
&\ls \sum_{k_1=0}^\fz b_1^{v_1(\gz_1(R)-k_1)(r-p)}
[w(R)]^{1-p/r}
 b_1^{p[k_1-\gz_1(R)]v_1(s_1+1)\zeta_{1,\,-}}\nonumber\\
&\hs\times\|a_R\|_{L^q_w(\rnm)}^p [w(R)]^{p/r-p/q} \nonumber\\
&\ls  [w(R)]^{1-p/q}\|a_R\|^p_{L^q_w(\rnm)}b_1^{v_1\gz_1(R)
[r-p-p(s_1+1)\zeta_{1,\,-}]}. \nonumber
\end{align}

{\bf Step 6. Estimate for $\sum_{\wh R\in m(\wz\Oz)}({\mathrm
K_1}+{\mathrm K_2})$.}

Observe that the integral in ${\mathrm K_2}$ is on the domain $(\bar
R_1)^\complement\times (\bar R_2^\ast)^\complement$ and the integral
in ${\mathrm K_1}$ is on the domain $(\bar R_1)^\complement\times
\bar R_2^\ast$. Thus, applying the ideas used in the estimate of
${\mathrm K_1}$ on the first variable to both variables of ${\mathrm
K_2}$, we also have
\begin{eqnarray*}
{\mathrm K_2}\ls
[w(R)]^{1-p/q}\|a_R\|^p_{L^q_w(\rnm)}b_1^{v_1\gz_1(R)
[r-p-p(s_1+1)\zeta_{1,\,-}]}.
\end{eqnarray*}

Take $h_1(t)\equiv t^{a_1}$ for $t\in (0,\,1)$ and $a_1\equiv
p+p(s_1+1)\zeta_{1,\,-}-r$. Then, by $a_1>0$, we obtain that
$\sum_{j\ge 0}jh_1(t^j)^{q/(q-p)}<\fz$. By Lemma \ref{l2.1}(iv), we
have
$$b_1^{v_1\gz_1(R) [r-p-p(s_1+1)\zeta_{1,\,-}]}\sim
h_1\bigg(\frac{|R_1|}{|\wh R_1|}\bigg).$$ From this, Definition
\ref{d4.2}(II), H\"{o}lder's inequality, Lemma \ref{l4.4} and
Proposition \ref{p2.2}(ii) with $w\in\ca_q(\vec A)$, it follows that
\begin{eqnarray*}
\sum_{R\in m(\wz\Oz)}({\mathrm K_1}+{\mathrm K_2}) &&\ls\sum_{R\in
m(\wz\Oz)}
\|a_R\|^p_{L^q_w(\rnm)}[w(R)]^{1-p/q}h_1\bigg(\frac{|R_1|} {|\wh
R_1|}\bigg)\nonumber
\\
&&\ls \bigg\{\sum_{R\in m(\wz\Oz)}\|a_R\|^q_{L^q_w(\rnm)}\bigg\}^{p/q}
\bigg\{ \sum_{R\in m(\wz\Oz)}w(R)h_1\lf(\frac{|R_1|}
{|\wh R_1|}\r)^{\frac{q}{q-p}} \bigg\}^{1-p/q}\nonumber
\\
&&\ls [w(\Oz)]^{p/q(1-q/p)}[w(\wz\Oz)]^{1-p/q}\ls1.\nonumber
\end{eqnarray*}

{\bf Step 7. Estimate for $\sum_{\wh R\in m(\wz\Oz)}({\mathrm
K_3}+{\mathrm K_4})$.}

To estimate ${\mathrm K_3}$ and ${\mathrm K_4}$, notice that if
$I_i^{(1)}\times R_2\in m(\wz\Oz)$ for $i=1,\,2$,  then either
${\mathrm I_1}^{(1)}={\mathrm I_2}^{(2)}$ or ${\mathrm
I_1}^{(1)}\cap {\mathrm I_2}^{(1)}=\emptyset$. Recall that for any
$R_1\times R_2\in m(\wz\Oz)$, then $\wh R_2=\wh R_2(\wz R_1)$, where
$\wz R_1\times R_2\in m_1(\wz\Oz')$ and $\wz R_1\supset R_1$. Thus,
we have
\begin{eqnarray*}
\sum_{R\in m(\wz\Oz)}w(R)h_2\bigg(\frac{|R_2|}{|\wh
R_2|}\bigg)^{\frac{q}{q-p}}
&\le&\sum_{\wz R_1\times R_2\in
m_1(\wz\Oz')}\sum_{\gfz{R=R_1\times R_2\in
m(\wz \Oz)}{  R_1\subset \wz R_1  }}w(R)h_2\bigg(\frac{|R_2|}{|\wh R_2|}\bigg)^{\frac{q}{q-p}}
\\
&\le&\sum_{\wz R_1\times R_2\in m_1(\wz\Oz')}w(\wz R_1\times
R_2)h_2\bigg(\frac{|R_2|}{|\wh R_2|}\bigg)^{\frac{q}{q-p}},
\end{eqnarray*}
where $h_2(t)\equiv t^{a_2}$ for $t\in (0,\,1)$ and $a_2\equiv
p+p(s_2+1)\zeta_{2,\,-}-r$. From this, Lemma \ref{l4.4} and an
argument similar to the estimate for $\sum_{R\in m(\wz\Oz)}({\mathrm
K_1}+{\mathrm K_2})$, we deduce
 $\sum_{R\in m(\wz\Oz)} ({\mathrm K_3}+{\mathrm K_4})\ls 1.$
This together with \eqref{e4.41} implies \eqref{e4.39} and thus
completes the proof of Lemma  \ref{l4.3}.
\end{proof}

\begin{lem}\label{l4.5}
Let the assumptions be as in Theorem \ref{t4.1}. Then
$\hpa\subset\cs'_{\fz,\,w}(\rnm).$
\end{lem}

To prove Lemma \ref{l4.5}, for $i=1,\,2$ and $N_i\in\zz_+$, we let
 $\vec N\equiv(N_1,\,N_2)$.
Set
\begin{eqnarray*}
\mathscr{S}_{N_i}(\rr^{n_i})\equiv\Big\{&&\vz^{(i)}\in\cs(\rr^{n_i}):\\
&&\ \|\vz^{(i)}\|_{\mathscr{S}_{N_i}(\rr^{n_i})}\equiv
\dsup_{x_i\in\rr^{n_i}}\dsup_{|\az_i|\le N_i}|\partial^{\az_i}
\vz^{(i)}(x_i)|[1+\rho_i(x_i)]^{N_i}\le 1\Big\},
\end{eqnarray*}
and denote by $ \mathscr{S}_{\vec N}(\rnm)$ the collection of all
$\vz$ such that $\vz(x)=\vz^{(1)}(x_1)\vz^{(2)}(x_2)$ for  all
$x=(x_1,\,x_2)\in\rnm$ and all
$\vz^{(i)}\in\mathscr{S}_{N_i}(\rr^{n_i}).$

For any $f\in\cs'(\rnm)$ and $x\in\rnm$, we define the {\it grand
maximal function} $\cm_{\vec N}(f)(x)$ of $f$ by
$$\cm_{\vec N}(f)(x) \equiv\sup_{\vz\in\mathscr{S}_{\vec N}(\rnm)}
\sup_{k_1,\, k_2\in\zz}|f\ast\vz_{k_1,\,k_2}(x)|.
$$
Notice that if $N_1,\,N_2\ge 2$, then for all locally integrable
functions $f$ and $x\in\rnm$, $\cm_{\vec N}(f)(x)\ls\cm_s(f)(x)$.
Thus if $w\in \ca_p(\vec A)$ with $p\in(1,\, \fz)$, then $\cm_{\vec
N}$ is bounded on $L^p_w(\rnm)$. Moreover, we have the following
proposition.

\begin{prop}\label{p4.1}
Let the assumptions be as in Theorem \ref{t4.1}. If $N_i\ge s_i+2$
for $i=1,\,2$, then $\cm_{\vec N}$ is bounded from $\hpa$ to
$L^p_w(\rnm)$.
\end{prop}

Then Lemma \ref{l4.5} follows from Proposition \ref{p4.1}.

\begin{proof}[Proof of Lemma \ref{l4.5}]
Fix $\vz\in\cs(\rnm)$. Let $\vz_y(x)=\vz(x+y)$ for all $x\in\rnm$
and $y\in B^{(1)}_0\times B^{(2)}_0$. Notice that there exists a
positive constant $C$, depending on $\vz$, such that
$C\vz_y\in\mathscr{S}_{\vec N}(\rnm)$ for all $y\in B^{(1)}_0\times
B^{(2)}_0$. If $a$ is $(p,\,q,\,\vec s)_w$-atom, then  for
$j_1,\,j_2\in\nn$ and $w\in \ca_q(\vec A)$, by Proposition
\ref{p4.1} and Proposition \ref{p2.2}(i), we  have
\begin{eqnarray*}
|a\ast\vz_{j_1,\,j_2}(x)|^p&&\ls\inf_{y\in
B^{(1)}_{j_1}\times B^{(2)}_{j_2}}[\cm_{\vec N}(a)(x-y)]^p\\
&&\ls \frac 1{w(x+B^{(1)}_{j_1}\times
B^{(2)}_{j_2})}\int_{x+B^{(1)}_{j_1}\times B^{(2)}_{j_2}}
[\cm_{\vec N}(a)(y)]^p w(y)\, dy\le C_{x,\,w}
b_1^{\frac{-j_1}q}b_2^{-\frac{j_2}q},
\end{eqnarray*}
where $C_{x,\,w}$ is a positive constant independent of $j_1$ and
$j_2$, and the atom $a$. If $f=\sum_{k\in\zz}\lz_ka_k$ in
$\cs'(\rnm)$, where $a_k$ is $(p,\,q,\,\vec s)_w$-atom and
$\sum_{k\in\zz}|\lz_k|^p<\fz$, then
$$|f\ast\vz_{j_1,\,j_2}(x)|^p
\le C_{x,\,w}
b_1^{\frac{-j_1}q}b_2^{-\frac{j_2}q}\sum_{k\in\zz}|\lz_k|^p\to0$$ as
$j_1,\,j_2\to\fz$, which completes the proof of Lemma \ref{l4.5}.
\end{proof}

Finally we prove Proposition \ref{p4.1}.

\begin{proof}[Proof of Proposition \ref{p4.1}]
The proof of Proposition \ref{p4.1} is similar to that of Lemma
\ref{l4.3}. By a reason similar to that used in Step 1
of the proof of Lemma \ref{l4.3}, it
suffices to show that $\|\cm_{\vec N}(a)\|_{L^p_w(\rnm)}\ls 1$ for
all $(p,\,q,\,\vec s)_w$-atoms $a$.

Assuming that $a=\sum_{R\in m(\wz\Oz)}a_R$ is a $(p,\ q,\ \vec
s)_w$-atom associated with open set $\Oz$ with $w(\Oz)<\fz$ as in
Definition \ref{d4.2}. Let all the notation be as in the proof of
Lemma \ref{l4.3}. Similarly to the proof of \eqref{e4.41}, using the
$L^q_w(\rnm)$-boundedness of $\cm_{\vec N}$ (see Proposition
\ref{p2.2}(ii)), we have
\begin{eqnarray*}
\int_{\cup_{R\in m(\wz\Oz)}\bar R^\ast} [\cm_{\vec N}(a)(x)]^p
w(x)\, dx\ls 1.
\end{eqnarray*}
And similarly to the proof of \eqref{e4.42}, we write
\begin{eqnarray*}
&&\iint_{\lf( \bigcup_{R\in m(\wz\Oz)}\bar R^\ast\r)^\complement}
[\cm_{\vec N}(a)(x)]^p w(x)\, dx\\
&&\hs\le \sum_{R\in
m(\wz\Oz)}\bigg[\iint_{(\bar R_1^\ast)^\complement\times \bar R_2}+
\iint_{(\bar R_1^\ast)^\complement\times (\bar
R_2)^\complement}+\iint_{\bar R_1\times
(\bar R_2^\ast)^\complement}\\
&&\hs\hs+\iint_{(\bar R_1)^\complement\times (\bar
R_2^\ast)^\complement} \bigg] [\cm_{\vec N}(a_R)(x)]^p w(x) dx
\equiv \sum_{R\in m(\wz\Oz)}({\mathrm J_1}+{\mathrm J_2}+{\mathrm
J_3}+{\mathrm J_4}).
\end{eqnarray*}

For any $\psi\in\mathscr{S}_{\vec N}(\rnm)$, $x_1\in \bar
R_{1,\,k_1+1}^\ast\setminus\bar R_{1,\,k_1}^\ast$ with $k_1\in
\zz_+$, $x_2\in \bar R_2$ and $y\in B^{(1)}_{j_1}\times
B^{(2)}_{j_2}$ with $j_1,\,j_2\in \zz$, similarly to the proofs of
Case I and Case II in the proof of Lemma \ref{l4.3}, we have that

(I) if $j_1>L(R_1)$, then $|a_R\ast\psi_{j_1,\,j_2}(x-y)|$ has the
same upper estimate as in \eqref{e4.46} and
\begin{eqnarray*}
&&\sup_{j_1>L(R_1)} b_1^{-j_1}
(1+b_1^{L(R_1)+v_1(\gz_1(R)-k_1)-j_1})^{-N_1}
b_1^{(s_1+1)[L(R_1)-j_1]
\zeta_{1,\,-}}\\
&&\hs\ls b_1^{-v_1(\ell(\wh R_1)-k_1)}
b_1^{[k_1-\gz_1(R)]v_1(s_1+1)\zeta_{1,\,-}}.
\end{eqnarray*}
Here, unlike the calculation of \eqref{e4.47}, we only need
$N_1\ge s_1+2$;

(II) if $j_1\le L(R_1)$, then $|a_R\ast\psi_{j_1,\,j_2}|$ has the
same upper estimate as in \eqref{e4.49} and
\begin{eqnarray*}
&&\sup_{j_1\le L(R_1)} b_1^{-v_1(\ell(\wh R_1)-k_1)}
b_1^{-[L(R_1)+v_1(\gz_1(R)-k_1)-j_1](s_1+1)\zeta_{1,\,-}}\\
&&\hs\ls b_1^{-v_1(\ell(\wh R_1)-k_1)}
b_1^{[k_1-\gz_1(R)]v_1(s_1+1)\zeta_{1,\,-}}.
\end{eqnarray*}
Then similarly to the estimate of \eqref{e4.43}, by $N_1\ge s_1+2$,
we have
\begin{eqnarray*}
\cm_{\vec N}(a_R)(x)\ls b_1^{-v_1(\ell(\wh R_1)-k_1)}
b_1^{[k_1-\gz_1(R)]v_1(s_1+1)\zeta_{1,\,-}}
\int_{R_1''}\cm_{N_2}^{(2)}(a_R(z_1,\,\cdot)) (x_2)\, dz_1,
\end{eqnarray*}
where
$$\cm_{N_2}^{(2)}(g)(x_2)=
\sup_{\psi^{(2)}\in\mathscr{S}_{N_2}(\rr^{n_2})}\sup_{k_2\in\zz}\lf|
\lf(\psi^{(2)}_{k_2}\ast g\r)(x_2)\r|.$$ Observing that for
$s\in(1,\,\fz)$ and $\nu\in\ca_s(A_2)$, $\cm_{N_2}^{(2)}$ is bounded
on $L^s_{\nu}(\rn)$. Then similarly to the estimate of
\eqref{e4.51}, we obtain
\begin{eqnarray*}
{\mathrm J_1}\ls
[w(R)]^{1-p/q}\|a_R\|^p_{L^q_w(\rnm)}b_1^{v_1\gz_1(R)
[r-p-p(s_1+1)\zeta_{1,\,-}]}.\nonumber
\end{eqnarray*}
Also, similarly to the proof in Step 6 of the proof of Lemma
\ref{l4.3}, we also have
\begin{eqnarray*}
{\mathrm J_2}\ls
[w(R)]^{1-p/q}\|a_R\|^p_{L^q_w(\rnm)}b_1^{v_1\gz_1(R)
[r-p-p(s_1+1)\zeta_{1,\,-}]},
\end{eqnarray*}
and $\sum_{R\in m(\wz\Oz)}({\mathrm J_1}+{\mathrm J_2})\ls1.$
Finally, similarly to the proof in Step 7 of the proof of Lemma
\ref{l4.3}, we obtain $\sum_{R\in m(\wz\Oz)}({\mathrm J_3}+{\mathrm
J_4})\ls1,$ which completes the proof of Proposition \ref{p4.1}.
\end{proof}

\begin{rem}\label{r4.2}
Let $w\in\ca_\fz(\vec A)$ and $(p,\,q,\,s)_w$ be an admissible
triplet. By Proposition \ref{p4.1} and Theorem \ref{p4.1}, for
$N_i\ge s_i+2$ with $i=1,\,2$, we obtain the boundedness of
$\cm_{\vec N}$ from $\hp$ to $L^p_w(\rnm)$.
\end{rem}




\section{Weighted finite atomic Hardy spaces}\label{s5}

\hskip\parindent In this section we establish finite atomic decomposition of the anisotropic product Hardy spaces.

\begin{defn}\label{d5.1}
Let $w\in \ca_\fz(\vec A)$, $q_w$ be as in \eqref{e2.7} and $(p,\
q,\, \vec s)_w$ be an admissible triplet as in Definition
\ref{d4.2}. Let $a$ be a $(p,\ q,\, \vec s)_w$-atom associated with an
open set $\Oz$. We say $a$ is a $(p,\, q,\, \vec s)_w^\ast$-atom if
$a\in\cs(\rnm)$,   $\Oz$ is bounded, and there exist only finitely many
$R\in m(\wz\Oz)$ such that $a_R\ne0$.

The {\it weighted finite Hardy space} $\hpaf$
is defined to be the space of all functions $f
=\sum_{j=1}^k\lz_ja_j$, where $k\in\nn$,
$\{a_j\}_{j=1}^k$ are $(p,\ q,\, \vec s)^\ast_w$-atoms
and $\{\lz_j\}_{j=1}^k\subset\cc$. The norm of
$f$ is defined by
$\|f\|_\hpaf
\equiv\inf\{(\sum_{j=1}^k|\lz_j|^p)^{1/p}\},$
where the infimum is taken over all the above finite decompositions of $f$.
\end{defn}

The main result of this section is as follows.

\begin{thm}\label{t5.1}
Let $w\in \ca_\fz(\vec A)$, $q_w$ be as in \eqref{e2.7}, $(p,\ q,\,
\vec s)_w$ be an admissible triplet as in Definition \ref{d4.2}. Then,
\begin{enumerate}

\item[(i)] $\hpaf$ is dense in $\hp$.

\item[(ii)] Moreover, if $\vec s\equiv(s_1,\,s_2)$ satisfies
\begin{align}\label{e5.1}
s_1 &>[(q_w/p)-1+(q_w/p)(v_2/v_1) (\log_{b_1}b_2)]\zeta_{1,\,-}^{-1}-1\
and\\
\label{e5.2}
s_2 &>[(q_w/p)-1 +(q_w/p)(v_1/v_2)(\log_{b_2}b_1)]\zeta_{2,\,-}^{-1}-1,
\end{align}
then $\|f\|_\hpaf\sim \|f\|_\hp$ for all $f\in \hpaf$.
\end{enumerate}
\end{thm}

\begin{rem}\label{r5.1}
Notice that comparing with the non-product case (see \cite{blyz,msv,gly}),
we need additional assumptions \eqref{e5.1} and \eqref{e5.2}
on vanishing moments of atoms in Theorem \ref{t5.1}(ii).
This is due to the fact that the product Hardy space is not just a product
of one-parameter Hardy spaces.
\end{rem}

To prove Theorem \ref{t5.1}, we need the following auxiliary lemma,
which generalizes Lemma 2 and Lemma 4 in
Appendix (III) of \cite{fjw}. Lemma \ref{l5.1} below can be also deduced
with some effort from \cite[Lemma 6.3]{bh}.

\begin{lem}\label{l5.1}
Let $A $ be a dilation on $\rn$, $s\in \zz_+$ and $M\in [0,\,\fz)$.

(i) If $g\in \cs(\rn)$ and $\psi\in \cs_s(\rn)$, then there exists a
positive constant $C$ such that for all $k\in\zz\setminus\nn$ and
all $x\in\rn$,
$|(g\ast\psi_k)(x) |\le Cb^{k(s+1)\zeta_-}[1+\rho(x)]^{-M}.$

(ii) If $g\in \cs_s(\rn)$ and $\psi\in \cs(\rn)$, then there exists
a positive constant $C$ such that for all $k\in\zz_+$ and all
$x\in\rn$,
$|(g\ast\psi_k)(x)|\le C
b^{-k[(s+1)\zeta_-+1]}[1+b^{-k}\rho(x)]^{-M}.$
\end{lem}

\begin{proof}
To prove (i), let $k\in\zz\setminus\nn$. Since $\psi\in \cs_s(\rn)$,
for all $x\in\rn$, we have
\begin{eqnarray*}
(g\ast \psi_k)(x)
&&=\bigg[\int_{\rho(y)\le
\rho(x)/(2b^\sz)}+\int_{\rho(y)>\rho(x)/(2b^\sz)}\bigg]
\psi_k(y)\bigg[g(x-y)-\sum_{|\az|\le s}\frac{\partial^\az
g(x)}{\az!}(-y)^\az\bigg]\,
dy\\
&&\equiv {\mathrm I_1}+{\mathrm I_2}.
\end{eqnarray*}
For ${\mathrm I_1}$, since $g\in\cs(\rn)$, by Taylor's remainder
theorem, we have
\begin{eqnarray*}
\lf| g(x-y)-\sum_{|\az|\le s}\frac1{\az!}\partial^\az
g(x)(-y)^\az\r| &&\ls |y|^{s+1}
\sup_{|\az|=s+1,\ \rho(z)\le \rho(x)/(2b^\sz)}|\partial^\az g(x-z)|\\
&&\ls |y|^{s+1}\sup_{\rho(z)\le \rho(x)/(2b^\sz)}[1+\rho(x-z)]^{-M}.
\end{eqnarray*}
This together with $\rho(x-z)\ge \rho(x)/b^\sz-\rho(z)\ge
\rho(x)/(2b^\sz)$, \eqref{e2.3}, \eqref{e2.4}, $k\le 0$ and
$\psi\in\cs_s(\rn)$, yields
\begin{eqnarray*}
|{\mathrm I_1}|&&\ls [1+\rho(x)]^{-M}\lf\{\int_{\rho(y)\le
1}\rho(y)^{(s+1)\zeta_-}|\psi_k(y)|\,
dy+\int_{\rho(y)>1}\rho(y)^{(s+1)\zeta_+}|\psi_k(y)|\, dy\r\}\\
&& \ls
b^{k(s+1)\zeta_-}[1+\rho(x)]^{-M}\lf\{\int_{\rn}[\rho(y)^{(s+1)\zeta_-}+
\rho(y)^{(s+1)\zeta_+}]|\psi(y)|\, dy\r\}\\
&&\ls b^{k(s+1)\zeta_-}[1+\rho(x)]^{-M}.
\end{eqnarray*}
For ${\mathrm I_2}$, if $\rho(x)>1$, since $g\in\cs(\rn)$ and $k\le
0$, by Taylor's remainder theorem, \eqref{e2.3} and \eqref{e2.4}, we
have
\begin{eqnarray*}
|{\mathrm I_2}|&&\ls
\int_{\rho(y)>\rho(x)/(2b^\sz)}|y|^{s+1}\sup_{|\az|=s+1}\|\partial^\az
g\|_{L^\fz(\rn)}|\psi_k(y)|\,dy\\
&&\ls
\int_{\rho(x)/(2b^\sz)<\rho(y)<1}\rho(y)^{(s+1)\zeta_-}|\psi_k(y)|\,dy
+\int_{\rho(x)/(2b^\sz)<\rho(y),\,\rho(y)>1}\rho(y)^{(s+1)\zeta_+}
|\psi_k(y)|\,dy\\
&&\ls
b^{k(s+1)\zeta_-}\int_{\rho(y)>b^{-k}\rho(x)/(2b^\sz)}[\rho(y)^{(s+1)\zeta_-}+
\rho(y)^{(s+1)\zeta_+}]|\psi(y)|\,dy\ls b^{k(s+1)\zeta_-}\rho(x)^{-M}.
\end{eqnarray*}
By this and $\rho(x)>1$, we have
$|{\mathrm I_2}|\ls b^{k(s+1)\zeta_-}[1+\rho(x)]^{-M}.$
For $\rho(x)\le 1$, similarly to the above estimate, we obtain
$|{\mathrm I_2}|\ls b^{k(s+1)\zeta_-}\ls b^{k(s+1)\zeta_-}[1+\rho(x)]^{-M}.$
Combining above estimates for ${\mathrm I_1}$ and ${\mathrm I_2}$
completes the proof Lemma \ref{l5.1}(i).

To prove (ii), we observe the identity $g* \psi_k = (g_{-k} * \psi)_k$.
Thus, if $k\in \zz_+$, then (i) with the roles of $g$ and $\psi$ exchanged yields
$$
|g* \psi_k (x)|= |(g_{-k} * \psi)_k(x)| \ls b^{-k(s+1) \zeta_-}
[1+ \rho(A^{-k}x)]^{-M} b^{-k},
$$
which completes the proof of Lemma \ref{l5.1}.
\end{proof}

By Lemma \ref{l5.1} and an argument similar to the proof of  \cite[Lemma 2.2]{cyz}, we have the following estimates. We leave the details to the reader.

\begin{lem}\label{l5.2}
Let $i=1,\ 2$, $A_i$ be a dilation on $\rr^{n_i}$, $s_i\in\zz_+$ and
$M_i\in [0,\ \fz)$. Suppose that $f\in\cs_{s_1,\, s_2}(\rnm)$,
$\vz^{(1)}\in\cs_{s_1}(\rn),\, \vz^{(2)}\in\cs_{s_2}(\rn)$ and
$\vz_{t_1,\, t_2}(x)\equiv\vz^{(1)}_{t_1}(x_1)\vz^{(2)}_{t_2}(x_2)$
for all $t_1,\, t_2\in\zz$ and $x=(x_1,\,x_2)\in\rnm$. Then there
exists $C>0$ such that
$|(\vz_{t_1,\ t_2}\ast f)(x)| / C$
is bounded for all $x\in\rnm$ by:
\begin{align*}
&
\prod_{i=1}^2b_i^{t_i(s_i+1)\zeta_{i,-}}[1+\rho_i(x_i)]^{-M_i}
\qquad\text{if }\ t_1,\ t_2 \le0,\\
& b_1^{t_1(s_1+1)\zeta_{1,-}}
b_2^{-t_2[(s_2+1)\zeta_{2,-}+1]}[1+\rho_1(x_1)]^{-M_1}
[1+b_2^{-t_2}\rho_2(x_2)]^{-M_2}
\qquad\text{if }\ t_1\le 0,\ t_2 \ge 0,\\
& b_1^{-t_1[(s_1+1)\zeta_{1,-}+1]}
b_2^{t_2(s_2+1)\zeta_{2,-}} [1+b_1^{-t_1}\rho_1(x_1)]^{-M_1}
[1+\rho_2(x_2)]^{-M_2}
\qquad\text{if }\ t_1\ge 0,\ t_2 \le 0,\\
& \prod_{i=1}^2
[1+b_i^{-t_i}\rho_i(x_i)]^{-M_i}b_i^{-t_i[(s_i+1)\zeta_{i,-}+1]}
\qquad\text{if }\ t_1,\ t_2\ge 0.
\end{align*}
\end{lem}

We now turn to the proof of Theorem \ref{t5.1}.

\begin{proof}[Proof of Theorem \ref{t5.1}]

We first show (i). Let the notation be as in the proof of the
Lemma \ref{l4.1}. For $f\in\hp$, by $\eqref{e4.8}$, we have
\begin{eqnarray}\label{e5.3}
\hs f=\sum_{k\in\zz}\lz_k a_k=\sum_{k\in\zz}\lz_k\sum_{P\in
m(\wz\Oz_k)} a_P= \sum_{k\in\zz}\lz_k\sum_{P\in
m(\wz\Oz_k)}\sum_{R\in\hr_k,\ R^\ast=P}\lz_k^{-1}e_R
\end{eqnarray}
in $\cs'(\rnm)$. For $N,\,L\in\nn$  and $k\in\zz$, let $\hr_{k,
L}\equiv \{R\in\hr_k:\, |\ell(R_i)|\le L,\ i=1,\ 2\} $ and
$f_{N,\,L}\equiv\sum_{|k|\le N}\lz_k a_{k,\,L},$
where $a_{k,\, L}\equiv\sum_{P\in m(\wz\Oz_k)} a_{P,\, L}$, $a_{P,\,
L}\equiv \sum_{R\in \hr_{k,L},\, R^\ast=P}\lz_k^{-1}e_R$ if
$\{R\in\hr_{k,\,L}:\, R^\ast=P\}\ne \emptyset$ and otherwise
$a_{P,\,L}=0$.

On the other hand, notice that $\Oz_k$ is a bounded set. In fact,
let $M_i>0$ satisfying that $(s_i+1)\zeta_{i,-}-M_i>0$. Observing
that $1+\rho_i(x_i)\le b_i^{t_i}+\rho_i(x_i) \sim
b_i^{t_i}+\rho_i(y_i)$ for $y_i\in x_i+B^{(i)}_{t_i}$ and
$t_i\in\zz_+$, by Lemma \ref{l5.2}, we have
\begin{eqnarray*}
&&[\vec S_\psi(f)(x_1,\,x_2)]^2\\
&&\quad\ls\lf\{\int_{-\fz}^0\int_{-\fz}^0 \int_{y_1\in
x_1+B^{(1)}_{t_1}}\int_{y_2\in x_2+B^{(2)}_{t_2}}
 \frac{b_1^{2t_1(s_1+1)\zeta_{1,-}}b_2^{2t_2(s_2+1)\zeta_{2,-}}}
 {[1+\rho_1(y_1)]^{2M_1}[1+|\rho_2(y_2)]^{2M_2}}\r.\\
&&\quad\quad+\int_0^\fz\int_{-\fz}^0 \int_{y_1\in
x_1+B^{(1)}_{t_1}}\int_{y_2\in x_2+B^{(2)}_{t_2}}
 \frac{b_1^{-2t_1(s_1+1)\zeta_{1,-}-2t_1}
 b_2^{2t_2(s_2+1)\zeta_{2,-}}}{[1+b_1^{-t_1}\rho_1(y_1)]^{2M_1}
 [1+|\rho_2(y_2)]^{2M_2}}\\
&&\quad\quad+\int_{-\fz}^0\int_0^\fz \int_{y_1\in
x_1+B^{(1)}_{t_1}}\int_{y_2\in x_2+B^{(2)}_{t_2}}
\frac{b_1^{2t_1(s_1+1)\zeta_{1,-}}b_2^{-2t_2(s_2+1)\zeta_{2,-}-2t_2}}
{[1+\rho_1(y_1)]^{2M_1}
[1+b_2^{-t_2}\rho_2(y_2)]^{2M_2}}\\
&&\quad\quad\lf.+\int_0^\fz\int_0^\fz \int_{y_1\in
x_1+B^{(1)}_{t_1}}\int_{y_2\in x_2+B^{(2)}_{t_2}} \frac{
b_1^{-2t_1(s_1+1)\zeta_{1,-}-2t_1}b_2^{-2t_2(s_2+1)
\zeta_{2,-}-2t_2}}{[1+b_1^{-t_1}\rho_1(y_1)]^{2M_1}
[1+b_2^{-t_2}\rho_2(y_2)]^{2M_2}}\r\}\\
&&\quad\quad\times \,dy_1\,dy_2\,
\frac{d\sz(t_1)}{b_1^{t_1}}\,\frac{d\sz(t_2)}{b_2^{t_2}}
\ls [1+\rho_1(x_1)]^{-2M_1}[1+\rho_2(x_2)]^{-2M_2}.
\end{eqnarray*}
Thus for any $k\in\zz$, $\Oz_k$ is a bounded set in $\rnm$ and so is
$\wz\Oz_k$.

Therefore, for any $N\in\nn$ and $k=-N,\,\ldots,\,N$, $a_{k,\,L}$ is
a $(p,\,q,\,\vec s)^\ast_w$-atom associated with the bounded open
set $\Oz_k$ and thus $f_{N,\,L}\in\hpaf$.

Observe that for any $\ez>0$, there exists an integer $N_\ez>0$
such that
$(\sum_{|k|>N_\ez}|\lz_k|^p)^{1/p}<\ez.$
Moreover, for $k=-N_\ez,\,\cdots,\,N_\ez$, similarly to the estimate
for \eqref{e4.11}, we have
\begin{eqnarray*}
\lf\|a_k-a_{k,\,L}\r\|_{L^q_w(\rnm)}
&&\ls \lz_k^{-1}\bigg\|\bigg(\sum_{P\in  m(\wz\Oz_k)} \sum_{\gfz{R\notin
\hr_{k,\,L}}{R^\ast=P}}c_R^2\chi_{R\cap(\Oz_k\setminus\Oz_{k+1})}\bigg)^{1/2}
\bigg\|_{L^q_w(\rnm)},\\
\end{eqnarray*}
which together with \eqref{e4.12} implies that $\lf\|a_k-a_{k,\,
L}\r\|_{L^q_w(\rnm)}\to0$ as $L\to\fz$. Similarly to the estimate of
\eqref{e4.14}, we also have $\sum_{P\in m(\wz\Oz_k)}
\lf\|a_P-a_{P,\, L} \r\|^q_{L^q_w(\rnm)}\to 0$ as $L\to\fz$. Thus
there exists an integer $L_\ez>0$ such that
$\frac{(2(N_\ez+1))^{1/p}}{\ez}(a_k-a_{k,\, L_\ez})$ is a
$(p,\,q,\,s)_w$-atom. Therefore,
\begin{eqnarray*}
&& \|f-f_{N_\ez,\,L_\ez}\|_\hp\\
&&\hs \ls
\bigg\{\sum_{|k|>N_\ez}|\lz_k|^p\bigg\}^{1/p}+
\bigg\|\sum_{|k|\le N_\ez}\lz_k(a_k-a_{k,\, L_\ez})\bigg\|_\hp\\
&&\hs\ls\ez+\|f\|_{\hp}\bigg[\sum_{|k|\le N_\ez}\|a_k-a_{k,\,
L_\ez}\|_\hp^p\bigg]^{1/p}
\ls\ez (1+\|f\|_{\hp}),
\end{eqnarray*}
 which gives (i).

Now we prove (ii). From Definition \ref{d5.1} and Theorem \ref{t4.1},
we automatically deduce
$\|f\|_\hp\ls \|f\|_\hpaf.$
Thus, to show (ii), it suffices to prove that for all $f\in\hpaf$,
$\|f\|_\hpaf\ls \|f\|_\hp.$

Let $f\in\hpaf$. Since $f\in \cs(\rnm)$, by Lemma \ref{l2.5}, we
know that \eqref{e5.3} also holds in $L^q(\rnm)$ and hence,
pointwise. Assume that $\supp f\subset B^{(1)}_{h_1}\times
B^{(2)}_{h_2}$ for certain $h_1,\ h_2\in\zz$. By homogeneity, we
further assume that $\|f\|_\hp\equiv\|\vec
S_\psi(f)\|_{L^p_w(\rnm)}=1$, where $\psi$ is as in Proposition
\ref{p2.5}.

Let $i=1,\,2$. For certain given $N\in\nn$ which will be determined
later, set $D_i\equiv-v_iN+u_i+\sz_i$; then we choose certain
$M_0\in\nn$, depending on $N$, such that $d_i\equiv
v_i(M_0N-1)+u_i+\sz_i$ satisfies $d_i(s_i+1)\zeta_{i,\,-}\le
-D_i[1+(s_i+1)\zeta_{i,\,-}].$  We first assume that $N$ is large
enough such that $D_i>h_i$. Then, by the definition of $R_+$ in
\eqref{e4.1}, we know that there exist finite dyadic rectangles $R$,
whose collection is denoted by $\hr^{N}$, such that
\begin{eqnarray}\label{e5.4}
R_+\cap\lf\{ B^{(1)}_{D_1+\sz_1}\times
B^{(2)}_{D_2+\sz_2}\times[d_1,\ D_1)\times[d_2,\
D_2)\r\}\not=\emptyset.
\end{eqnarray}
From now on, we adopt the notation in the proof of Lemma \ref{l4.1}
again. Observe that for each $R\in\hr^{N}$, there exists a unique
$k\in\zz$ such that $R\in\hr_k$, and we denote by $J_{N}$ the set of
all such $k$'s.

Let $\wz
a_{P,\,N}\equiv\lz_k^{-1}\sum_{R\in\hr_k\cap\hr^N,\,R^\ast=P}e_R $
if $\{R\in\hr_k\cap\hr^N:\, R^\ast=P\}\ne\emptyset$ and otherwise
$\wz a_{P,\,N}=0$. Let $\wz a_{k,N}\equiv \sum_{P\in m(\wz\Oz_k)}\wz
a_{P,\,N}$. Then similarly to the proof of Lemma \ref{l4.1}, we know
that $\wz a_{k,\,N}$ is a $(p,\,q,\,\vec s)_w$-atom which is a
finite linear combination of particles $\wz a_{P,\,N}$. Obviously,
$\wz a_{P,\,N}$ is also a finite linear combination of $e_R$ and
hence is smooth. This further implies that $\wz a_{k,\,N}$ is a
$(p,\,q,\,\vec s)^\ast_w$-atom. Let $f_{N}\equiv\sum_{k\in J_N}\lz_k
\wz a_{k,\, N}$ and  $g_{N}\equiv f-f_{N}$. Then $f_N\in\hpaf$ and
$\|f_N\|^p_\hpaf\le\sum_{k\in\zz}|\lz_k|^p\ls1.$

So it remains to prove $g_{N}\in \hpaf$ and $\|g_N\|^p_\hpaf\ls1.$
In fact, we will prove that there exists a positive constant $\wz C$,
independent of $f$ and $N$, such that $\wz Cg_{N} $ is a
$(p,\,q,\,s)^\ast_w$-atom, which implies $\|g_N\|^p_\hpaf\ls1.$

Obviously, $g_N\in \cs_{s_1,\, s_2}(\rnm)$. Noticing that if
$R\in\hr^N$, then by \eqref{e5.4}, $\ell(R_i)\in(-N,\ M_0N)$. By
this, \eqref{e2.1} and Lemma \ref{l2.1}(iv), we further obtain
$$R_i'\equiv x_{R_i}+B^{(i)}_{v_i(\ell(R_i)-1)+u_i+2\sz_i}\subset x_{R_i}
+B^{(i)}_{D_i+\sz_i},$$
which together with  $R\subset R'$  and
$\eqref{e5.4}$ yields that $(x_{R_i}+B^{(i)}_{D_i+\sz_i})\cap
B^{(i)}_{D_i+\sz_i}\not=\emptyset$. Then by \eqref{e4.9} and
\eqref{e2.1}, we obtain
\begin{eqnarray*}
\supp f_N\subset \bigcup_{R\in \hr^N} R'\subset\lf(
B^{(1)}_{D_1+3\sz_1}\times B^{(2)}_{D_2+3\sz_2}\r).
\end{eqnarray*}
From this, $\supp f\subset B^{(1)}_{h_1}\times B^{(2)}_{h_2}$ and
$D_i>h_i$, it follows that $\supp g_N\subset
(B^{(1)}_{D_1+3\sz_1}\times B^{(2)}_{D_2+3\sz_2})$.

We now claim that there exists an $N_0\in\nn$, depending on
$f,\,w,\,m,\,n,\,A_1$ and $A_2$, such that for all $N\ge N_0$,
\begin{eqnarray}\label{e5.5}
\|g_N\|_{L^q_w(\rnm)}\le \lf[w(B^{(1)}_{D_1}\times
B^{(2)}_{D_2})\r]^{1/q-1/p},
\end{eqnarray}
Now we prove that there exists a positive constant $\wz C$,
independent of $f$ and $N$, such that $\wz Cg_{N} $ is a
$(p,\,q,\,\vec s)^\ast_w$-atom.

In fact, by Lemma \ref{l2.1}(i), there exist certain
$P_i\in\cq^{(i)}$ and $x_{i,\,0}\in\rr^{n_i}$ satisfying that
$x_{i,\, 0}\in P_i\cap B^{(i)}_{D_i+\sz_i}$ and $v_i\ell(P_i)+u_i<
D_i+\sz_i\le v_i[\ell(P_i)-1]+u_i$. For this $P$, let $P''$ be as in
Definition \ref{d4.2}(I). Then $P\equiv P_1\times P_2 \subset
B^{(1)}_{D_1+3\sz_1}\times B^{(2)}_{D_2+3\sz_2}\subset P''$. To see
this, for any $x_i\in P_i$, since $x_{i,\,0}\in P_i\cap
B^{(i)}_{D_i+\sz_i}$ and $v_i\ell(P_i)+u_i< D_i+\sz_i$, using Lemma
\ref{l2.1}(iv) and \eqref{e2.1}, we obtain
\begin{eqnarray*}
&x_i\in x_{P_i}+B^{(i)}_{v_i\ell(P_i)+u_i}&\subset
x_{i,\,0}+B^{(i)}_{v_i\ell(P_i)+u_i}+B^{(i)}_{v_i\ell(P_i)+u_i}\\
&&\subset B^{(i)}_{D_i+\sz_i}+B^{(i)}_{D_i+\sz_i}+
B^{(i)}_{D_i+\sz_i}+\subset B^{(i)}_{D_i+3\sz_i},
\end{eqnarray*}
which implies that $P\subset B^{(1)}_{D_1+3\sz_1}\times
B^{(2)}_{D_2+3\sz_2}$. For any $x_i\in B^{(i)}_{D_i+3\sz_i}$, since
$D_i+\sz_i\le v_i[\ell(P_i)-1]+u_i$ and $x_{i,\,0}\in P_i\cap
B^{(i)}_{D_i+\sz_i}$, by Lemma \ref{l2.1}(iv) and \eqref{e2.1}, we
have
\begin{eqnarray*}
&x_i-x_{P_i}\in B^{(i)}_{D_i+3\sz_i}+x_{i,\,0}
+B^{(i)}_{v_i\ell(P_i)+u_i}&\subset
B^{(i)}_{v_i[\ell(P_i)-1]+u_i+2\sz_i}
+B^{(i)}_{v_i[\ell(P_i)-1]+u_i+\sz_i}\\
&&\subset B^{(i)}_{v_i[\ell(P_i)-1]+u_i+3\sz_i},\nonumber
\end{eqnarray*}
which implies that $B^{(1)}_{D_1+3\sz_1}\times
B^{(2)}_{D_2+3\sz_2}\subset P''$.

Let $\Oz\equiv B^{(1)}_{D_1+3\sz_1}\times B^{(2)}_{D_2+3\sz_2} $ and
$\wz\Oz$ be as in \eqref{e4.3}. Obviously, $\Oz$ is an open bounded
set. Noticing that $P\subset \Oz$, then we have $P\subset \wz \Oz$.
Thus, there exists a dyadic rectangle $P^\star\in m(\wz\Oz)$ such
that $P\subset P^\star$. Moreover, since $P\subset P^\star$,
similarly to the proof of $\Oz\subset P''$, we have that
 $\Oz\subset (P^\star)''$.
For $R\in\ m(\wz\Oz)$, let $a_R=g_N$ if $R=P^\star$ and $a_R=0$ if
$R\in m(\wz\Oz)$ and $R\ne P^\star$. By the vanishing moment
satisfied by $g_N$ and \eqref{e5.5} together with proposition
\ref{p2.2}(i), we know that $\wz Cg_N$ is a
$(p,\,q,\,s)_w^\ast$-atom associated with $\Oz$ for certain
positive constant $\wz C$ independent of $f$ and $N$.

Finally, we establish the estimate \eqref{e5.5}. Since $f\in
\cs(\rnm)$, by \eqref{e4.4}, \eqref{e5.3} and
$$
\lf|\lf\{ B^{(1)}_{D_1+\sz_1}\times B^{(2)}_{D_2+\sz_2}\times[d_1,\
D_1)\times[d_2,\ D_2)\r\}\setminus\lf(\bigcup_{R\in
\hr^N}R_+\r)\r|=0,$$together with the observation that for two
different rectangles $R$ and $S$, then
$R_+\cap S_+=\emptyset$, we have that for all $x\in
B^{(1)}_{D_1+3\sz_1}\times B^{(2)}_{D_2+3\sz_2}$,
\begin{eqnarray*}
|g_N(x)|&& =\lf|\sum_{R\in
\hr}e_R(x)-\sum_{R\in\hr^N}e_R(x)\r|\ls \bigg[ \int_{[d_1,\ D_1)^{\complement}\times \rr}
\int_{\rnm}
+\int_{\rr\times[d_2,\ D_2)^{\complement} }
\int_{\rnm}
\\
&&\hs +\int_{[d_1,\ D_1)\times[d_2,\ D_2)}\int_{ \rn\times
{\lf(B^{(2)}_{D_2+\sz_2}\r)}^{\complement}}
+\int_{[d_1,\ D_1)\times[d_2,\ D_2)}\int_{
{\lf(B^{(1)}_{D_1+\sz_1}\r)}^{\complement}\times\rrm}
\bigg]
\\
&&\hs\times|\tz_{t_1,\ t_2}(x-y)(\psi_{t_1,\ t_2}\ast f)(y)| \,
dy\,d\sz(t_1)\,d\sz(t_2) \equiv {\mathrm J_1}+{\mathrm J_2}+{\mathrm
J_3}+{\mathrm J_4}.
\end{eqnarray*}

Recall that $\tz=\tz^{(1)}\tz^{(2)}\in\cs_{s_1,\,s_2}(\rnm)$ with
$\supp\tz^{(1)}\subset B^{(1)}_0$ and $\supp\tz^{(2)}\subset
B^{(2)}_0$.
 Notice that if $x_i-y_i\in B^{(i)}_{t_i}$, then for all $t_i\in\zz$,
 we have
 \begin{eqnarray}\label{e5.6}
 1+b_i^{-t_i}\rho(x_i)\sim 1+b_i^{-t_i}\rho(y_i).
\end{eqnarray}
Let $M_i> 1$ for $i=1,\,2$. Since
$f,\,\psi=\psi^{(1)}\psi^{(2)},\,\tz=\tz^{(1)}\tz^{(2)}
\in\cs_{s_1,\,s_2}(\rnm)$, by Lemma \ref{l5.2} and \eqref{e5.6} , we
have
\begin{eqnarray*}
{\mathrm J_2}&&\ls\bigg(\int_0^\fz
\frac{b_1^{-t_1[(s_1+1)\zeta_{1,-}+1]}}{[1+b_1^{-t_1}
\rho_1(x_1)]^{M_1}}\,d\sz(t_1) +
\int_{-\fz}^0\frac{b_1^{t_1(s_1+1)\zeta_{1,\,-}}}
{[1+\rho_1(x_1)]^{M_1}}
\,d\sz(t_1)\bigg) \\
&&\hs\quad\times \bigg(\int_{-\fz}^{d_2}
\frac{b_2^{t_2(s_2+1)\zeta_{2,-}}}{[1+\rho_2(x_2)]^{M_2}}\,d\sz(t_2)
+\int_{D_2}^\fz
\frac{b_2^{-t_2[(s_2+1)\zeta_{2,-}+1]}}{[1+b_2^{-t_2}
\rho_2(x_2)]^{M_2}}\,d\sz(t_2)\bigg)\\
&&\ls b_2^{d_2(s_2+1)\zeta_{2,-}}
+ b_2^{-D_2 [(s_2+1)\zeta_{2,-}+1]}
\ls b_2^{ -D_2[1+(s_2+1)\zeta_{2,-}]}.
\end{eqnarray*}
The last inequality is a consequence of our stipulation that
$d_2(s_2+1)\zeta_{2,\,-}\le - D_2[1+(s_2+1)\zeta_{2,-}].$
Moreover, by the assumptions \eqref{e5.1}, \eqref{e5.2}
and that $(p,\,q,\,\vec s)_w$ is an admissible triplet,
there exists $\kappa>0$ such that
$(s_i+1)\zeta_{i,\,-}+1-(q_w+\kappa)/p>0$ for $i=1,\,2$, and
\begin{align}
\label{e5.7} b^\ast_1 & \equiv
b_1^{v_1[(s_1+1)\zeta_{1,\,-}+1-(q_w+\kappa)/p]}
b_2^{-v_2(q_w+\kappa)/p}<1,
\\
\label{e5.8} b^\ast_2 & \equiv b_1^{-v_1(q_w+\kappa)/p}
b_2^{v_2[(s_2+1)\zeta_{2,\,-}+1-(q_w+\kappa)/p]} <1.\end{align}
Thus, by \eqref{e5.8}, $\supp g_N\subset B^{(1)}_{D_1+3\sz_1}\times
B^{(2)}_{D_2+3\sz_2}$ and Proposition \ref{p2.2}(i) with
$w\in\ca_{q_w+\kappa}(\vec A)$, if we choose $N$ large enough, we
further obtain
\begin{eqnarray*}
\|{\mathrm J_2}\|_{L^q_w(\rnm)} &&\le C \lf[w\lf(B_{D_1}^{(1)}\times
B_{D_2}^{(2)}\r)\r]^{1/q}
b_2^{ - D_2[1+(s_2+1)\zeta_{2,-}]} \\
&&\le C\lf[w\lf(B_{D_1}^{(1)}\times B_{D_2}^{(2)}\r)\r]^{1/q-1/p}
(b_2^\ast)^N\le \lf[w\lf(B_{D_1}^{(1)}\times B_{D_2}^{(2)}\r)\r]^{1/q-1/p},
\end{eqnarray*}
where $C$ is a positive constant, which is the desired estimate.

For ${\mathrm J_4}$,
 observe that if $y_1\in (B^{(1)}_{D_1+\sz})^{\complement},\,
t_1\le D_1$ and $\tz^{(1)}_{t_1}(x_1-y_1)\not=0$, then by
\eqref{e2.2}, we have
\begin{eqnarray}\label{e5.9}
x_1\in y_1+B^{(1)}_{t_1}\subset \lf(B^{(1)}_{D_1}\r)^{\complement}
\end{eqnarray}
and thus $\rho_1(x_1)\ge b_1^{ D_1}$. Let
$M_1\in(1,\,(s_1+1)\zeta_{1,-}+1)$ and $M_2>1$. Then by \eqref{e5.9}
and an argument similar to the estimate of ${\mathrm J_2}$, we have
\begin{eqnarray*}
{\mathrm J_4}&&\ls\bigg[\int_0^{D_1}
\frac{b_1^{-t_1[(s_1+1)\zeta_{1,-}+1]}}{[1+b_1^{-t_1}
\rho_1(x_1)]^{M_1}}\,d\sz(t_1) +
\int_{d_1}^0\frac{b_1^{t_1(s_1+1)\zeta_{1,\,-}}}
{[1+\rho_1(x_1)]^{M_1}}
\,d\sz(t_1)\bigg]\\
&&\hs\times \bigg[\int_{d_2}^{0}
\frac{b_2^{t_2(s_2+1)\zeta_{2,-}}}{[1+\rho_2(x_2)]^{M_2}}\,d\sz(t_2)
+\int_0^{D_2}
\frac{b_2^{-t_2[(s_2+1)\zeta_{2,-}+1]}}{[1+b_2^{-t_2}
\rho_2(x_2)]^{M_2}}\,d\sz(t_2)\bigg]
\ls b_1^{ -D_1[(s_1+1)\zeta_{1,-}+1]}.
\end{eqnarray*}
Moreover, by \eqref{e5.7} and Proposition \eqref{p2.2}(i) with
$w\in\ca_{q_w+\kappa}(\vec A)$, similarly to the estimate of
$\|{\mathrm J_2}\|_{L^q_w}(\rnm)$, if we choose $N$ large enough, we
then have
\begin{eqnarray*}
\|{\mathrm J_4}\|_{L^q_w(\rnm)}\le \lf[w\lf(B_{D_1}^{(1)}\times
B_{D_2}^{(2)}\r)\r]^{1/q-1/p}.
\end{eqnarray*}

By symmetry, we have similar estimates for $\|{\mathrm
J_1}\|_{L^q_w(\rnm)}$ and $\|{\mathrm J_3}\|_{L^q_w(\rnm)}$, which
gives \eqref{e5.5} and hence completes the proof of Theorem
\ref{t5.1}.
\end{proof}

\section{Applications\label{s6}}

\hskip\parindent We first recall that a \textit{quasi-Banach space}
$\cb$ is a vector space endowed with a quasi-norm $\|\cdot\|_\cb$
which is non-negative, non-degenerate (i.\,e., $\|f\|_\cb=0$ if and
only if $f=0$), homogeneous, and obeys the quasi-triangle
inequality, i.\,e., there exists a positive constant $K$ no less
than $1$ such that for all $f,\, g\in\cb$, $\|f+g\|_\cb\le
K(\|f\|_\cb+\|g\|_\cb)$.

Recall that the following notion of $\gz$-quasi-Banach spaces was
first introduced in \cite{yz}.

\begin{defn}\label{d6.1}
Let $\gz\in(0,\,1]$. A quasi-Banach space $\cbgz$ with the
quasi-norm $\|\cdot\|_\cbgz$ is called a \textit{$\gz$-quasi-Banach
space} if $\|f+g\|^\gz_\cbgz\le \|f\|^\gz_\cbgz+\|g\|^\gz_\cbgz$ for
all $f,\, g\in\cbgz$.
\end{defn}

Notice that any Banach space is a $1$-quasi-Banach space, and the
quasi-Banach spaces $\ell^\gz$, $L_w^\gz(\rnm)$ and
$H_w^\gz(\rnm;\,\vec A)$ with $\gz\in(0,\,1)$ are typical
$\gz$-quasi-Banach spaces. Moreover,  according to the Aoki-Rolewicz
theorem (see \cite{a}, \cite[p. 66]{g1} or \cite{r84}), any
quasi-Banach space is essentially a $\gz$-quasi-Banach space,
where $\gz\equiv[\log_2(2K)]^{-1}$.

For any given $\gz$-quasi-Banach space $\cbgz$ with $\gz\in(0,\,1]$
and a linear space ${\mathcal Y}$, an operator $T$ from ${\mathcal
Y}$ to $\cbgz$ is called \textit{$\cbgz$-sublinear} if for all $f,\,
g\in{\mathcal Y}$ and $\lz,\,\nu\in\cc$, we have
$$\|T(\lz f+\nu g)\|_\cbgz
\le\lf(|\lz|^\gz\|T(f)\|^\gz_\cbgz+
|\nu|^\gz\|T(g)\|^\gz_\cbgz\r)^{1/\gz}$$ and $\|T(f)-T(g)\|_\cbgz\le
\|T(f-g)\|_\cbgz.$ The notion of $B_\gz$-sublinear operators was
first introduced in \cite{yz08}.

We remark that if $T$ is linear, then $T$ is $\cbgz$-sublinear.
Moreover, if $\cbgz$ is a space of functions, $T$ is sublinear in
the classical sense and $T(f)\ge 0$ for all $f\in{\mathcal Y}$, then
$T$ is also $\cbgz$-sublinear.

\begin{thm}\label{t6.1}
Let $w\in \ca_\fz(\vec A)$, $q_w$ as in \eqref{e2.7} and
$(p,\,q,\,\vec s)_w$ an admissible triplet.
Let $\gz\in [p,\,1]$ and $\cbgz$ be a
$\gz$-quasi-Banach space. Suppose that $T:\hpaf\to \cbgz$ is a
$\cbgz$-sublinear operator such that
\begin{eqnarray}\label{e6.1}
\sup\{ \|T(a)\|_\cbgz:\, a \ {\rm is\ any}\ (p, \,q, \,\vec
s)^\ast_w{\rm -atom}\}<\fz.
\end{eqnarray}
Then there exists a unique bounded $\cbgz$-sublinear operator $\wz
T$ from $\hp$ to $\cbgz$ which extends $T$.
\end{thm}

\begin{proof} Without loss of generality, we may also
assume that $\vec s$ satisfies \eqref{e5.1} and \eqref{e5.2}. For
every $f\in\hpaf$, by Theorem \ref{t5.1}(ii), there exist
$\{\lz_j\}_{j=1}^\ell\subset \cc$ and $(p, \,q, \,\vec
s)^\ast_w$-atoms $\{a_j\}_{j=1}^\ell$ such that
$f=\sum_{j=1}^\ell\lz_ja_j$ pointwise and
$\sum_{j=1}^\ell|\lz_j|^p\ls\|f\|^p_{\hp}.$ Then by
\eqref{e6.1}, we have
$$\|T(f)\|_{\cbgz}\ls\lf[ \sum_{j=1}^\ell
|\lz_j|^p \r]^{1/p}\ls\|f\|_\hp.$$
Since $\hpaf$ is dense in $\hp$
by Theorem \ref{t5.1}(i), a density argument gives the desired
result. This finishes the proof of Theorem \ref{t6.1}.
\end{proof}

\begin{rem}\label{r6.1}
If $T$ is a bounded $\cbgz$-sublinear operator from $\hp$ to
$\cbgz$, then it is clear that for all admissible triplet
$(p,\,q,\,\vec s)_w$, $T$ maps all $(p,\,q,\,\vec s)^\ast_w$-atoms
into uniformly bounded elements of $\cbgz$. Thus the condition
\eqref{e6.1} of Theorem \ref{t6.1} is also necessary.
\end{rem}

Motivated by Theorem 1 in \cite{f87}, we introduce the rectangular
atoms in the current setting and then derive the boundedness of
sublinear operators from their behavior on rectangular atoms.

\begin{defn}\label{d6.2}
Let $w\in \ca_\fz(\vec A)$ and $q_w$ be as in \eqref{e2.7} and $(p,
\ q,\ \vec s)_w$ be an admissible triplet as in Definition
\ref{d4.2}. For $R\in\hr$, a function $a_R$ is said to be a {\it
rectangular $(p,\,q,\,\vec s)_w$-atom} if

(i) $a_R$ is supported on $R''=R''_1\times R''_2$, where
$R''_i\equiv x_{R_i}+B^{(i)}_{v_i(\ell(R_i)-1)+u_i+3\sz_i},\ i=1,\
2$;

(ii) $\int_\rrm a_R(x_1,\ x_2)x_1^\az\,dx_1=0$ for all $|\az|\le
s_1$ and almost all $x_2\in \rrm$, and

\quad\ \  $\int_\rn a_R(x_1,\, x_2)x_2^\bz\,dx_2=0$ for all
$|\bz|\le s_2$ and almost all $ x_1\in \rn$;

(iii) $\|a\|_{L^q_w(\rnm)}\ls [w(R)]^{1/q-1/p}$.
\end{defn}

Let $i=1,\,2$. For any $R_i\in\cq^{(i)}$ and $k\in\zz_+$, set
$R_{i,\,k}\equiv x_{R_i}+B^{(i)}_{v_i(\ell(R_i)-1)+u_i+5\sz_i+k}.$

The following corollary is very useful in the study of boundedness
of  operators in $\hp$.

\begin{cor}\label{c6.1}
Let $w\in \ca_\fz(\vec A)$, $q_w$ as in \eqref{e2.7} and
$(p,\,q_1,\,\vec s)_w$ an admissible triplet.
Let  $T$ be a
bounded sublinear operator from $L^{q_1}_w(\rnm)$ to
$L^{q_0}_w(\rnm)$, where $q_0\in[q_1,\, \fz)$. Let $q\in [p,\, 2)$
be such that $1/q-1/p=1/q_0-1/q_1$. If there exist positive
constants $C,\, \ez$ such that for all $k\in\zz_+$ and all
rectangular $(p,\, q_1,\,\vec s)_w$-atoms $a_R$,
\begin{eqnarray}\label{e6.2}
\int_{(R_{1,k}\times R_{2,\,k})^\complement}|T(a_R)(x)|^q
w(x)\,dx\le C \min \{ b_1^{-k\ez},\, b_2^{-k\ez}\},
\end{eqnarray}
then  $T$ uniquely extends to a bounded operator from $\hp$ to
$L^q_w(\rnm)$.
\end{cor}

\begin{proof}
Let all the notation be as in the proof of Lemma \ref{l4.3}. To show
Corollary \ref{c6.1}, by Theorem \ref{t6.1}, we only need to show
that for all $(p,\,q_1,\,\vec s)^\ast_w$-atoms $a=\sum_{R\in
m(\wz\Oz)} a_R$, $\|T a\|_{L^q_w(\rnm)}\ls 1$.

Recall that $\eta_0\equiv b_1^{v_1-5\sz_1}b_2^{v_2-5\sz_2}$. For
each $R=R_1\times R_2\in m(\wz\Oz)$, let $\wh R_1\equiv\wh R_1(R_2)$
being the ``longest'' dyadic cube containing $R_1$ such that
$|(\wh R_1\times R_2)\cap\wz\Oz|>\eta_0|\wh R_1\times R_2|.$
Let
\begin{eqnarray*}
\wz\Oz'\equiv \{x\in\rnm:\,
\cm_s(\chi_{\wz\Oz})(x)>b^{-2u_1}_1b^{-2u_2}_2\eta_0\}.
\end{eqnarray*}
For any given $\wz R_1\times R_2\in m(\wz\Oz')$ and $\wz R_1\supset
R_1$, let $\wh R_2\equiv\wh R_2(\wz R_1)$ being the ``longest''
dyadic cube containing $R_2$ such that
$|(\wz R_1\times \wh R_2)\cap \wz\Oz'|>\eta_0|\wz R_1\times \wh R_2|.$

Let $\gz_1(R)\equiv\gz_1(R,\,\wz\Oz)\equiv\ell(\wh R_1)-\ell(R_1)$
and $\gz_2(\wz R_1\times R_2)\equiv\gz_2(\wz R_1\times
R_2,\,\wz\Oz')\equiv\ell(\wh R_2)-\ell(R_2)$. Then by Lemma
\ref{l4.4}, for any $\dz>0$, we obtain
\begin{eqnarray}\label{e6.3}
\sum_{R\in m(\wz\Oz)}b_1^{\gz_1(R)\dz}w(R)\ls w(\Oz)
\end{eqnarray}
and
\begin{eqnarray}\label{e6.4}
\sum_{R\in m_1(\wz\Oz')}b_2^{\gz_2(R_1\times R_2)\dz}w(R)\ls w(\Oz).
\end{eqnarray}

Set $$\bar R^\ast\equiv \bar R_1^\ast\times \bar R_2^\ast\equiv
(x_{R_1}+B^{(1)}_{v_1(\ell(\wh R_1)-1)+u_1+5\sz_1})\times
(x_{R_2}+B^{(2)}_{v_2(\ell(\wh R_2)-1)+u_2+5\sz_2}).$$

By the argument for \eqref{e4.40} and the
$L^{q_0}_w(\rnm)$-boundedness of $\cm_s$ (see Proposition
\ref{p2.2}(ii)), we have
\begin{eqnarray}\label{e6.5}
\bigcup_{R\in m(\wz\Oz)} \bar R^\ast\subset \wz\Oz'''\ {\rm and} \
w(\wz\Oz''')\ls w(\Oz).
\end{eqnarray}
From this, $1/q-1/p=1/q_0-1/q_1$ and the size condition of $a$,
together with H\"older's inequality and the boundedness of $T$ from
$L^{q_1}_w(\rnm)$ to $L^{q_0}_w(\rnm)$, we deduce that
\begin{eqnarray*}
&\bigg\{\dint_{\wz\Oz'''} |T(a)(x)|^q w(x)\, dx\bigg\}^{1/q}&\ls
\bigg\{\int_{\wz\Oz'''} |T(a)(x)|^{q_0} w(x)\,
dx\bigg\}^{1/q_0}w(\wz\Oz''')^{1/q-1/q_0}\\
&&\ls \|a\|_{L^{q_1}_w(\rnm)}w(\Oz)^{1/p-1/q_1}\ls 1.
\end{eqnarray*}

It remains to prove that $\int_{ \lf(\wz\Oz'''\r)^\complement}
|T(a)(x)|^q\, dx\ls 1$. Without loss of generality, we may assume
that $q\le 1$. The proof of the case $q\in (1,\, 2)$ is similar and
we omit the details.
Since $q\le 1$ and $a=\sum_{R\in m(\wz\Oz)}a_R$, by \eqref{e6.5}, we
obtain
\begin{eqnarray*}
&&\int_{ \lf(\wz\Oz'''\r)^\complement} |T(a)(x)|^qw(x)\,dx\\
&&\hs\le \sum_{R\in m(\wz\Oz)}\int_{ \lf(\wz\Oz'''\r)^\complement}
|T(a_R)(x)|^qw(x)\,dx\\
&&\hs \le \sum_{R\in m(\wz\Oz)} \bigg[\int_{(\rn\setminus \bar
R_1^\ast)\times \rrm} + \int_{\rn\times (\rrm\setminus \bar
R_2^\ast)} \bigg] |T(a_R)(x)|^qw(x)\,dx \equiv {\mathrm
E_1}+{\mathrm E_2}.
\end{eqnarray*}
Since $a_R[w(R)]^{1/q_1-1/p}\|a_R\|_{L^{q_1}_w(\rnm)}^{-1}$ is a
rectangular $(p,\,q_1,\,\vec s)_w$-atom, by \eqref{e6.2}, we have
\begin{eqnarray*}
\int_{(\rrm\setminus \bar R_1^\ast)\times
\rn}|T(a_R)(x)|^qw(x)\,dx\ls
\|a_R\|^q_{L^{q_1}_w(\rnm)}[w(R)]^{1-q/q_0}b_1^{\ez\gz_1(R)}.
\end{eqnarray*}
From this, $1/q_1-1/p=1/q_0-1/q$, H\"older's inequality, the size
condition of $a$ and \eqref{e6.3}, it follows that
\begin{eqnarray*}
&{\mathrm E_1} &\ls \bigg\{\sum_{R\in
m(\wz\Oz)}\|a_R\|^{q_1}_{L^{q_1}_w(\rnm)}\bigg\}^{q/q_1} \bigg\{
\sum_{R\in m(\wz\Oz)} [w(R)]^{\frac{q_1(q_0-q)}{q_0(q_1-q)}}
b_1^{\frac{q_1\ez\gz_1(R)}{q_1-q}}\bigg\}^{1-q/q_1}
\\
&&\ls [w(\Oz)]^{q(1/q_1-1/p)}[w(\Oz)]^{q(1/q_1-1/q_0)}
\bigg\{
\sum_{R\in m(\wz\Oz)} w(R)b_1^{q_1\ez\gz_1(R)/(q_1-q)}
\bigg\}^{1-q/q_1}\ls 1.
\end{eqnarray*}
Similarly, by \eqref{e6.4}, we obtain ${\mathrm E_2}\ls 1$. This
finishes the proof of Corollary \ref{c6.1}.
\end{proof}

\appendix
\renewcommand{\thesection}{\appendixname}
\section{}

\renewcommand{\thesection}{\Alph{section}}

\hskip\parindent In this appendix, we give the proof of Proposition
\ref{p3.1} by establishing a more general version, namely,
Theorem \ref{ta.1} below. Let $\cb$ be a Banach space and
$L^{\fz}_c(\rn,\, \cb)$ the set of $f\in
L^\fz(\rn,\,\cb)$ with compact support. Through the whole appendix,
we use $\cb_1$ and $\cb_2$ to denote two {\it Banach spaces}.

\begin{defn}\label{da.1}
An operator $\ct$ is called a Calder\'on-Zygmund operator if $\ct$
is bounded from $L^r(\rn,\,\cb_1)$ to $L^r(\rn,\,\cb_2)$ for certain
fixed $r\in(1,\,\fz)$, and $\ct$ has a distributional
$\cl(\cb_1,\,\cb_2)$-valued kernel $\ck$ such that for all $f\in
L_{c}^\fz(\rn,\,\cb_1)$ and $x\not\in\supp f$,
$$\ct(f)(x)=\int_\rn\ck(x,\,y)f(y)\,dy,$$
where  $\ck$ is a standard kernel in the following sense: there
exist positive constants $C$ and $\ez$ such that for all
$x,\,y,\,z\in\rn$ satisfying $\rho(z-y)\le b^{-2\sz} \rho(x-y)$,
\begin{eqnarray}\label{ea.1}
\|\ck(x,\,y)\|_{L(\cb_1,\,\cb_2)}\le C/\rho(x-y)
\end{eqnarray}
and
\begin{eqnarray}\label{ea.2}
\quad\|\ck(y,\, x)-\ck(z,\, x)\|_{L(\cb_1,\,\cb_2)}+ \|\ck(x,\,
y)-\ck(x,\, z)\|_{L(\cb_1,\,\cb_2)}\le
C\frac{\rho(z-y)^\ez}{\rho(x-y)^{1+\ez}}.
\end{eqnarray}
\end{defn}

Let $L^{1,\,\fz}(\rn,\,\cb)$ be the set of all
$\cb$-measurable functions $f$ on $\rn$ such that
$$\|f\|_{L^{1,\,\fz}(\rn,\,\cb)}\equiv
\sup_{\az>0} \az |\lf\{x\in\rn: \|f\|_\cb>\az
\r\}|  <\fz.$$

Then by \cite[Theorem 1.1]{gly1}, we have the following result.

\begin{lem}\label{la.1}
Let $p\in(1,\, \fz)$. Suppose that $\ct$ is a Calder\'on-Zygmund
operator. Then $\ct$ is bounded from $L^p(\rn,\, \cb_1)$ to
$L^p(\rn,\, \cb_2)$ and bounded from $L^1(\rn,\, \cb_1)$ to $L^{1,\,
 \fz}(\rn,\, \cb_2)$
\end{lem}

The following theorem is the main result of this appendix,
which is a weighted version of Lemma \ref{la.1}. This theorem
extends \cite[Theorems 7.11 and 7.12]{d01} to the weighted
anisotropic settings and also has an independent interest.

\begin{thm}\label{ta.1}
 Suppose that $\ct$ is Calder\'on-Zygmund
operator. If $p\in(1,\,\fz)$ and $w\in\ca_p(A)$, then
 $\ct$ is bounded from $L^p_w(\rn,\, \cb_1)$ to
$L^p_w(\rn,\, \cb_2)$, and if $w\in\ca_1(A)$, then
 $\ct$ is bounded from
 $L^1_w(\rn,\,\cb_1)$ to $L_w^{1,\, \fz}(\rn,\, \cb_2)$.
\end{thm}

The proof of  Theorem \ref{ta.1} follows from the procedure in \cite{d01}. Here we present some details for the
convenience of readers.

To this end, we first introduce the dyadic maximal function in this
setting. For any given $\cb$-measurable function $f\in
L^1_{\loc}(\rn,\,\cb)$ and $x\in\rn$, we define the dyadic maximal
function by $M_d(f)(x)\equiv \sup_{k\in\zz} E_k(f)(x),$ where
$$E_k(f)(x)\equiv
\sum_{Q\in\cq_k}\lf(\frac{1}{|Q|}\int_Q\|f(y)\|_\cb\, dy\r)
\chi_Q(x) $$ and $\cq_k\equiv\{Q^k_\az: \az\in I_k\}$ denotes the
set of dyadic cubes as in Lemma \ref{l2.1}.

In fact, $E_k(f)$ is a discrete analog of an approximation of the
identity. The following Proposition \ref{pa.1} makes this precise,
whose proof is similar to that of \cite[Theorem 2.10]{d01} and we
omit the details.

\begin{prop}\label{pa.1}
(i) Let $p\in(1,\,\fz]$. The dyadic maximal function $M_d$ is bounded from
$L^1(\rn,\,\cb)$ to $L^{1,\,\fz}(\rn)$ and bounded from
$L^p(\rn,\,\cb)$ to $L^p(\rn)$.

(ii) If $f\in L^1_{\loc}(\rn,\,\cb)$, then $\lim_{k\to\fz}
E_k(f)(x)=\|f(x)\|_\cb $ and $\|f(x)\|_\cb\le M_d(f)(x)$ almost
everywhere.
\end{prop}

The following proposition provides the Calder\'on-Zygmund decomposition in our setting with a non-typical assumption on $f$ instead of the usual $f\in L^1$. This adds an extra layer of difficulty to the standard arguments as in \cite[Theorem 2.11]{d01}.

\begin{prop}\label{pa.2}
 Given a $\cb$-measurable
function  $f\in L^p_w(\rn,\,\cb)$ for certain $p\in[1,\,\fz)$ and
$w\in\ca_p(A)$, and a positive number $\lz$, then exists a sequence
$\{Q_j\}_j\subset \cq$ of disjoint dyadic cubes such that

\noindent (i) $\cup_jQ_j=\{x\in\rn: M_d(f)(x)>\lz\}$;

\noindent (ii) $\|f(x)\|_\cb<\lz$ for almost every $x\not\in \cup_j
Q_j$;

\noindent (iii) $\lz<\frac{1}{|Q_j|}\int_{Q_j}\|f(x)\|_\cb\, dx\le
C\lz$, where $C\ge1$ is a constant
 independent of $f$ and $\lz$;

\noindent (iv) for any $Q\in\{ Q_j\}_j$, there exists unique $\wz
Q\in \cq$ such that $Q\subset\wz Q$, $\ell(\wz Q)=\ell(Q)-1$ and
$\frac{1}{|\wz Q|}\int_{\wz Q} \|f(x)\|_\cb\, dx<\lz$.

\end{prop}

\begin{proof}

Let $p\in [1,\,\fz)$,  $w\in\ca_p(A)$ and $f\in L^p_w(\rn)$. It is
easy to see that $f\in L^1_\loc(\rn,\,\cb)$. In fact, if $p>1$, by
$w\in\ca_p(A)$, we have $w^{-p'/p}=w^{1-p'}\in \ca_{p'}(A)$, which
implies that $w\in L_{\loc}^{-p'/p}(\rn)$, where $p'\in\rn$
satisfying $1/p'+1/p=1$. Then for any $k\in\zz$ and $B_k$, by
H\"oler's inequality, we have
$$ \int_{B_k}\|f(x)\|_\cb\, dx\le \|f\|_{L^p_w(\rn,\cb)}
\lf\{\int_{B_k} [w(x)]^{-p/p'}\, dx\r\}^{1/p'}<\fz.$$ If $p=1$,
observing that
$\sup_{B}\frac1{|B|}\int_Bw(x)\,dx\sup_B[w(x)]^{-1}\ls1$, we have
$$ \int_{B_k}\|f(x)\|_\cb\, dx\le
 \sup_{B_k}[w(x)]^{-1}\int_{B_k} \|f(x)\|_\cb w(x)\, dx<\fz.$$

Moreover, we claim that for almost all $y\in\rn$, we have
$E_k(f)(y)\to 0, \ {\rm as} \ k\to\ -\fz.$ To see this, notice that
for almost all $y\in\rn$, by Lemma 2.1 (i), there exists an unique dyadic cube
$Q_{k,\,y}\in\cq_k$ for each $k\in \zz$ such that $y\in Q_{k,\,y}$.

Then for sufficient small $k\in\zz$, by $Q_{0,\,y}\subset Q_{k\,y}$,
Proposition 2.1 (i), Lemma 2.1 (iii) and (iv), we have
\begin{eqnarray*}
\frac{w(Q_{k,\,y})}{w(Q_{0,\,y})}\gtrsim
\frac{w(x_{Q_{k,\,y}}+B_{vk-u})}{w(x_{Q_{0,\,y}}+B_{u})}\gtrsim
\frac{|B_{vk-u}|^{1/p}}{|B_{u}|^{1/p}}\gtrsim b^{vk/p}.
\end{eqnarray*}
From this, H\"older's inequality, $w\in\ca_p(A)$ and $v<0$,
it follows that
\begin{eqnarray*}
&E_k(f)(y)&\le \|f\|_{L^p_w(\rn)}\frac1{|Q_{k,\,y}|}
\lf\{\int_{Q_{k,\,y}}[w(x)]^{-p'/p}\,dx\r\}^{1/p'}
\\
&&\ls \|f\|_{L^p_w(\rn)}[w(Q_{k,\,y})]^{-1/p}\to 0 \hs\hs {\rm as} \
k\to -\fz.
\end{eqnarray*}
Thus, the claim holds.

For each $k\in\zz$, set $$\Omega_k\equiv \{x\in\rn:\
E_k(f)(x)>\lz,\ {\rm and} \ \forall\ j<k,
 E_j(f)(x)\le \lz \}.$$
Then we have
$$
\{x\in\rn: M_d(f)(x)>\lz\}=\bigcup_k \Omega_k.
$$
Indeed, obviously, we have
$$\bigcup_k \Omega_k\subset \{x\in\rn: M_d(f)(x)>\lz\}.$$
On the other hand,  for almost all $y\in\rn$ such that
$M_d(f)(y)>\lz$, since
$E_k(f)(y)\to 0$ as $k\to -\fz,$
there exists a minimal $k_0\in\zz$ such that $E_{k_0}(f)(y)>\lz$ and
any $j< k_0$ $E_j(f)(y)\le \lz$. Thus, we obtain $y\in
\Omega_{k_0}$.

Moreover, observe that $\Omega_k$ can be covered by disjoint dyadic
cubes for each $k\in\zz$. In fact, if $Q\cap\Omega_k\ne\emptyset$,
then $Q\subset\Omega_k$ by the definition of $E_k f$. Also notice that
$\{\Omega_k\}_k$ are disjoint with each other. By this and
$\{x\in\rn: M_d(f)(x)>\lz\}=\cup_k \Omega_k$, we get (i).

By the definition of $\Omega_k$ and $\cup_k \Omega_k=\cup_j Q_j$, we
obtain that (ii), (iv) and the first inequality of (iii) holds.
Furthermore, for any $Q\in\{Q_j\}_j$, by (iv) and Lemma 2.1(iv),
there exists an unique dyadic cube $\wz Q\supset Q$ such that
$\ell(\wz Q)=\ell(Q)-1$ and
$$\frac1{|Q|}\int_{Q}\|f(x)\|_{\cb}\, dx\le \frac{|\wz Q|}{|Q|}\frac 1{|\wz Q|}
\int_{\wz Q}\|f(x)\|_\cb\, dx\le C \lz, $$
where $C\ge 1$ is a constant independent of $f$ and $\lz$. Thus, the second inequality of
(iii) holds. This finishes the proof of Proposition \ref{pa.2}.
\end{proof}

For any $f\in L^1_{\loc}(\rn,\, \cb)$ and $E\subset\rn$, set
$f_{E}\equiv \frac1{|E|}\int_E  f(x)\, dx,$ and define the sharp
maximal function associated with dilation $A$ by setting, for all
$x\in\rn$,
$$\cm^\sharp f(x)\equiv \sup_{k\in\zz,\, y\in \rn}\sup_{x\in y+B_k}
b^{-k}\int_{y+B_k} \|f(z)-f_{y+B_k}\|_\cb \, dz.$$

Then by a similar argument to that used in \cite[Proposition
6.4]{d01}, we have the following result. We omit the details.

\begin{prop}\label{pa.3}
For any $f\in L^1_\loc(\rn,\,\cb)$ and all $x\in\rn$, $\cm^\sharp
f(x)\le M(\|f\|_\cb)(x)$, and
$$\frac12 \cm^\sharp f(x)
\le
\sup_{k\in\zz,\,y\in\rn}\sup_{x\in y+B_k}\inf_{a\in
\cb}b^{-k}\int_{y+B_k} \|f(z)-a\|_\cb\, dz \le
\cm^\sharp f(x).
$$
\end{prop}

Based on this, we have the following conclusion.

\begin{lem}\label{la.2}
If $p_0,\,p\in [1,\,\fz)$, $p_0\le p$, $w\in \ca_p(A)$ and $f\in
L_{\loc}^1(\rn,\, \cb)$ such that $M_d(f)\in L^{p_0}_w(\rn)$, then
there exists a positive constant $C$, independent of $f$, such that
$$\int_\rn [M_d (f)(x)]^pw(x)\,dx\le C \int_\rn [\cm^\sharp(f)(x)]^p\,
w(x)\,dx.$$
\end{lem}

The proof of Lemma \ref{la.2} needs the following generalized
``good-$\lz$'' inequality,  which is a extension of
\cite[Lemma 7.10]{d01}.

\begin{lem}\label{la.3}
Let $p_0\in [1,\, \fz)$ and $w\in \ca_{p_0}(A)$. Then there exists a
positive constant $C_0$ such that for all $f\in L^{p_0}_w(\rn,\,
\cb)$, $\gz>0$ and $\lz>0$, $$w(\{x\in\rn:\ M_d(f)(x)>2\lz,\,
\cm^\sharp (f)(x)\le \gz\lz\})\le C_0\gz^{1/p}w(\{x\in\rn:\
 M_d(f)(x)>\lz\}).$$
\end{lem}
\begin{proof}
Fix $\lz,\, \gz>0$. Since $f\in L^{p_0}_w(\rn,\,
\cb)$, by  Proposition \ref{pa.2} the set
$\{x\in\rn: M_d(f)(x)>\lz\}$ can be written as the union of disjoint
dilated cubes. To show Lemma \ref{la.3}, it suffices to prove that
if $Q$ is one of such cubes, then $w(E)\ls\gz^{1/p}w(Q)$, where
$E\equiv \{x\in Q: M_d(f)(x)>2\lz,\, \cm^\sharp(f)(x)<\gz\lz\}$. By
Lemma \ref{l2.1} and Proposition \ref{p2.1}(i), we have
$$\frac{w(E)}{w(Q)}\le \frac{w(E)}{w(x_Q+B_{v\ell(Q)-u})}\ls \frac{|E|^{1/p}}
{|x_Q+B_{v\ell(Q)-u}|^{1/p}}\ls \frac{|E|^{1/p}}{|Q|^{1/p}},$$ where
$u$ and $v$ are the same as in Lemma \ref{l2.1}(iv). Therefore, to
finish the proof of Lemma \ref{la.3}, we only need to prove $|E|\ls
\gz |Q|$. By Proposition \ref{pa.2}(iv), there exists $\wz
Q\in\cq$ such that $\ell(\wz Q)=\ell(Q)-1$, $Q\subset \wz Q$ and
\begin{eqnarray}\label{ea.3}
\frac{1}{|\wz Q|}\int_{\wz Q} \|f(x)\|_\cb \, dx<\lz.
\end{eqnarray}
Furthermore, if $x\in Q$ and $M_d(f)(x)>2\lz$, then there exist
certain $k_0\in\zz$ and $Q_{k_0}\in\cq_{k_0}$ such that
$E_{k_0}(f)(x)>2\lz$, namely, $\frac{1}{|Q_{k_0}|}\int_{Q_{k_0}}
\|f(y)\|_\cb\, dx>2\lz,$ Proposition \ref{pa.2}(iv) further implies
that $Q_{k_0}\subset Q$. Therefore, for such $x$, we have
$$M_d(f\chi_Q)(x)\ge
E_{k_0}(f\chi_{Q})(x)=\frac{1}{|Q_{k_0}|}\int_{Q_{k_0}}
\|f\chi_Q(y)\|_\cb\, dy>2\lz,$$  from which and \eqref{ea.3}, it follows that
\begin{eqnarray*}
M_d((f-f_{\wz Q})\chi_Q)(x)
&&\ge
M_d(f\chi_Q)(x)-M_d(f_{\wz Q}\chi_Q)(x)
\\
&&\ge M_d(f\chi_Q)(x)-\frac{1}{|\wz Q|}\int_{\wz Q}
\|f(y)\|_\cb\,
dy>\lz,
\end{eqnarray*}
where we used the fact that
\begin{equation}\label{ea.4}
\lf\|\int_\Omega f(x)\, dx\r\|_\cb\le\int_\Omega\|f(x)\|_\cb\, dx
\end{equation}
for all measurable sets $\Omega$ and integrable functions $f$ on
$\Omega$; see \cite{g1,y95}. Therefore, $E\subset\{x\in\rn:
M_d((f-f_{\wz Q})\chi_Q)(x)>\lz\}$.

Moreover, by $\ell(\wz Q)=\ell(Q)-1$, Proposition
\ref{pa.1}(i) and Lemma 2.1, we have
\begin{eqnarray}\label{ea.5}
&& |\{x\in\rn:\ M_d((f-f_{\wz Q})\chi_Q)(x)>\lz\}|\\
&&\hs\ls
\frac{1}{\lz}\int_Q \|f(x)-f_{\wz Q}\|_\cb \, dx\nonumber\\
&&\hs\ls \frac{1}{\lz}\int_Q \|f(x)-f_{x_{\wz Q}+B_{v\ell(\wz
Q)+u}}\|_\cb\, dx+\frac{1}{\lz}\|f_{\wz Q}-f_{x_{\wz Q}+B_{v\ell(\wz
Q)+u}}\|_\cb\nonumber\\
&&\hs\ls \frac{|\wz Q|}{\lz}\frac{1}{b^{v\ell(\wz Q)+u}}\int_{x_{\wz
Q}+B_{v\ell(\wz Q)+u}} \|f(x)-f_{x_{\wz Q}+B_{v\ell(\wz
Q)+u}}\|_\cb\, dx \ls\frac{|Q|}{\lz}\inf_{x\in Q} \cm^\sharp
(f)(x).\nonumber
\end{eqnarray}
If the set $E$ is empty, there is nothing to prove. Otherwise, there
exists certain $x\in Q$ such that $\cm^\sharp(f)(x)<\gz\lz$, which
together with \eqref{ea.5} further implies that $|E|\ls \gz |Q|$.
This finishes the proof of Lemma \ref{la.3}.
\end{proof}

\begin{proof}[Proof of Lemma \ref{la.2}]
For $N>0$, let $${\mathrm {I_N}}\equiv\int_0^N p\lz^{p-1}
w(\{x\in\rn: M_d(f)(x)>\lz\})\, d\lz.$$ The assumptions that $p_0\le
p$ and $M_d(f)\in L^{p_0}_w(\rn)$ imply that $\mathrm{I_N}<\fz$.
Then, by Lemma \ref{la.3},
\begin{eqnarray*}
&{\mathrm {I_N}}&=2^p\int_0^{N/2} p \lz^{p-1} w(\{x\in\rn:
M_d(f)(x)>2\lz\})\,
d\lz\\
&& \le 2^p\int_0^{N/2} p\lz^{p-1}\bigg[w(\{x\in\rn:
M_d(f)(x)>2\lz,\,
\cm^\sharp(f)(x)\le \gz\lz\})\\
&&\hs +w(\{x\in\rn:
\cm^\sharp(f)(x)>\gz\lz\})\bigg]\, d\lz\\
&&\le C_0 2^p\gz^{1/p} {\mathrm
{I_N}}+\frac{2^p}{\gz^p}\int_0^{N\gz/2} p\lz^{p-1}w(\{x\in\rn:
\cm^\sharp(f)(x)>\lz\})\, d\lz.
\end{eqnarray*}
Choose $\gz$ such that $C_02^p\gz^{1/p}=1/2$. Thus, we obtain
$${\mathrm {I_N}}\le \frac{2^{p+1}}{\gz^p}\int_0^{ N\gz/2} p\lz^{p-1}w(\{x\in\rn:
\cm^\sharp(f)(x)>\lz\})\, d\lz,$$ which implies the desired
conclusion of the lemma. This finishes the proof of Lemma
\ref{la.2}.
\end{proof}

\begin{lem}\label{la.4}
If $\ct$ is a Calder\'on-Zygmund operator as in Definition
\ref{da.1}, then for each $s\in(1,\,\fz)$, there exists a positive
constant $C_s$ such that for all $f\in L^\fz_c(\rn,\,\cb_1)$ and
$x\in\rn$,
$$\cm^\sharp (\ct(f))(x)\le C_s [M(\|f\|_\cb^s)(x)]^{1/s}, $$
where $M$ is the Hardy-Littlewood maximal operator.
\end{lem}

\begin{proof}
Fix $s\in(1,\,\fz)$. For any given $x\in\rn$, pick $y\in\rn$ and
$k\in\zz$ such that $x\in y+B_k$. By Proposition \ref{pa.3}, to
complete the proof of Lemma \ref{la.4},  it suffices to find an
element $a\in\cb_2$ such that
$$b^{-k}\int_{y+B_k} \|\ct(f)(z)-a\|_{\cb_2}\, dz\ls
[M(\|f\|_{\cb_1}^s)(x)]^{1/s}.$$ Decompose $f$ as $f=f_1+f_2$, where
$f_1=f\chi_{y+B_{k+2\sz}}$. Now let $a\equiv\ct(f_2)(x)$. By
Definition \ref{da.1} and $f\in L^\fz_c(\rn,\,\cb_1)$, we have that
$a\in\cb_2$ and
\begin{eqnarray*}
&&b^{-k}\int_{y+B_k} \|\ct(f)(z)-a\|_{\cb_2}\, dz\\
&&\hs \le b^{-k}\int_{y+B_k}
\|\ct(f_1)(z)\|_{\cb_2}\,dz+b^{-k}\int_{y+B_k}
\|\ct(f_2)(z)-a\|_{\cb_2}\,dz \equiv {\mathrm I}+{\mathrm
{II}}.
\end{eqnarray*}

By H\"older's inequality and the boundedness of $\ct$ from
$L^s(\rn,\, \cb_1)$ to $L^s(\rn,\, \cb_2)$ (see Lemma \ref{la.1}), we
then have
\begin{eqnarray*}
&&{\mathrm I}\ls \lf\{b^{-k}\int_{y+B_k} \|\ct(f_1)(z)\|_{\cb_2}^s\,
dz
\r\}^{1/s}\\
&&\hs \ls \lf\{b^{-k-2\sz}\int_{y+B_{k+2\sz}} \|f(z)\|_{\cb_1}^s\,
dz \r\}^{1/s}\ls [M(\|f\|_{\cb_1}^s)(x)]^{1/s}.
\end{eqnarray*}

Moreover, if $x-y\in B_k$, $y-z\in B_k$  and $\az-y\in
B_{k+2\sz}^\complement$, by \eqref{e2.1} and \eqref{e2.2}, we obtain
$\rho(z-x)\le b^{-\sz}\rho(x-\az)$ and $\rho(x-\az)\ge b^{k+\sz}$.
From this, \eqref{ea.2} and H\"older's inequality, it
follows that
\begin{eqnarray*}
&{\mathrm {II} }&\ls
b^{-k}\int_{y+B_k}\int_{\rn\setminus(y+B_{k+2\sz})} \|\ck(z,
\, \az)-\ck(x,\,\az)\|_{L(\cb_1,\,\cb_2)} \|f(\az)\|_{\cb_1}\, d\az\,dz\\
&&\ls
b^{-k}\int_{y+B_k}\int_{\rho(x-\az)\ge b^{k+\sz}}\frac{[\rho(z-x)]^\ez}
{[\rho(x-\az)]^{1+\ez}} \|f(\az)\|_{\cb_1}\, d\az\, dz\\
&&\ls b^{-k}\int_{y+B_k}
b^{k\ez}\sum_{j=0}^\fz\int_{b^{k+\sz+j+1}\le \rho(x-\az)<
b^{k+2\sz+j}} b^{-(k+j)(1+\ez)}\|f(\az)\|_{\cb_1}\,d\az\, dz\\
&&\ls \sum_{j=0}^\fz
b^{-j\ez}b^{k+2\sz+j+1}\int_{y+B_{k+2\sz+j+1}}\|f(\az)\|_{\cb_1}\,
d\az\ls [M(\|f\|_{\cb_1}^s)(x)]^{1/s}.
\end{eqnarray*}
Combining the estimates of I and II yields the desired result and
thus finishes the proof Lemma \ref{la.4}.
\end{proof}

\begin{proof}[Proof of Theorem \ref{ta.1}]
We first prove that $\ct$ is bounded from $L_w^p(\rn,\,\cb_1)$ to
$L_w^p(\rn,\,\cb_2)$ when $p\in(1,\, \fz)$ and $w\in \ca_p(A)$. By
\cite[Lemma 8,\, p. 5]{st89}, there exists $r\in(1,\, p)$ such that
$w\in \ca_{p/r}(A)$.
Since $L^\fz_c(\rn,\, \cb_1)$ is dense in $L^p_w(\rn,\, \cb_1)$ (see
\cite[Remark 2.2]{gly1}), then we only need to prove the conclusions of
Theorem \ref{ta.1} by
assuming that $f\in L^\fz_c(\rn,\, \cb_1)$. Observe that if
$\ct(f)\in L^p_w(\rn,\,\cb_2)$, then by Proposition \ref{pa.1}(ii),
Lemma \ref{la.2}, Lemma \ref{la.4} and Proposition \ref{p2.1}(ii),
we have
\begin{eqnarray*}
&\dint_\rn \|\ct(f)(x)\|_{\cb_2}^pw(x)\,dx&\le \int_{\rn}
[M_d(\ct(f))(x)]^pw(x)\,dx\ls \int_{\rn} [\cm^\sharp(\ct(f))(x)]^pw(x)\,dx\\
&&\ls \int_\rn [M(\|f\|_{\cb_1}^r)(x)]^{p/r}w(x)\,dx
\ls \int_\rn \|f(x)\|_{\cb_1}^pw(x)\,dx.
\end{eqnarray*}

Now we turn to prove $\ct(f)\in L^p_w(\rn,\,\cb_2)$. Since $f\in
L^\fz_c(\rn,\, \cb_1)$, we assume that $\supp f\subset B_{k_0}$ for
certain $k_0\in\zz$. Write
\begin{eqnarray*}
\|\ct(f)\|_{L^p_w(\rn,\,\cb_2)}^p=
\lf\{\int_{B_{k_0+\sz}}+\int_{B_{k_0+\sz}^\complement}\r\}
\|\ct(f)(x)\|_{\cb_2}^pw(x)\,dx={\mathrm I}+{\mathrm {II}}.
\end{eqnarray*}

By \cite[p. 7]{st89}, there exists $\eta\in(1,\,\fz)$ such that $w$
satisfies the reverse H\"older's inequality, which implies that
$w\in L^{\eta}_{\loc}(\rn)$. This combined with H\"older's inequality
and Lemma \ref{ta.1} yields that $\mathrm{I}<\fz$.

For $x\in (B_{k_0+\sz})^\complement$ and $y\in B_{k_0}$, we have
$x-y\in B_{k_0}^\complement$ and $\rho(x)\ls \rho(x-y)+\rho(y)\ls
\rho(x-y)$. By this, $f\in L^\fz_c(\rn,\, \cb_1)$, \eqref{ea.4} and
\eqref{ea.1}, we have
\begin{eqnarray*}
\|\ct(f)(x)\|_{\cb_2} &&\le \int_\rn \|f(y)\|_{\cb_1}\|\ck(x,\,
y)\|_{L(\cb_1,\,\cb_2)}\, dy
\ls\int_{B_{k_0}}\frac{1}{\rho(x-y)}\,dy\ls \rho(x)^{-1}.
\end{eqnarray*}
Therefore,
\begin{eqnarray*}
{\mathrm {II}}\ls \sum_{j=k_0}^\fz \int_{B_{\sz+j+1}\setminus
B_{\sz+j}} \rho(x)^{-p}w(x)\,dx\ls \sum_{j=k_0}^\fz
b^{-jp}w(B_{\sz+j+1}).
\end{eqnarray*}
By $w\in\ca_{p/s}(A)$ and Proposition \ref{p2.1}(i), we have
$w(B_{\sz+j+1})\ls b^{jp/s}w(B_{k_0})$, which together with
$s\in(1,\,\fz)$ implies that ${\mathrm {II}}$ is finite. Thus,
$T(f)\in L^p_w(\rn,\,\cb_2)$, which completes the proof of the boundedness of
$T$ from $L^p_w(\rn,\,\cb_1)$ to $L^p_w(\cb_2)$.

Finally, we prove that $\ct$ is bounded from $L_w^1(\rn,\,\cb_1)$ to
$L_w^{1,\,\fz}(\rn,\,\cb_2)$. Fix $\lz>0$ and $f\in
L^\fz_c(\rn,\,\cb_1)$. By Proposition \ref{pa.2}, there exists a
sequence $\{Q_j\}_j$ of disjoint dilated cubes such that the conclusions
(i)-(iv) of Proposition \ref{pa.2} hold. Then we write $f=g+b$, where
$$g(x)\equiv\lf\{\begin{array}{ll}f(x)\hs\hs & x\in\rn\setminus\bigcup_j
Q_j\\
\dfrac{1}{|Q_j|}\int_{Q_j}f(y)\, dy\hs\hs & x\in Q_j\end{array}\r.$$
and $b(x)\equiv \sum_j b_j(x)$ with
$$b_j(x)=\lf\{f(x)-\frac{1}{|Q_j|}\int_{Q_j}f(y)\, dy\r\}\chi_{Q_j}(x).$$
Thus by Proposition \ref{pa.2} and \eqref{ea.4}, we obtain
\begin{eqnarray}\label{ea.6}
\|g(x)\|_{\cb_1}\le \lz \hs \mathrm{for\ almost\ all}\ x\in\rn,\, \supp b\subset
\bigcup_j Q_j\ \mathrm{and}\ \int_{Q_j} b(x)\,dx=0.
\end{eqnarray}
So the estimate of $w(\{x\in\rn: \|\ct(f)(x)\|_{\cb_2}>2\lz\})$ is
reduced to the estimates of $w(\{x\in\rn: \|\ct(g)(x)\|_{\cb_2}>\lz\})$ and
$w(\{x\in\rn: \|\ct(b)(x)\|_{\cb_2}>\lz\}).$ Notice that
$w\in\ca_1(A)$ implies $w\in\ca_2(A)$ and thus, $\ct$ is bounded
from $L^2_w(\rn,\, \cb_1)$ to $L^2_w(\rn,\,\cb_2)$ as already proved
above in this proof. Then by \eqref{ea.6}, we have
\begin{eqnarray*}
&w(\{x\in\rn: \|\ct(g)(x)\|_{\cb_2}>\lz\})&\le
\frac{1}{\lz^2}\int_\rn
\|\ct(g)(x)\|_{\cb_2}^2w(x)\,dx\\
&&\ls \frac{1}{\lz^2}\int_\rn\|g(x)\|_{\cb_1}^2w(x)\,dx
\ls \frac{1}{\lz} \int_\rn \|g(x)\|_{\cb_1}w(x)\,dx.
\end{eqnarray*}
To obtain a desired estimate for $\ct(g)$, we still need to show that
$$\int_\rn\|g(x)\|_{\cb_1}w(x)\,dx\ls \int_\rn\|f(x)\|_{\cb_1}\,
w(x)\,dx.$$ Notice that for all $x\in\rn\setminus\cup_j Q_j$, we
have $g(x)=f(x)$. On each $Q_j$, by \eqref{ea.4} and $w\in\ca_1(A)$,
we have
\begin{eqnarray*}
\dint_{Q_j} \|g(x)\|_{\cb_1} w(x)\,dx\le
\int_{Q_j}\frac{1}{|Q_j|}\int_{Q_j} \|f(y)\|_{\cb_1}\, dy\, w(x)\,dx
\ls \int_{Q_j} \|f(y)\|_{\cb_1} w(y)\,dy.
\end{eqnarray*}
Since $\{Q_j\}_j$ are disjoint, we further have $$w(\{x\in\rn: \|\ct
(g)(x)\|_{\cb_2}>\lz\})\ls \int_{Q_j} \|f(y)\|_{\cb_1} w(y)\,dy,$$
which completes the estimate for $\ct(g)$.

On the other hand, set $Q_j^\ast\equiv x_{Q_j}+B_{v\ell(Q_j)+u+2\sz}$,
where $u,\, x_Q,\, v$ and $\ell(Q_j)$ are as in Lemma
\ref{l2.1}. Then we obtain
$$w(\{x\in\rn: \|\ct(b)(x)\|_{\cb_2}>\lz\})\le w\lf(\cup_j Q_j^\ast\r)
+w(\{x\in\rn\setminus\cup_j Q_j^\ast: \|\ct(b)(x)\|_{\cb_2}>\lz\}).$$
Since $w\in\ca_1(A)$, by Proposition \ref{p2.1}, Lemma \ref{l2.1},
Proposition \ref{pa.3}(iv) and the definition of $\ca_1(A)$, we
have
\begin{eqnarray*}
w\lf(\bigcup_j Q_j^\ast\r)&&\le \sum_j w(Q_j^\ast)\ls \sum_j
|Q_j|\frac{w(Q^\ast_j)}{|Q^\ast_j|}
\ls \frac{1}{\lz}\sum_j\int_{Q_j} \|f(y)\|_{\cb_1}\, w(y)\,dy\\
&&\ls \frac{1}{\lz}\|f\|_{L^1_w(\rn,\, \cb_1)}.
\end{eqnarray*}
Moreover, from the fact that $b_j$ has zero average on $Q_j$, and
\eqref{ea.4}, it follows that
\begin{eqnarray}\label{ea.7}
&&w\lf(\{x\in\rn\setminus \cup_j Q^\ast_j: \|
\ct(b)(x)\|>\lz\}\r)\\
&&\hs\le \frac{1}{\lz}\sum_j \int_{\rn\setminus \cup_j Q^\ast_j}
\|\ct(b_j)(x)\|_{\cb_2}\,
w(x)\,dx\nonumber\\
&&\hs= \frac{1}{\lz}\sum_j\int_{\rn\setminus \cup_j Q^\ast_j} \lf|
\int_{Q_j}[\ck(x,\,y)-\ck(x,\, x_{Q_j})]b_j(y)\, dy\r|\,w(x)\,dx\nonumber\\
&&\hs\le \frac{1}{\lz}\sum_j \int_{Q_j}\int_{\rn\setminus \cup_j
Q^\ast_j}\|\ck(x,\,y)-\ck(x,\,
x_{Q_j})\|_{L(\cb_1,\,\cb_2)}\,w(x)\,dx\|b_j(y)\|_{\cb_1}\, dy.\nonumber
\end{eqnarray}
Observe that $x-x_{Q_j}\in B^\complement_{v\ell(Q_j)+u+2\sz}$ and
$y-x_{Q_j}\in B_{v\ell(Q_j)+u}$ imply that $\rho(y-x_{Q_j})\le
b^{-3\sz}\rho(x-x_{Q_j})$. Then by \eqref{ea.2}, for all $y\in Q_j$, we have
\begin{eqnarray*}
&&\int_{\rn\setminus \cup_j Q^\ast_j}\|\ck(x,\,y)-\ck(x,\,
x_{Q_j})\|_{L(\cb_1,\,\cb_2)}\,w(x)\,dx\\
&&\hs \ls\int_{\rn\setminus \cup_j
Q^\ast_j}\frac{\rho(y-x_{Q_j})^\ez}
{\rho(x-x_{Q_j})^{1+\ez}}\,w(x)\,dx\\
&&\hs \ls\sum_{k=0}^\fz b^{-k\ez}\frac{1}{b^{v\ell(Q_j)+u+2\sz+k+1}}
\int_{B_{v\ell(Q_j)+u+2\sz+k+1}} w(x)\, dx\ls M(w)(y).
\end{eqnarray*}
From this, $w\in\ca_1(A)$,
\begin{eqnarray*}
\int_{Q_j}\|b_j(y)\|_{\cb_1}\, w(y)\,dy=\int_{Q_j}\|b(y)\|_{\cb_1}\,w(y)
dy\le \int_{Q_j}(\|f(y)\|_{\cb_1}+\|g(y)\|_{\cb_1})\, w(y)\,dy
\end{eqnarray*}
and \eqref{ea.7}, it follows that
\begin{eqnarray*}
&& w(\{x\in\rn\setminus \cup_j Q^\ast_j: \|
\ct(b)(x)\|_{\cb_2}>\lz\})\\
&&\hs\ls \frac1{\lz}\sum_j \int_{Q_j}\|b_j(y)\|_{\cb_1}\, M(w)(y)\,dy\\
&&\hs \ls \frac1{\lz}\sum_j \int_{Q_j}(\|f(y)\|_{\cb_1}+\|g(y)\|_{\cb_1})\,
w(y)\,dy\\
&&\hs\ls \frac1{\lz}\sum_j \int_{Q_j}\|f(y)\|_{\cb_1}\, w(y)\,dy
\ls \frac1{\lz}\int_\rn \|f(y)\|_{\cb_1}\, w(y)\,dy.
\end{eqnarray*}
This finishes the proof of Theorem \ref{ta.1}.
\end{proof}

\bigskip

\noindent Marcin Bownik

 \medskip

\noindent Department of Mathematics, University of Oregon, Eugene,
OR 97403-1222, USA

\smallskip

\noindent {\it E-mail address}: \texttt{mbownik@uoregon.edu}

\bigskip

\noindent Baode Li, Dachun Yang (Corresponding author) and Yuan Zhou

 \medskip

\noindent School of Mathematical Sciences, Beijing Normal
University, Laboratory of Mathematics and Complex Systems, Ministry
of Education, Beijing 100875, People's Republic of China

\smallskip

\noindent{\it E-mail addresses}: \texttt{baodeli@mail.bnu.edu.cn,
dcyang@bnu.edu.cn} and

\hspace{2.45cm} \texttt{yuanzhou@mail.bnu.edu.cn}

\end{document}